\documentclass[10pt]{amsart}
\usepackage{amsmath,amstext,amssymb,amsopn,amsthm,mathrsfs}
\newtheorem{theo}{Theorem}[section]
\newtheorem{prop}{Proposition}[section]

\numberwithin{equation}{section}
\newcommand{\dem}{\medskip \par \noindent \mbox{\bf Proof }}

\begin{document}

\title[Transference result Jacobi Littlewood-Paley $g-$function] {A transference result of the $L^p$ continuity of the Jacobi Littlewood-Paley $g-$function
to the Gaussian and Laguerre Littlewood-Paley $g-$function.}

\author{Eduard Navas}
\address{Departamento de Matem\'aticas, Universidad Nacional Experimental Francisco de Miranda, Punto Fijo, Venezuela.}
\email{[Eduard Navas]enavas@correo.unefm.edu.ve}
\author{Wilfredo O. Urbina}
\address{
 Department of Mathematics and Actuarial Sciences, Roosevelt University  Chicago, Il, 60605, USA.}
\email{[Wilfredo Urbina]wurbinaromero@roosevelt.edu}
\thanks{\emph{2000 Mathematics Subject Classification} Primary 42C10; Secondary 26A24}
\thanks{\emph{Key words and phrases:} Transference, Littlewood-Paley $g$  function, Orthogonal polynomials.}

\begin{abstract}
In this paper we develop a transference method to obtain the
$L^{p}$-continuity of the Gaussian-Littlewood-Paley $g$ function
and the $L^{p}$-continuity of the Laguerre-Littlewood-Paley $g$
function from the $L^{p}$-continuity of the
Jacobi-Littlewood-Paley $g$ function, in dimension one, using the
well known asymptotic relations between Jacobi polynomials and
Hermite and Laguerre polynomials.
\end{abstract}

 \maketitle
\section{Preliminaries}
It is well knownn, in the theory of classical orthogonal polynomials, the asymptotic relations between Jacobi polynomials and
Hermite and Laguerre polynomials. Using those asymptotic relations
we develop a transference method to obtain $L^{p}$-continuity for
the Gaussian-Littlewood-Paley $g$ function and the $L$-continuity
of the Laguerre-Littlewood-Paley $g$ function from the $L^{p}$-
continuity of the Jacobi-Littlewood-Paley $g$ function, in
dimension one. We are going to
use the normalizations given in G. Szeg\"o's book \cite{sz}, for all classical polynomials.

\begin{itemize}
\item {\em Jacobi polynomials:} For  $\alpha,\beta > -1,$ the
Jacobi polynomials $\{P^{(\alpha,\beta)}_n\}_{n \in \mathbb N }$
are defined as the orthogonal  polynomials associated with the
Jacobi measure $\mu_{\alpha,\beta}$ (or beta measure) in $(-1,1)$,
defined as

 \begin{eqnarray}
\mu_{\alpha,\beta}(dx) &=& \omega_{\alpha,\beta}(x) dx = \nonumber
\eta_{\alpha,\beta} \chi_{(-1,1)}(x)(1-x)^{\alpha}(1+x)^{\beta}
dx,
\end{eqnarray}
where $ \eta_{\alpha,\beta} = \frac{1}{2^{\alpha+\beta+1}
B(\alpha+1,\beta+1)} =
\frac{\Gamma(\alpha+\beta+2)}{2^{\alpha+\beta+1}
\Gamma(\alpha+1)\Gamma(\beta+1)}$.

The function  $\omega_{\alpha,\beta}$ is called the (normalized)
{\em Jacobi weight}.\\

The Jacobi polynomials can be obtained from the polynomial
canonical basis $\{ 1, x, x^2,\cdots,x^n,\cdots \} $ using the
Gram-Schmidt orthogonalization process with respect to the inner
product in $L^2(\mu_{\alpha,\beta})$. Thus we have the {\em
orthogonality property} of Jacobi polynomials  with respect to
$\mu_{\alpha,\beta} $,
\begin{equation} \label{proJacOrt}
\int^{\infty}_{-\infty} P^{(\alpha,\beta)}_n(y)
P^{(\alpha,\beta)}_m(y) \, \mu_{\alpha,\beta}(dy)
=\eta_{\alpha,\beta}  h_{n}^{\left( \alpha ,\beta
\right)}\delta_{n,m} = \hat{h_{n}}^{\left( \alpha ,\beta
\right)}\delta_{n,m},
\end{equation}
$n,m =0,1,2, \cdots$, where
\begin{equation}\label{l2norm}
h_{n}^{\left( \alpha ,\beta \right)}=\frac{ 2^{\alpha +\beta
+1}}{(2n+\alpha +\beta +1)}\frac{\Gamma \left( n+\alpha +1\right)
\Gamma \left( n+\beta +1\right) }{\Gamma \left( n+1\right) \Gamma
\left( n+\alpha +\beta +1\right) },
\end{equation}
and
\begin{eqnarray*}
 \hat{h_{n}}^{\left( \alpha ,\beta \right)}&=&\frac{1}{(2n+\alpha
+\beta +1)}\frac{\Gamma(\alpha+\beta+2)\Gamma \left( n+\alpha
+1\right) \Gamma \left( n+\beta +1\right)
}{ \Gamma(\alpha+1)\Gamma(\beta+1)\Gamma \left( n+1\right) \Gamma \left( n+\alpha +\beta +1\right) }\\
&=&  \| P_n^{\left( \alpha
,\beta\right)}\|_{2,(\alpha,\beta)}^{2}.
\end{eqnarray*}

On the other hand, the Jacobi polynomial of parameter $(\alpha, \beta)$ of degree $n$, $P_n^{\left( \alpha
,\beta\right)},$ is a polynomial solution of the {\em Jacobi
differential equation}, with parameters $\alpha ,\beta, n$,
\begin{equation}\label{ecuadif}
\left( 1-x^{2}\right) y^{\prime \prime }+\left[ \beta -\alpha
-\left( \alpha +\beta +2\right) x\right] y^{\prime }+n\left(
n+\alpha +\beta +1\right) y=0,
\end{equation}
i.e.  $P_n^{\left( \alpha ,\beta\right)}$ is an eigenfunction of
the (one-dimensional) second order diffusion operator
\begin{equation}\label{jacobiop}
{\mathcal L}^{\alpha,\beta}= -(1-x^2)  \frac{d^2}{dx^2}
-(\beta-\alpha-(\alpha+\beta+2)x)\frac{d}{dx},
\end{equation}
associated with the eigenvalue $ \lambda^{\alpha+\beta}_n=
n(n+\alpha+\beta+1)$. ${\mathcal L}^{\alpha,\beta}$ is  called the
{\em Jacobi differential operator.} Observe that if we choose $
\delta_{\alpha,\beta} = \sqrt{1-x^2} \frac{d}{dx}, $ and consider
its formal $L^2(\mu_{\alpha,\beta})$-adjoint,
$$
\delta^*_{\alpha,\beta} = -\sqrt{1-x^2} \frac{d}{dx} +
[(\alpha+\frac{1}{2}) \sqrt{\frac{1+x}{1-x}} - (\beta+\frac{1}{2})
\sqrt{\frac{1-x}{1+x}}] I,
$$
then $ {\mathcal L}^{\alpha,\beta} = \delta^*_{\alpha,\beta}
\delta_{\alpha,\beta}.$
The differential operator $\delta_{\alpha,\beta}$ is considered the ``natural" notion of derivative in the Jacobi case.\\

The operator semigroup associated to the Jacobi polynomials is
defined for positive or bounded measurable Borel functions of
 $(-1,1)$, as
\begin{equation}\label{Jac1}
T^{\alpha,\beta}_t f(x)=    \int_{-1}^{1}
p^{\alpha,\beta}(t,x,y)f(y) \mu_{\alpha,\beta}(dy),
\end{equation}
where
$$p^{\alpha,\beta}(t,x,y) =   \sum_k \frac{e^{- k(k+\alpha+\beta+1) t}}{ \hat{h_{k}}^{\left( \alpha ,\beta \right)}} P^{(\alpha,\beta)}_k(x) P^{(\alpha,\beta)}_k(y).$$
There is not a simple explicit representation of $p^{\alpha,\beta}(t,x,y)$ since  the eigenvalues $\lambda_n$  are not
linearly distributed; there is one obtained by G. Gasper \cite{gasp} which is analog of Bailey's $F_4$ representation of the kernel of Abel summability for Jacobi series, also called the Jacobi-Poisson integral, see \cite{bai}. From that form, taking $x=-y=1$, it can be proved that $p^{\alpha,\beta}(t, x, y) $ is a positive kernel.\\

$\{ T^{\alpha,\beta}_t \}$ is called  the  {\em  Jacobi semigroup}
and can be proved that is a Markov semigroup, for details see
\cite{ur2}. The Jacobi-Poisson semigroup $\{ P^{\alpha,\beta}_t
\}$ can be defined, using Bochner's {subordination} formula,
\begin{eqnarray*}
e^{-\lambda^{1/2}t} = \frac{1}{\sqrt{\pi}} \int_0^\infty
\frac{e^{-u} }{\sqrt{u}} e^{-{\frac{\lambda t^2}{4u}}}du.
\end{eqnarray*}
as the subordinated semigroup of the Jacobi semigroup,
\begin{eqnarray*}
P_t^{\alpha,\beta} f(x) = \frac{1}{\sqrt \pi} \int_0^{\infty}
\frac{e^{-u}}{\sqrt u} T^{\alpha,\beta}_{t^2/4u}f(x) du.
\end{eqnarray*}

For a function $f\in
L^{2}\left(\left[-1,1\right],\mu_{(\alpha,\beta)}\right)$ let us
consider its Fourier- Jacobi expansion
\begin{equation}\label{jacobides}
f= \sum _{k=0}^\infty\frac{\langle f,P_{k}^{(\alpha,\beta)}
\rangle}{ \hat{h_{k}}^{\left( \alpha ,\beta \right)}}
P_{k}^{(\alpha,\beta)},
\end{equation}
where
$$\langle f,P_{k}^{(\alpha,\beta)} \rangle= \int_{-1}^{1}f(y) P_{k}^{(\alpha,\beta)}(y) \mu_{\alpha,\beta}(dy).$$
 Then, the action of $T_t$ and $P_t$ can be expressed as
\begin{equation*}
T_t^{\alpha,\beta} f =  \sum _{k=0}^\infty \frac{\langle
f,P_{k}^{(\alpha,\beta)} \rangle}{ \hat{h_{k}}^{\left( \alpha
,\beta \right)}} e^{-\lambda_k t} P_{k}^{(\alpha,\beta)},
\end{equation*}
and
\begin{equation*}
P_t^{\alpha,\beta} f = \sum _{k=0}^\infty \frac{\langle
f,P_{k}^{(\alpha,\beta)} \rangle}{ \hat{h_{k}}^{\left( \alpha
,\beta \right)}} e^{-\sqrt{\lambda_k} t} P_{k}^{(\alpha,\beta)}.
\end{equation*}
Following the classical case, the Jacobi-Littlewood-Paley $g$
function can be define as

\begin{equation}\label{fung-jacob}
g^{(\alpha,\beta)}f(x)=\int_{0}^{\infty}t|\nabla_{(\alpha,\beta)}
P^{(\alpha,\beta)}_{t}f(x)|^{2}dt
\end{equation}
where $\nabla_{(\alpha,\beta)}=\left(\frac{\partial}{\partial
t},\delta_{\alpha,\beta}\right)=\left(\frac{\partial}{\partial
t},\sqrt{1-x^{2}}\frac{\partial}{\partial x}\right)$

The $L^{p}$-continuity of the Jacobi-Littlewood-Paley $g$-function
$g(\alpha;\beta)$, was proved by A. Nowak and P. Sj\"ogren in
\cite{NowakSjogren} .
\begin{theo} \label{contfuctgJacob}
Assume that $1< p < \infty$ and $\alpha, \beta \in
[-1/2,\infty)^d$. There exists a constant $c_p$ such that
\begin{equation}
\| g^{(\alpha,\beta)} f \|_{p,(\alpha, \beta)} \leq c_p \| f
\|_{p,(\alpha, \beta)}.
\end{equation}
\end{theo}
\vspace{0.3cm}

\item {\em Hermite polynomials}:  The {\em Hermite polynomials}
$\{H_n\}_{n}$, are defined as the orthogonal polynomials
associated with the Gaussian measure in $\mathbb R$, $ \gamma (dx)
=  \frac{e^{-x^2}}{\sqrt{\pi}} dx ,$ i.e.
\begin{equation} \label{proHerOrt}
\int^{\infty}_{-\infty} H_n(y) H_m(y) \, \gamma (dy) =   2^n
n!\delta_{n,m},
\end{equation}
$n,m =0,1,2, \cdots$, with the {\em normalization}
\begin{equation*}
H_{2n+1}(0) = 0, \quad H_{2n}(0)= (-1)^n \frac{(2n)!}{n!}.
\end{equation*}
We have
\begin{eqnarray}
H_n^{'}(x) &=& 2n H_{n-1}(x),\\
 H_n^{''}(x) &-& 2xH_n^{'}(x) + 2n H_n(x) = 0, \label{eqHerm}
\end{eqnarray}
thus $H_n$ is an eigenfunction of the one dimensional {\em
Ornstein-Uhlenbeck operator} (or harmonic oscillator  operator),
\begin{equation}
L=  -\frac{1}{2}\frac{d^2}{dx^2} + x \frac{d}{dx},
\end{equation}
associated with the eigenvalue $\lambda_n= n$. Observe that if we
choose $ \delta_\gamma = \frac{1}{\sqrt{2}}\frac{d}{dx}, $ and
consider its formal $L^2(\gamma )$-adjoint,
$$
\delta^*_\gamma = - \frac{1}{\sqrt{2}}\frac{d}{dx} + \sqrt{2} x I
$$
then $ L= \delta^*_\gamma \delta_\gamma.$
The differential operator $\delta_\gamma$ is considered the ``natural" notion of derivative in the Hermite case.\\

The Gaussian-Littlewood-Paley $g$ function can be defined as

\begin{equation}\label{fung-Hermite}
g^{\gamma}f(x)=\int_{0}^{\infty}t|\nabla_{\gamma}
P^{\gamma}_{t}f(x)|^{2}dt
\end{equation}
where $\nabla_{\gamma}=\left(\frac{\partial}{\partial
t},\delta_{\gamma}\right)=\left(\frac{\partial}{\partial
t},\frac{1}{\sqrt{2}}\frac{\partial}{\partial x}\right)$\\

The $L^p$ continuity of the of the Gaussian-Littlewood-Paley $g$
function was proved by C. Guti{\'e}rrez in \cite{gu},\\

\begin{theo} \label{contfuctgHermit}
Assume that $1< p < \infty$. There exists a constant $c_p$ such
that
\begin{equation}
\| g^{\gamma} f \|_{p,\gamma} \leq c_p \| f \|_{p,\gamma}.\\
\end{equation}
\end{theo}
\vspace{0.3cm}

\item  {\em Laguerre polynomials}: For $\alpha>-1$, the {\em
Laguerre polynomials} $\{L^{\alpha} _k\}$  are defined as the
orthogonal  polynomials  associated with the Gamma measure on
$(0,\infty)$, $ \mu_{\alpha} (dx) = \chi_{(0,\infty)}(x)
\frac{x^{\alpha}e^{-x}}{\Gamma(\alpha+1)} dx ,$ i.e.
\begin{equation} \label{proLagrOrt}
\int^{\infty}_{0} L^{\alpha}_n(y) L^{\alpha}_m(y) \,\mu_{\alpha}
(dy) =\dbinom{n+\alpha}{n} \delta_{n,m} =
\frac{\Gamma(n+\alpha+1)}{\Gamma(\alpha+1) n!} \delta_{n,m},
\end{equation}
$n,m =0,1,2, \cdots$ We have
\begin{equation} \label{eq:aPri}
 (L_k^{\alpha}(x))' = -L_{k-1}^{\alpha+1}(x).
\end{equation}
\begin{equation}\label{eq:aecdif}
x(L_k^{\alpha}(x))''
+(\alpha+1-x)(L_k^{\alpha}(x))'+kL_k^{\alpha}(x)=0.
\end{equation}
thus $L_k^{\alpha}$ is an eigenfunction of the (one-dimensional)
{\em Laguerre differential operator}
\begin{equation}
{\mathscr L}^{\alpha}= -x \frac{d^2}{dx^2} -(\alpha +1- x)
\frac{d}{dx},
\end{equation}
associated with the eigenvalue $\lambda_k=k$. Observe that if we
choose $ \delta_{\alpha} = \sqrt{x} \frac{d}{dx}, $ and consider
its formal $L^2(\alpha )$-adjoint,
$$
\delta^*_{\alpha} = - \sqrt{x}\frac{d}{dx} +
[\frac{\alpha+1/2}{\sqrt{x}}+\sqrt{x}] I
$$
then $ {\mathscr L}^\alpha= \delta^*_{\alpha} \delta_{\alpha}.$
The differential operator $\delta_{\alpha}$ is considered the
``natural" notion of derivative in the Laguerre case.\\

The Laguerre-Littlewood-Paley $g$ function can be defined as
\begin{equation}\label{fung-Lagu}
g^{\alpha}f(x)=\int_{0}^{\infty}t|\nabla_{\alpha}
P^{\alpha}_{t}f(x)|^{2}dt
\end{equation}
where $\nabla_{\alpha}=\left(\frac{\partial}{\partial
t},\delta_{\alpha}\right)=\left(\frac{\partial}{\partial
t},\sqrt{x}\frac{\partial}{\partial x}\right)$.\\

The $L^{p}$ continuity of the Laguerre-Littlewood-Paley $g$
function was proved by A. Nowak in \cite{nowak}.\\

\begin{theo} \label{contfuctgLague}
Assume that $1< p < \infty$  and $\alpha\in [1/2,\infty)^{d}$.
There exists a constant $c_p$ such that
\begin{equation}
\| g^{\alpha} f \|_{p,\alpha} \leq c_p \| f \|_{p,\alpha}.
\end{equation}
\end{theo}

\item  Finally, let us consider the asymptotic relations between
Jacobi polynomials and other classical orthogonal polynomials (see
\cite{sz}, (5.3.4) and (5.6.3)),

\begin{enumerate}
\item[i)] For Hermite polynomials,
\begin{equation}\label{jacoherm}
\lim_{\lambda \rightarrow \infty} \lambda^{-n/2} C^{\lambda}_n
(x/\sqrt{\lambda})= \frac{H_n(x)}{n!},
\end{equation}
\item[ii)] For Laguerre polynomials,
\begin{equation}\label{jacolag}
\lim_{\beta \rightarrow \infty} P_{n}^{\left(\alpha ,\beta
\right)}(1-2x/\beta) = L^{\alpha}_n(x).\\
\end{equation}
\end{enumerate}
Both relations holds uniformly in every closed interval of
${\mathbb R}$.\\
\end{itemize}

Actually these relations are expression of deeper relations
between the measures and operators involved. As a consequence of those relations we have the following
technical results, that were proved in \cite{NaUrb}, and are
needed to prove Theorem \ref{L2resultfunctg}.\\
\vspace{0.2cm}

\begin{prop} \label{normrel}(norm relations)
\begin{enumerate}
\item[i)] Let  $f\in L^{2}(\mathbb{R},\gamma)$ and define
$f_{\lambda}(x)=f(\sqrt{\lambda}x)\chi_{[-1,1]}(x)$, then $
f_{\lambda}\in L^{2}([-1,1],\mu_{\lambda})$ and
\begin{equation}
\lim_{\lambda\to
\infty}\|f_{\lambda}\|_{2,\lambda}=\|f\|_{2,\gamma}
\end{equation}

\item[ii)] Let $f\in L^{2}(\mathbb{R},\mu_{\alpha})$ and define $
f_{\beta}(x)=f\left(\frac{\beta}{2}(1-x)\right)\chi_{[-1,1]}(x)$,
then $f_{\beta}\in L^{2}([-1,1],\mu_{(\alpha,\beta)})$ and
\begin{equation}
\lim_{\beta\to
\infty}\|f_{\beta}\|_{2,(\alpha,\beta)}=\|f\|_{2,\alpha}
\end{equation}
\end{enumerate}
\end{prop}
\begin{prop} (inner product relations) With the same notation as in Lemma \ref{normrel},
\begin{enumerate}
\item[i)] Let  $f\in L^{2}(\mathbb{R},\gamma),$ then
\begin{equation}
\displaystyle\lim_{\lambda\to\infty}\langle
f_{\lambda},\lambda^{-k/2}C^{\lambda}_{k}\rangle=\langle
f,\frac{H_{k}}{k!}\rangle.
\end{equation}
\item[ii)] Let  $f\in L^{2}(\mathbb{R},\mu_{\alpha}),$ then
\begin{equation}
\displaystyle\lim_{\beta\to\infty}\langle
f_{\beta},P_{k}^{(\alpha,\beta)} \rangle =\langle
f,L_{k}^{\alpha}\rangle. \\
\end{equation}
\end{enumerate}
\end{prop}

The results of this paper follows the same scheme of the proof given in \cite{NaUrb}, where this transference method was used to
obtain the the $L^p$ boundedness for the Riesz transform in the Hermite and Laguerre case from the $L^p$ boundedness for the Riesz transform in the Jacobi case. Unfortunately due to the non-linearity of the Littlewood-Paley g-function the computations for the  case $p \neq 2$ are more involved.

\section{Main Results}
We want to obtain the $L^p$-continuity for the
Gaussian-Littlewood-Paley $g$  and the $L^p$-continuity for the
Laguerre-Littlewood-Paley $g$  from the $L^p$-continuity of the
Jacobi-Littlewood-Paley $g$ , using a transference method based on
the asymptotic relations between Jacobi polynomials and Hermite
and Laguerre polynomials. We will start considering the case
$p=2$; more precisely we want to prove\\

\begin{theo} \label{L2resultfunctg}The $L^2(\mu_{\alpha, \beta})$ boundedness for the Jacobi-Littlewood-Paley $g$
\begin{equation}\label{JacobRieszL2ineq}
\| g^{(\alpha, \beta)} f \|_{2,(\alpha, \beta)} \leq C_2 \| f
\|_{2,(\alpha, \beta)}
\end{equation}
implies

\begin{enumerate}
\item[i)] the $L^2(\gamma)$  boundedness for the
Gaussian-Littlewood-Paley $g$
\begin{equation}\label{HermfunctgL2ineq}
\| g^{\gamma} f \|_{2,\gamma} \leq C_2 \| f \|_{2,\gamma}.
\end{equation}
\item[ii)] the $L^2(\mu_\alpha)$  boundedness for the
Laguerre-Littlewood-Paley $g$
 \begin{equation}\label{LagfunctgL2ineq}
\| g^{\alpha} f \|_{2,\alpha} \leq C_2 \| f \|_{2,\alpha}.
\end{equation}
\end{enumerate}
\end{theo}

\dem \\
\begin{enumerate}
\item[i)] Let $f\in L^{2}(\mathbb{R},\gamma)$ and define
$f_{\lambda}(x)=f(\sqrt{\lambda}x)$; then
\begin{eqnarray*}
 \left\|g^{(\alpha,\beta)}f_{\lambda}\right\|_{2,(\alpha,\beta)}^{2}&=&\int_{-1}^{1}\left(\int_{0}^{\infty}t|\nabla_{\alpha,\beta}
P^{(\alpha,\beta)}_{t}f_{\lambda}(x)|^{2}dt\right)d\mu_{\alpha,\beta}(x)\\
&=&\int_{0}^{\infty}t\left(\int_{-1}^{1}\left(\left(\frac{\partial}{\partial
t}P^{(\alpha,\beta)}_{t}f_{\lambda}(x)\right)^{2}
+(1-x^{2})\left(\frac{\partial}{\partial
x}P^{(\alpha,\beta)}_{t}f_{\lambda}(x)\right)^{2}\right)d\mu_{\alpha,\beta}(x)\right)dt
\end{eqnarray*}
Let us study the first term, by Parseval's identity,
\begin{eqnarray*}
\int_{-1}^{1}\left(\frac{\partial}{\partial
t}P^{(\alpha,\beta)}_{t}f_{\lambda}(x)\right)^{2}d\mu_{\alpha,\beta}(x) &=&\sum_{k=1}^{\infty}\left|\left\langle
f_{\lambda},\frac{P^{(\alpha,\beta)}_{k}}{\left\|P^{(\alpha,\beta)}_{k}\right\|_{2,(\alpha,\beta)}}\right\rangle\right|^{2}
\lambda_{k}e^{-2t\lambda_{k}^{1/2}},
\end{eqnarray*}
Therefore, interchanging the integral with the series, using
Lebesgue's domi\-na\-ted convergence theorem, we get
\begin{eqnarray*}
\int_{0}^{\infty}t\left(\int_{-1}^{1}\left(\frac{\partial}{\partial
t}P^{(\alpha,\beta)}_{t}f_{\lambda}(x)\right)^{2}d\mu_{\alpha,\beta}(x)\right)dt&=&\int_{0}^{\infty}t\left(\sum_{k=1}^{\infty}\left|\left\langle
f_{\lambda},\frac{P^{(\alpha,\beta)}_{k}}{\left\|P^{(\alpha,\beta)}_{k}\right\|_{2,(\alpha,\beta)}}\right\rangle\right|^{2}
\lambda_{k}e^{-2t\lambda_{k}^{1/2}}\right)dt\\
&=&\frac{1}{4}\sum_{k=1}^{\infty}\left|\left\langle
f_{\lambda},\frac{P^{(\alpha,\beta)}_{k}}{\left\|P^{(\alpha,\beta)}_{k}\right\|_{2,(\alpha,\beta)}}\right\rangle\right|^{2}
= \frac{1}{4} \| f_\lambda\|_{2,(\alpha,\beta)}^2
\end{eqnarray*}
In particular, taking $\alpha=\beta=\lambda-1/2$, we obtain
\begin{eqnarray*}
\int_{0}^{\infty}t\left(\int_{-1}^{1}\left(\frac{\partial}{\partial
t}P^{(\lambda-1/2,\lambda-1/2)}_{t}f_{\lambda}(x)\right)^{2}d\mu_{\lambda}(x)\right)dt&=&\frac{1}{4}\int_{-1}^{1}|f_{\lambda}(x)|^{2}d\mu_{\lambda}(x)
=\frac{1}{4}\left\|f_{\lambda}\right\|_{2,\lambda}^{2}
\end{eqnarray*}
Thus, we get using Proposition 1.1
\begin{eqnarray*}
&&\lim_{\lambda \to
\infty}\int_{0}^{\infty}t\left(\int_{-1}^{1}\left(\frac{\partial}{\partial
t}P^{(\lambda-1/2,\lambda-1/2)}_{t}f_{\lambda}(x)\right)^{2}d\mu_{\lambda}(x)\right)dt
=\lim_{\lambda \to
\infty}\frac{1}{4}\left\|f_{\lambda}\right\|_{2,\lambda}^{2}=\frac{1}{4}\left\|f\right\|_{2,\gamma}^{2}\\&=&
\frac{1}{4}\int_{-\infty}^{\infty}|f(x)|^{2}d\gamma(x)= \frac{1}{4}\sum_{k=1}^{\infty}\left|\left\langle
f,\frac{H_{k}}{\left\|H_{k}\right\|_{2,\gamma}}\right\rangle\right|^{2}
=\sum_{k=1}^{\infty}\left|\left\langle
f,\frac{H_{k}}{\left\|H_{k}\right\|_{2,\gamma}}\right\rangle\right|^{2}k\int_{0}^{\infty}te^{-2tk^{1/2}}dt\\
&=&\int_{0}^{\infty}t\left(\sum_{k=1}^{\infty}\left|\left\langle
f,\frac{H_{k}}{\left\|H_{k}\right\|_{2,\gamma}}\right\rangle\right|^{2}ke^{-2tk^{1/2}}\right)dt=
\int_{0}^{\infty}t\left(\int_{-\infty}^{\infty}\left(\frac{\partial}{\partial
t}P^{\gamma}_{t}f(x)\right) d\gamma(x)\right) dt
\end{eqnarray*}

Hence,
\begin{equation*}
\lim_{\lambda \to
\infty}\int_{0}^{\infty}t\left(\int_{-1}^{1}\left(\frac{\partial}{\partial
t}P^{(\lambda-1/2,\lambda-1/2)}_{t}f_{\lambda}(x)\right)^{2}d\mu_{\lambda}(x)\right)dt=\int_{0}^{\infty}t\left(\int_{-\infty}^{\infty}\left(\frac{\partial}{\partial
t}P^{\gamma}_{t}f(x)\right)^{2}d\gamma(x)\right)dt.
\end{equation*}
Now, let us consider the second term. First of all, observe that
\begin{eqnarray*}
&&\left\|\sqrt{1-x^{2}}P^{(\alpha+1,\beta+1)}_{k-1}\right\|^{2}_{2,(\alpha,\beta)}=\hspace{15cm}\\
&&=\frac{4(\alpha+1)(\beta+1)}{2^{\alpha+\beta+3}(\alpha+\beta+3)(\alpha+\beta+2)B(\alpha+2,\beta+2)}\\
&&\quad \quad \quad \quad \quad \quad \quad \quad \quad \quad \times \int^{1}_{-1}(1-x)^{\alpha+1}(1+x)^{\beta+1}\left[P^{(\alpha+1,\beta+1)}_{k-1}(x)\right]^{2}dx \hspace{7.1cm}\\
&&=\frac{4(\alpha+1)(\beta+1)}{(\alpha+\beta+3)(\alpha+\beta+2)}\left\|P^{(\alpha+1,\beta+1)}_{k-1}\right\|^{2}_{2,(\alpha+1,\beta+1)}=\frac{4k}{\left(k+\alpha+\beta+1\right)}\left\|P^{(\alpha,\beta)}_{k}\right\|^{2}_{2,(\alpha,\beta)}.
\end{eqnarray*}
Then, using Parseval's identity, we get
\begin{eqnarray*}
&&\int_{-1}^{1}\left(\sqrt{1-x^{2}}\frac{\partial}{\partial
x}P^{(\alpha,\beta)}_{t}f_{\lambda}(x)\right)^{2}d\mu_{\alpha,\beta}(x)=\left\|\sqrt{1-x^{2}}\frac{\partial}{\partial
x}P^{(\alpha,\beta)}_{t}f_{\lambda}\right\|_{2,(\alpha,\beta)}^{2}\\
   &=&\sum_{k=1}^{\infty}\left|\left\langle f_{\lambda},\frac{P^{(\alpha,\beta)}_{k}}{\left\|P^{(\alpha,\beta)}_{k}\right\|_{2,(\alpha,\beta)}^{2}}\right\rangle\right|^{2}
 e^{-2t\lambda_{k}^{1/2}}\frac{(k+\alpha+\beta+1)^{2}}{4}\left\|\sqrt{1-x^{2}}P^{(\alpha+1,\beta+1)}_{k-1}\right\|_{2,(\alpha,\beta)}^{2}\\
 &=&\sum_{k=1}^{\infty}\left|\left\langle f_{\lambda},\frac{P^{(\alpha,\beta)}_{k}}{\left\|P^{(\alpha,\beta)}_{k}\right\|_{2}^{2}}\right\rangle\right|^{2}
 e^{-2t\lambda_{k}^{1/2}}\lambda_{k}\left\|P^{(\alpha,\beta)}_{k}\right\|_{2,(\alpha,\beta)}^{2}=\sum_{k=1}^{\infty}\left|\left\langle f_{\lambda},\frac{P^{(\alpha,\beta)}_{k}}{\left\|P^{(\alpha,\beta)}_{k}\right\|_{2}}\right\rangle\right|^{2}
\lambda_{k} e^{-2t\lambda_{k}^{1/2}}.
\end{eqnarray*}
Now as before, using Lebesgue's dominated convergence theorem,
interchanging the integral with the series, we get
\begin{eqnarray*}
&&\int_{0}^{\infty}t\left(\int_{-1}^{1}\left(\sqrt{1-x^{2}}\frac{\partial}{\partial
x}P^{(\alpha,\beta)}_{t}f_{\lambda}(x)\right)^{2}d\mu_{\alpha,\beta}(x)\right)dt\hspace{6cm}\\
&=&\int_{0}^{\infty}t\left(\sum_{k=1}^{\infty}\left|\left\langle
f_{\lambda},\frac{P^{(\alpha,\beta)}_{k}}{\left\|P^{(\alpha,\beta)}_{k}\right\|_{2,(\alpha,\beta)}}\right\rangle\right|^{2}
\lambda_{k}e^{-2t\lambda_{k}^{1/2}}\right)dt
=\frac{1}{4}\int_{-1}^{1}|f_{\lambda}(x)|^{2}d\mu_{\alpha,\beta}(x).
\end{eqnarray*}
Then, for the Gegenbauer case, $\alpha=\beta=\lambda-1/2$  we get
\begin{equation*}
\int_{0}^{\infty}t\left(\int_{-1}^{1}\left(\sqrt{1-x^{2}}\frac{\partial}{\partial
x}P^{(\lambda-1/2,\lambda-1/2)}_{t}f_{\lambda}(x)\right)^{2}d\mu_{\lambda}(x)\right)dt=\frac{1}{4}\left\|f_{\lambda}\right\|_{2,\lambda}^{2}.
\end{equation*}
On the other hand, again using Parserval's identity, we have
\begin{eqnarray*}
&&\int_{-\infty}^{\infty}\left(\frac{1}{\sqrt{2}}\frac{\partial}{\partial
x}P^{\gamma}_{t}f(x)\right)^{2}d\gamma(x)=\sum_{k=1}^{\infty}\left|\left\langle
f,\frac{H_{k}}{\left\|H_{k}\right\|_{2,\gamma}^{2}}\right\rangle\right|^{2}2k^{2}e^{-2tk^{1/2}}\left\|H_{k-1}\right\|_{2,\gamma}^{2}\\&=&\sum_{k=1}^{\infty}\left|\left\langle
f,\frac{H_{k}}{\left\|H_{k}\right\|_{2,\gamma}^{2}}\right\rangle\right|^{2}2k^{2}e^{-2tk^{1/2}}\frac{1}{2k}\left\|H_{k}\right\|_{2,\gamma}^{2}
=\sum_{k=1}^{\infty}\left|\left\langle
f,\frac{H_{k}}{\left\|H_{k}\right\|_{2,\gamma}}\right\rangle\right|^{2}ke^{-2tk^{1/2}}.
\end{eqnarray*}
Hence,
\begin{eqnarray*}
\int_{0}^{\infty}t\left(\int_{-\infty}^{\infty}\left(\frac{1}{\sqrt{2}}\frac{\partial}{\partial
x}P^{\gamma}_{t}f(x)\right)^{2}d\gamma(x)\right)dt&=&\int_{0}^{\infty}t\left(\sum_{k=1}^{\infty}\left|\left\langle
f,\frac{H_{k}}{\left\|H_{k}\right\|_{2,\gamma}}\right\rangle\right|^{2}ke^{-2tk^{1/2}}\right)dt\\&=&\frac{1}{4}\left\|f\right\|_{2,\gamma}^{2}.
\end{eqnarray*}
Therefore, using Proposition 1.1,
\begin{eqnarray*}
&&\lim_{\lambda \to
\infty}\int_{0}^{\infty}t\left(\int_{-1}^{1}\left(\sqrt{1-x^{2}}\frac{\partial}{\partial
x}P^{(\lambda-1/2,\lambda-1/2)}_{t}f_{\lambda}(x)\right)^{2}d\mu_{\lambda}(x)\right)dt\\
&& \quad \quad =\lim_{\lambda \to
\infty}\frac{1}{4}\left\|f_{\lambda}\right\|_{2,\lambda}^{2}=\frac{1}{4}\left\|f\right\|_{2,\gamma}^{2}=\int_{0}^{\infty}t\left(\int_{-\infty}^{\infty}\left(\frac{1}{\sqrt{2}}\frac{\partial}{\partial
x}P^{\gamma}_{t}f(x)\right)^{2}d\gamma(x)\right)dt.
\end{eqnarray*}
Now, taking $\alpha=\beta=\lambda-1/2$,
\begin{eqnarray*}
\lim_{\lambda \to
\infty}\left\|g^{(\lambda-1/2,\lambda-1/2)}f_{\lambda}\right\|_{2,\lambda}^{2}&=&\lim_{\lambda
\to
\infty}\int_{0}^{\infty}t\left(\int_{-1}^{1}\left(\frac{\partial}{\partial
t}P^{(\lambda-1/2,\lambda-1/2)}_{t}f_{\lambda}(x)\right)^{2}d\mu_{\lambda}(x)\right)dt\\
&&+\lim_{\lambda \to
\infty}\int_{0}^{\infty}t\left(\int_{-1}^{1}\left(\sqrt{1-x^{2}}\frac{\partial}{\partial
x}P^{(\lambda-1/2,\lambda-1/2)}_{t}f_{\lambda}(x)\right)^{2}d\mu_{\lambda}(x)\right)dt\\
&=&\int_{0}^{\infty}t\left(\int_{-\infty}^{\infty}\left(\frac{\partial}{\partial
t}P^{\gamma}_{t}f(x)\right)^{2}d\gamma(x)\right)dt\\
&&\quad \quad
+\int_{0}^{\infty}t\left(\int_{-\infty}^{\infty}\left(\frac{1}{\sqrt{2}}\frac{\partial}{\partial
x}P^{\gamma}_{t}f(x)\right)^{2}d\gamma(x)\right)dt =
\left\|g^{\gamma}f\right\|_{2,\gamma}^{2}.
\end{eqnarray*}
now, by the $L^{2}$ continuity of $g^{(\lambda-1/2,\lambda-1/2)},$
we have
$$\left\|g^{(\lambda-1/2,\lambda-1/2)}f_{\lambda}\right\|_{2,\lambda}\leq C_{2}\left\|f_{\lambda}\right\|_{2,\lambda}$$
Finally,
$$\left\|g^{\gamma}f\right\|_{2,\gamma}=\lim_{\lambda \to
\infty}\left\|g^{(\lambda-1/2,\lambda-1/2)}f_{\lambda}\right\|_{2,\lambda}\leq
C_{2}\lim_{\lambda \to
\infty}\left\|f_{\lambda}\right\|_{2,\lambda}=C_{2}\left\|f\right\|_{2,\gamma}.\hspace{0.5cm\blacksquare}$$
 \item[ii)] Let $f\in L^{2}((0,\infty),\mu_{\alpha}),\;$ define
 $f_{\beta}(x)=f\left(\frac{\beta}{2}(1-x)\right)$, then analogously as in the Hermite case
\begin{equation*}
\left\|g^{(\alpha,\beta)}f_{\beta}\right\|_{2,(\alpha,\beta)}^{2}=\int_{0}^{\infty}t\left(\int_{-1}^{1}\left(\left(\frac{\partial}{\partial
t}P^{(\alpha,\beta)}_{t}f_{\beta}(x)\right)^{2}
+(1-x^{2})\left(\frac{\partial}{\partial
x}P^{(\alpha,\beta)}_{t}f_{\beta}(x)\right)^{2}\right)d\mu_{\alpha,\beta}(x)\right)dt,
\end{equation*}
where,
\begin{equation*}
\int_{0}^{\infty}t\left(\int_{-1}^{1}\left(\frac{\partial}{\partial
t}P^{(\alpha,\beta)}_{t}f_{\beta}(x)\right)^{2}d\mu_{\alpha,\beta}(x)\right)dt=\frac{1}{4}\int_{-1}^{1}|f_{\beta}(x)|^{2}d\mu_{\alpha,\beta}(x)
=\frac{1}{4}\left\|f_{\beta}\right\|_{2,({\alpha,\beta})}^{2}
\end{equation*}
and
\begin{equation*}
\int_{0}^{\infty}t\left(\int_{-1}^{1}\left(\sqrt{1-x^{2}}\frac{\partial}{\partial
x}P^{(\alpha,\beta)}_{t}f_{\lambda}(x)\right)^{2}d\mu_{\alpha,\beta}(x)\right)dt=
\frac{1}{4}\int_{-1}^{1}|f_{\beta}(x)|^{2}d\mu_{\alpha,\beta}(x)
=\frac{1}{4}\left\|f_{\beta}\right\|_{2,(\alpha,\beta)}^{2},
\end{equation*}
Let us study the first term. Using Proposition 1.1 and  Lebesgue's
dominated convergence theorem, interchanging the integral with the
series, we get
\begin{eqnarray*}
&&\lim_{\beta
\to\infty}\int_{0}^{\infty}t\left(\int_{-1}^{1}\left(\frac{\partial}{\partial
t}P^{(\alpha,\beta)}_{t}f_{\beta}(x)\right)^{2}d\mu_{\alpha,\beta}(x)\right)dt=\lim_{\beta
\to\infty}\frac{1}{4}\left\|f_{\beta}\right\|_{2,(\alpha,\beta)}^{2}=\frac{1}{4}\left\|f\right\|_{2,\alpha}^{2}\\
&&\hspace{0.5cm}=\int_{0}^{\infty}t\left(\sum_{k=1}^{\infty}\left|\left\langle
f,\frac{L_{k}^{\alpha}}{\left\|L_{k}^{\alpha}\right\|_{2,\alpha}}\right\rangle\right|^{2}ke^{-2tk^{1/2}}\right)dt=\int_{0}^{\infty}t\left(\int_{0}^{\infty}\left(\frac{\partial}{\partial
t}P^{\alpha}_{t}f(x)\right)^{2}d\mu_{\alpha}(x)\right)dt.
\end{eqnarray*}
Let us look now the second term. First of all observe that,
\begin{eqnarray*}
\left\|\sqrt{x}L_{k-1}^{\alpha+1}\right\|_{2,\alpha}^{2}&=&\int_{0}^{\infty}\left|\sqrt{x}L_{k-1}^{\alpha+1}(x)\right|^{2}\frac{x^{\alpha}e^{-x}}{\Gamma(\alpha+1)}dx
=(\alpha+1)\int_{0}^{\infty}\left|L_{k-1}^{\alpha+1}(x)\right|^{2}\frac{x^{\alpha+1}e^{-x}}{\Gamma(\alpha+2)}dx\\&=&(\alpha+1)\left\|L_{k-1}^{\alpha+1}\right\|_{2,(\alpha+1)}^{2}
=(\alpha+1)\frac{\Gamma(k-1+\alpha+1+1)}{\Gamma(\alpha+2)(k-1)!}=k\left\|L_{k}^{\alpha}\right\|_{2,\alpha}^{2},
\end{eqnarray*}
then we get, using Parseval's identity
\begin{eqnarray*}
\int_{0}^{\infty}\left(\sqrt{x}\frac{\partial}{\partial
x}P^{\alpha}_{t}f(x)\right)^{2}d\mu_{\alpha}(x)
&=&\sum_{k=1}^{\infty}\left|\left\langle
f,\frac{L_{k}^{\alpha}}{\left\|L_{k}^{\alpha}\right\|_{2,\alpha}^{2}}\right\rangle\right|^{2}
e^{-2tk^{1/2}}\left\|\sqrt{x}L_{k-1}^{\alpha+1}\right\|_{2,\alpha}^{2}\\
&=&\sum_{k=1}^{\infty}\left|\left\langle
f,\frac{L_{k}^{\alpha}}{\left\|L_{k}^{\alpha}\right\|_{2,\alpha}}\right\rangle\right|^{2}
ke^{-2tk^{1/2}}.\\
\end{eqnarray*}
Hence,
\begin{eqnarray*}
&&\lim_{\beta
\to\infty}\int_{0}^{\infty}t\left(\int_{-1}^{1}\left(\sqrt{1-x^{2}}\frac{\partial}{\partial
x}P^{(\alpha,\beta)}_{t}f_{\beta}(x)\right)^{2}d\mu_{\alpha,\beta}(x)\right)dt=\lim_{\beta
\to\infty}\frac{1}{4}\left\|f_{\beta}\right\|_{2,(\alpha,\beta)}^{2}
=\frac{1}{4}\left\|f\right\|_{2,\alpha}^{2}\\&&=\int_{0}^{\infty}t\left(\sum_{k=1}^{\infty}\left|\left\langle
f,\frac{L_{k}^{\alpha}}{\left\|L_{k}^{\alpha}\right\|_{2,\alpha}}\right\rangle\right|^{2}ke^{-2tk^{1/2}}\right)dt=\int_{0}^{\infty}t\left(\int_{0}^{\infty}\left(\sqrt{x}\frac{\partial}{\partial
x}P^{\alpha}_{t}f(x)\right)^{2}d\mu_{\alpha}(x)\right)dt.
\end{eqnarray*}
Therefore,
\begin{eqnarray*}
\lim_{\beta
\to\infty}\left\|g^{(\alpha,\beta)}f_{\beta}\right\|_{2,(\alpha,\beta)}^{2}&=&\lim_{\beta
\to\infty}\int_{0}^{\infty}t\left(\int_{-1}^{1}\left(\frac{\partial}{\partial
t}P^{(\alpha,\beta)}_{t}f_{\beta}(x)\right)^{2}d\mu_{\alpha,\beta}(x)\right)dt\\
&&\hspace{0.5cm}+\lim_{\beta
\to\infty}\int_{0}^{\infty}t\left(\int_{-1}^{1}\left(\sqrt{1-x^{2}}\frac{\partial}{\partial
x}P^{(\alpha,\beta)}_{t}f_{\beta}(x)\right)^{2}d\mu_{\alpha,\beta}(x)\right)dt\\
&=&\int_{0}^{\infty}t\left(\int_{0}^{\infty}\left(\frac{\partial}{\partial
t}P^{\alpha}_{t}f(x)\right)^{2}d\mu_{\alpha}(x)\right)dt\\
&&\hspace{0.5cm}+\int_{0}^{\infty}t\left(\int_{0}^{\infty}\left(\sqrt{x}\frac{\partial}{\partial
x}P^{\alpha}_{t}f(x)\right)^{2}d\mu_{\alpha}(x)\right)dt=\left\|g^{\alpha}f\right\|_{2,\alpha}^{2}
\end{eqnarray*}
and now, by the $L^{2}$ continuity of $g^{(\alpha,\beta)},$ we
have
$$\left\|g^{(\alpha,\beta)}f_{\beta}\right\|_{2,(\alpha,\beta)}\leq C_{2}\left\|f_{\beta}\right\|_{2,(\alpha,\beta)}$$
Thus finally,
\begin{equation*}
\left\|g^{\alpha}f\right\|_{2,\alpha}=\lim_{\beta
\to\infty}\left\|g^{(\alpha,\beta)}f_{\beta}\right\|_{2,(\alpha,\beta)}\leq
C_{2} \lim_{\beta
\to\infty}\left\|f_{\beta}\right\|_{2,(\alpha,\beta)}=C_{2}\left\|f\right\|_{2,\alpha}.\hspace{0.5cm}\blacksquare\\
\end{equation*}
\end{enumerate}

Now we are going to consider the general case $p\neq 2.$ For the
proof we will follow the argument given by Betancour et al in \cite{Bet}.\\

\begin{theo}\label{Lpresult}
Let $\alpha,\beta>-1 $ and $1<p<\infty$, then the
$L^p(\mu_{\alpha, \beta})$ boundedness for the
Jacobi-Littlewood-Paley $g$ function
\begin{equation}
\| g^{(\alpha, \beta)} f \|_{p,(\alpha, \beta)} \leq C_p \| f
\|_{p,(\alpha, \beta)}
\end{equation}
implies
\begin{enumerate}
\item[i)] The $L^p(\gamma)$-boundedness for the
Gaussian-Littlewood-Paley $g$ function
\begin{equation}
\| g^{\gamma} f \|_{p,\gamma} \leq C_p \| f \|_{p,\gamma}.
\end{equation}
and\\
\item[ii)] The $L^p(\mu_\alpha)$-boundedness for the
Laguerre-Littlewood-Paley $g$ function
 \begin{equation}
\| g^{\alpha} f \|_{p,\alpha} \leq C_p \| f \|_{p,\alpha}.
\end{equation}
\end{enumerate}

\end{theo}

\dem \\
\begin{enumerate}
\item [i)] Assume that the operator
$g^{(\lambda-1/2,\lambda-1/2)}$ is bounded in
$L^{p}\left([-1,1],\mu_{\lambda}\right)$. Let $\phi\in
C^{\infty}_{0}(\mathbb{R}),$ for each $\lambda>0$ define the
function
 $$\phi_{\lambda}(x)=\phi(\sqrt{\lambda}x),\;$$
$x\in \mathbb{R}.\;$ Let $\lambda$ big enough such that
$supp\,\phi_{\lambda}$ is contained in
$[-1,1].$ In what follows  $\;\lambda\;$ will be taken satisfying that condition. \\

Now, from the boundedness of $g^{(\lambda-1/2,\lambda-1/2)}$ we
have
\begin{eqnarray*}
\|g^{(\lambda-1/2,\lambda-1/2)}\phi_{\lambda}\|_{p, \lambda}\leq
C\|\phi_{\lambda}\|_{p, \lambda},
\end{eqnarray*}
that is to say,
\begin{eqnarray*}
\left\|\left[\int_{0}^{\infty}t|\nabla_{\lambda}
P^{(\lambda-1/2,\lambda-1/2)}_{t}\phi_{\lambda}|^{2}dt
\right]^{1/2}\right\|_{p, \lambda}\leq C\|\phi_{\lambda}\|_{p,
\lambda}
\end{eqnarray*}
where $\nabla_{\lambda}=\nabla_{(\lambda-1/2,\lambda-1/2)}.$ Now
making the change of variables $x=\frac{y}{\sqrt{\lambda}}$ and
taking
$Z(\lambda)=\frac{\lambda^{1/2}\left[\Gamma(\lambda)\right]^{2}2^{2\lambda}}{2\pi\Gamma\left(2\lambda\right)}$
we have
\begin{eqnarray*}
\left\{\int_{-\sqrt{\lambda}}^{\sqrt{\lambda}}\left|\left[\int_{0}^{\infty}t|\nabla_{\lambda}
P^{(\lambda-1/2,\lambda-1/2)}_{t}\phi_{\lambda}(\frac{y}{\sqrt{\lambda}})|^{2}dt
\right]^{1/2}\right|^{p}Z(\lambda)(1-\frac{y^{2}}{\lambda})^{\lambda-1/2}dy\right\}^{1/p}\hspace{3cm}
\end{eqnarray*}
$$\hspace{9cm} \leq
C\|\phi_{\lambda}\|_{p, \lambda},$$ which implies
\begin{eqnarray*}
\left\{\int_{-\sqrt{\lambda}}^{\sqrt{\lambda}}\left|\left[\int_{0}^{\infty}t|\nabla_{\lambda}
P^{(\lambda-1/2,\lambda-1/2)}_{t}\phi_{\lambda}(\frac{y}{\sqrt{\lambda}})|^{2}dt
\right]^{1/2}(1-\frac{y^{2}}{\lambda})^{\lambda/p-1/2p}e^{\frac{y^{2}}{p}}\right|^{p}\frac{e^{-y^{2}}}{\sqrt{\pi}}dy\right\}^{1/p}\hspace{3cm}
\end{eqnarray*}
$$\hspace{9cm} \leq
C(Z(\lambda))^{-1/p}\|\phi_{\lambda}\|_{p, \lambda}.$$
On the
other hand, analogously we have from the case $p = 2$,
\begin{eqnarray*}
\left\{\int_{-\sqrt{\lambda}}^{\sqrt{\lambda}}\left|\left[\int_{0}^{\infty}t|\nabla_{\lambda}
P^{(\lambda-1/2,\lambda-1/2)}_{t}\phi_{\lambda}(\frac{y}{\sqrt{\lambda}})|^{2}dt
\right]^{1/2}(1-\frac{y^{2}}{\lambda})^{\lambda/2-1/4}e^{\frac{y^{2}}{2}}\right|^{2}\frac{e^{-y^{2}}}{\sqrt{\pi}}dy\right\}^{1/2}\hspace{3cm}
\end{eqnarray*}
$$\hspace{9cm} \leq
C(Z(\lambda))^{-1/2}\|\phi_{\lambda}\|_{2, \lambda},$$ Now, define
for any $K\in\mathbb{N}\;\;$ and $\lambda>0$ such that
$\sqrt{\lambda}>K,\;$ the functions
\begin{eqnarray*}
F_{\lambda ,K}(y)=\left\{\begin{array}{lcl}
 \left[\int_{0}^{\infty}t|\nabla_{\lambda}
P^{(\lambda-1/2,\lambda-1/2)}_{t}\phi_{\lambda}(\frac{y}{\sqrt{\lambda}})|^{2}dt
\right]^{1/2}(1-\frac{y^{2}}{\lambda})^{\lambda/2-1/4}e^{\frac{y^{2}}{2}} & \mbox{if}& |y|\leq K\\
   \\
  0&\mbox{if}& |y|>K\\
\end{array}
\right.
\end{eqnarray*}
and
\begin{eqnarray*}
f_{\lambda ,K}(y)=\left\{\begin{array}{lcl}
 \left[\int_{0}^{\infty}t|\nabla_{\lambda}
P^{(\lambda-1/2,\lambda-1/2)}_{t}\phi_{\lambda}(\frac{y}{\sqrt{\lambda}})|^{2}dt
\right]^{1/2}(1-\frac{y^{2}}{\lambda})^{\lambda/p-1/2p}e^{\frac{y^{2}}{p}} &\mbox{if}& |y|\leq K\\
   \\
  0&\mbox{if}& |y|>K.\\
\end{array}
\right.
\end{eqnarray*}
Observe that  $F_{\lambda , K}=f_{\lambda, K}\Omega_{\lambda},\;$
where
$$\Omega_{\lambda}(y)=e^{\frac{y^{2}}{2}-\frac{y^{2}}{p}}\left(1-\frac{y^{2}}{\lambda}\right)^{\lambda/2-\lambda/p-1/4+1/2p},$$
for all $K\in\mathbb{N}\;\;
 $ and $\sqrt{\lambda}>K.\;$ Moreover,
 \begin{eqnarray*}
|\Omega_{\lambda}(y)|\leq
e^{\frac{y^{2}}{2}-\frac{y^{2}}{p}}e^{-\frac{y^{2}}{\lambda}(\lambda/2-\lambda/p-1/4+1/2p)}=e^{y^{2}(\frac{1}{4\lambda}-\frac{1}{2p\lambda})}.\\
\end{eqnarray*}
Therefore, if  $p\leq 2$ $\hspace{0.3cm}|\Omega_{\lambda}(y)|\leq
1$ and for $p>2$
$$|\Omega_{\lambda}(y)|\leq
e^{K^{2}(\frac{1}{4\lambda}-\frac{1}{2p\lambda})},$$
 hence $\Omega_{\lambda}$ is bounded in $[-K,K].$ Now
\begin{eqnarray*}
(Z(\lambda))^{-1/p}\|\phi_{\lambda}\|_{p,
\lambda}&=&(Z(\lambda))^{-1/p}\left\{\int_{-1}^{1}|\phi_{\lambda}(x)|^{p}\frac{\lambda\left[\Gamma(\lambda)\right]^{2}2^{2\lambda}}{2\pi\Gamma\left(2\lambda\right)}
\left(1-x^{2}\right)^{\lambda-1/2}\right\}^{1/p}dx.
\end{eqnarray*}
Thus, making the change of variables $x=\frac{y}{\sqrt{\lambda}}$
we get,
\begin{eqnarray*}
(Z(\lambda))^{-1/p}\|\phi_{\lambda}\|_{p,
\lambda}&=&(Z(\lambda))^{-1/p}\left\{\int_{-\sqrt{\lambda}}^{\sqrt{\lambda}}|\phi_{\lambda}(\frac{y}{\sqrt{\lambda}})|^{p}Z(\lambda)
\left(1-\frac{y^{2}}{\lambda}\right)^{\lambda-1/2}dy\right\}^{1/p}\\
&=&\left\{\int_{-\sqrt{\lambda}}^{\sqrt{\lambda}}|\phi(y)|^{p}\left(1-\frac{y^{2}}{\lambda}\right)^{\lambda-1/2}dy\right\}^{1/p}\leq
C\left\{\int_{-\sqrt{\lambda}}^{\sqrt{\lambda}}|\phi(y)|^{p}\frac{e^{-y^{2}}}{\sqrt{\pi}}dy\right\}^{1/p}.\\
\end{eqnarray*}
Therefore,
\begin{eqnarray*}
\lim_{\lambda\to \infty}(Z(\lambda))^{-1/p}\|\phi_{\lambda}\|_{p,
\lambda}&\leq & \lim_{\lambda\to
\infty}C\left\{\int_{-\sqrt{\lambda}}^{\sqrt{\lambda}}|\phi(y)|^{p}\frac{e^{-y^{2}}}{\sqrt{\pi}}dy\right\}^{1/p}\\
&=&C\left\{\int_{-\infty}^{\infty}|\phi(y)|^{p}\frac{e^{-y^{2}}}{\sqrt{\pi}}dy\right\}^{1/p}=C\|\phi\|_{p,\gamma},
\end{eqnarray*}
and moreover,
\begin{eqnarray}\label{desLpfi}
(Z(\lambda))^{-1/p}\|\phi_{\lambda}\|_{p, \lambda}\leq
C\|\phi\|_{p,\gamma}
\end{eqnarray}
On the other hand,
\begin{eqnarray}\label{norF1gaus}
\nonumber \|F_{\lambda
,K}\|_{2,\gamma}&=&\left\{\int_{\mathbb{R}}|F_{\lambda
,K}(y)|^{2}\frac{e^{-y^{2}}}{\sqrt{\pi}}dy\right\}^{1/2}=\left\{\int_{-K}^{K}|F_{\lambda
,K}(y)|^{2}\frac{e^{-y^{2}}}{\sqrt{\pi}}dy\right\}^{1/2}\\
 &\leq
&\left\{\int_{-\sqrt{\lambda}}^{\sqrt{\lambda}}|F_{\lambda
,K}(y)|^{2}\frac{e^{-y^{2}}}{\sqrt{\pi}}dy\right\}^{1/2}\leq
C(Z(\lambda))^{-1/2}\|\phi_{\lambda}\|_{2,\lambda}.
\end{eqnarray}
Then, from (\ref{desLpfi}) and (\ref{norF1gaus}), we have
\begin{eqnarray*}
\|F_{\lambda ,K}\|_{2,\gamma}\leq C\|\phi\|_{2,\gamma}.
\end{eqnarray*}
Analogously,
\begin{eqnarray*}
\nonumber \|F_{\lambda
,K}\|_{p,\gamma}&=&\left\{\int_{\mathbb{R}}|F_{\lambda
,K}(y)|^{p}\frac{e^{-y^{2}}}{\sqrt{\pi}}dy\right\}^{1/p}=\left\{\int_{-K}^{K}|f_{\lambda
,K}\Omega_{\lambda}(y)|^{p}\frac{e^{-y^{2}}}{\sqrt{\pi}}dy\right\}^{1/p}.\\
\end{eqnarray*}
Since $\Omega_{\lambda}$ is bounded in $[-k,k],$ we get
\begin{eqnarray}\label{norf1gaus}
\nonumber \|F_{\lambda ,K}\|_{p,\gamma} \nonumber \leq
C\left\{\int_{-\infty}^{\infty}|f_{\lambda
,K}(y)|^{p}\frac{e^{-y^{2}}}{\sqrt{\pi}}dy\right\}^{1/p}\leq
C(Z(\lambda))^{-1/p}\|\phi_{\lambda}\|_{p, \lambda}
\end{eqnarray}
Then, from (\ref{desLpfi}) and (\ref{norf1gaus}) we have
\begin{eqnarray*}
\|F_{\lambda ,K}\|_{p,\gamma}\leq C\|\phi\|_{p,\gamma}
\end{eqnarray*}
for all $\sqrt{\lambda}>K.\;$ Thus, $\;F_{\lambda ,K}\;\;$ is a
bounded sequence in $\;L^{2}\left(\mathbb{R},\gamma\right)\;$ and
$\;L^{p}\left(\mathbb{R},\gamma\right).\;$ By Bourbaki-Alaoglu's
theorem,there exists a subsequence
$\;(\lambda_{j})_{j\in\mathbb{N}}\;$ such that $\lim_{j\to
\infty}\lambda_{j}=\infty$ and functions $\;F_{K}\in
L^{2}\left(\mathbb{R},\gamma\right)$ and
$f_{K}\in{L^{p}\left(\mathbb{R},\gamma\right)}$ satisfying
\begin{itemize}
\item $ F_{\lambda_j,K}\to  F_{K},$ as $j\to\infty,$ in the weak
topology on ${L^{2}\left(\mathbb{R},\gamma\right)}$ \item
$F_{\lambda_j,K}\to  f_{K},$ as $j\to\infty,$ in the weak topology
on ${L^{p}\left(\mathbb{R},\gamma\right)}$
\end{itemize}
Moreover,  $\mbox{sop} F_{K}\cup \mbox{sop} f_{K}\subseteq
[-K,K],\;$ and
\begin{eqnarray}\label{AnormFk}
\|F_{K}\|_{2,\gamma}\leq\lim_{j\to
\infty}\|F_{\lambda_{j},K}\|_{2,\gamma}\leq C\|\phi\|_{2,\gamma}.
\end{eqnarray}
Analogously one gets,
\begin{eqnarray}\label{Anormfk}
\|f_{K}\|_{p,\gamma}\leq C\|\phi\|_{p,\gamma}
\end{eqnarray}
 Observe that, defining for every $k\in\mathbb{N}$,
$$\tau_K(g)=\int_{-\infty}^{\infty}g(x)\chi_{[-K,K]}(x)dx,$$
 then, by Cauchy-Schwartz inequality $\tau_K\in\left(L^{2}\left(\mathbb{R},\gamma\right)\right)^{\ast}$,
and therefore,
\begin{eqnarray*}
\int_{-K}^{K}F_{K}(x)dx&=&\int_{-\infty}^{\infty}F_{K}(x)\chi_{[-K,K]}(x)dx
= \tau_K(F_K)=\lim_{j\to \infty}  \tau_K(F_{\lambda_{j},K})
\\
&=&\lim_{j\to
\infty}\int_{-\infty}^{\infty}F_{\lambda_{j},K}(x)\chi_{[-K,K]}(x)dx= \int_{-K}^{K}f_{K}(x)dx; \\
\end{eqnarray*}
thus, $F_{K}=f_{K}\;\;a.e.\;\mbox{on} \;\;[-K,K],$ for all $K \in
\mathbb{N},$ so
 $F_{K}=f_{K}\;\;a.e.\;$ on $\mathbb{R}.$
Then from (\ref{Anormfk}), we get
\begin{eqnarray}\label{AnormFk-p}
\|F_{K}\|_{p,\gamma}\leq C\|\phi\|_{p,\gamma},
\end{eqnarray}
and therefore, from (\ref{AnormFk}) and (\ref{AnormFk-p}), there
exists an increasing sequence $\{\lambda_{j}\}_{j\in
\mathbb{N}}\subset (0,\infty)$, with $\lim_{j\to
\infty}\lambda_{j}=\infty,$ and a function $F\in
L^{p}\left(\mathbb{R},\gamma\right)\cap
L^{2}\left(\mathbb{R},\gamma\right),$ such that
\begin{itemize}
\item For each $K\in \mathbb{N},\;$  $ F_{\lambda_{j},K}\to F,$ as
$j\to\infty,$ in the weak topology of
${L^{2}\left(\mathbb{R},\gamma\right)}$ and in the weak topology
of ${L^{p}\left(\mathbb{R},\gamma\right)}$, and \item
$\|F\|_{p,\gamma}\leq C\|\phi\|_{p,\gamma}.$\\
\end{itemize}

On the other hand,
\begin{eqnarray*}
F_{\lambda ,K}(y) &=&\chi_{[-K,K]}(y)
\left(\int_{0}^{\infty}t|\nabla_{\lambda}
P^{(\lambda-1/2,\lambda-1/2)}_{t}\phi_{\lambda}(\frac{y}{\sqrt{\lambda}})|^{2}(1-\frac{y^{2}}{\lambda})^{\lambda-1/2}e^{y^{2}}dt
\right)^{1/2}\\
&=&\chi_{[-K,K]}(y)
\left[\int_{0}^{\infty}t\left(\frac{\partial}{\partial
t}P^{(\lambda-1/2,\lambda-1/2)}_{t}\phi_{\lambda}(\frac{y}{\sqrt{\lambda}})\right)^{2}dt(1-\frac{y^{2}}{\lambda})^{\lambda-1/2}e^{y^{2}}\right.\\&&\hspace{0.5cm}+\left.\int_{0}^{\infty}t(1-\frac{y^{2}}{\lambda})\left(\frac{\partial}{\partial
x}P^{(\lambda-1/2,\lambda-1/2)}_{t}\phi_{\lambda}(\frac{y}{\sqrt{\lambda}})\right)^{2}dt(1-\frac{y^{2}}{\lambda})^{\lambda-1/2}e^{y^{2}}
\right]^{1/2}.\\
\end{eqnarray*}
Now, define
$$g_{1,\lambda}\phi(y)=\chi_{[-K,K]}(y)\left(\int_{0}^{\infty}t\left(\frac{\partial}{\partial
t}P^{(\lambda-1/2,\lambda-1/2)}_{t}\phi_{\lambda}(\frac{y}{\sqrt{\lambda}})\right)^{2}dt\right)(1-\frac{y^{2}}{\lambda})^{\lambda-1/2}e^{y^{2}}$$
and
$$g_{2,\lambda}\phi(y)=\chi_{[-K,K]}(y)\left(\int_{0}^{\infty}t(1-\frac{y^{2}}{\lambda})\left(\frac{\partial}{\partial
y}P^{(\lambda-1/2,\lambda-1/2)}_{t}\phi_{\lambda}(\frac{y}{\sqrt{\lambda}})\right)^{2}dt\right)(1-\frac{y^{2}}{\lambda})^{\lambda-1/2}e^{y^{2}}.$$
Let us study first the function $g_{1,\lambda}$. For a function
$\phi\in L^{2}(\mathbb{R},\gamma)$,
\begin{eqnarray*}
\sqrt{t}\frac{\partial}{\partial
t}P^{(\alpha,\beta)}_{t}\phi_{\lambda}(\frac{y}{\sqrt{\lambda}})&=&-\sum_{k=1}^{\infty}\left\langle
\phi_{\lambda},\frac{P^{(\alpha,\beta)}_{k}}{\left\|P^{(\alpha,\beta)}_{k}\right\|_{2,
(\alpha,\beta)}^{2}}\right\rangle
(t\lambda_{k})^{1/2}e^{-t\lambda_{k}^{1/2}}P^{(\alpha,\beta)}_{k}(\frac{y}{\sqrt{\lambda}})
\end{eqnarray*}
which converges absolutely. Then, taking the Cauchy product, we
obtain,
\begin{eqnarray*}
&&\left(\sqrt{t}\frac{\partial}{\partial
t}P^{(\alpha,\beta)}_{t}\phi_{\lambda}(\frac{y}{\sqrt{\lambda}})\right)^{2} =\hspace{2cm}\\
&&\sum_{k=1}^{\infty}\sum_{n=1}^{k}\left\langle
\phi_{\lambda},\frac{P^{(\alpha,\beta)}_{n}}{\left\|P^{(\alpha,\beta)}_{n}\right\|_{2}^{2}}\right\rangle\left\langle
\phi_{\lambda},\frac{P^{(\alpha,\beta)}_{k-n}}{\left\|P^{(\alpha,\beta)}_{k-n}\right\|_{2}^{2}}\right\rangle
t\lambda_{n}^{1/2}\lambda_{k-n}^{1/2}e^{-t(\lambda_{k}^{1/2}+\lambda_{k-n}^{1/2})}P^{(\alpha,\beta)}_{n}(\frac{y}{\sqrt{\lambda}})P^{(\alpha,\beta)}_{k-n}(\frac{y}{\sqrt{\lambda}}).
\end{eqnarray*}
Then, for the Gegenbauer case, $\alpha=\beta=\lambda-1/2$, taking
$\hat{\lambda_{k}}=k(k+2\lambda)$
\begin{eqnarray*}
&&\left(\sqrt{t}\frac{\partial}{\partial
t}P^{(\lambda-1/2,\lambda-1/2)}_{t}\phi_{\lambda}(\frac{y}{\sqrt{\lambda}})\right)^{2} =\hspace{2cm}\\
&&\sum_{k=1}^{\infty}\sum_{n=1}^{k}\left\langle
\phi_{\lambda},C^{\lambda}_{n}\right\rangle\left\langle
\phi_{\lambda},C^{\lambda}_{k-n}\right\rangle
t\hat{\lambda}_{n}^{1/2}\hat{\lambda}_{k-n}^{1/2}e^{-t(\hat{\lambda}_{k}^{1/2}+\hat{\lambda}_{k-n}^{1/2})}C^{\lambda}_{n}(\frac{y}{\sqrt{\lambda}})C^{\lambda}_{k-n}(\frac{y}{\sqrt{\lambda}})\\
&&\;\times\left[\frac{\Gamma(2\lambda)\Gamma(n+\lambda+1/2)}{\Gamma(\lambda+1/2)\Gamma(n+2\lambda)}\frac{\Gamma(2\lambda)\Gamma(k-n+\lambda+1/2)}{\Gamma(\lambda+1/2)\Gamma(k-n+2\lambda)}\right]^{2}
\frac{2(n+\lambda)[\Gamma(\lambda+1/2)]^{2}n!\Gamma(n+2\lambda)}{\Gamma(2\lambda+1)\Gamma(n+\lambda+1/2)}\\
&& \hspace{5cm}\times\frac{2(k-n+\lambda)[\Gamma(\lambda+1/2)]^{2}(k-n)!\Gamma(k-n+2\lambda)}{\Gamma(2\lambda+1)\Gamma(k-n+\lambda+1/2)}\\
&=&\sum_{k=1}^{\infty}\sum_{n=1}^{k}\left\langle
\phi_{\lambda},C^{\lambda}_{n}\right\rangle\left\langle
\phi_{\lambda},C^{\lambda}_{k-n}\right\rangle
t\hat{\lambda}_{n}^{1/2}\hat{\lambda}_{k-n}^{1/2}e^{-t(\hat{\lambda}_{k}^{1/2}+\hat{\lambda}_{k-n}^{1/2})}C^{\lambda}_{n}(\frac{y}{\sqrt{\lambda}})C^{\lambda}_{k-n}(\frac{y}{\sqrt{\lambda}})\\
&&\hspace{5cm}\times\frac{(n+\lambda)n!(k-n)!(k-n+\lambda)}{\lambda^{2}(2\lambda)_{n}(2\lambda)_{k-n}}.\\
\end{eqnarray*}
Thus,
\begin{eqnarray*}
&&\left(\sqrt{t}\frac{\partial}{\partial
t}P^{(\lambda-1/2,\lambda-1/2)}_{t}\phi_{\lambda}(\frac{y}{\sqrt{\lambda}})\right)^{2}=\hspace{13cm}\\
&& \sum_{k=1}^{\infty}\sum_{n=1}^{k}\left\langle
\phi_{\lambda},\frac{C^{\lambda}_{n}}{\|C^{\lambda}_{n}\|^{2}}\right\rangle\left\langle
\phi_{\lambda},\frac{C^{\lambda}_{k-n}}{\|C^{\lambda}_{k-n}\|^{2}}\right\rangle
t\hat{\lambda}_{n}^{1/2}\hat{\lambda}_{k-n}^{1/2}e^{-t(\hat{\lambda}_{n}^{1/2}+\hat{\lambda}_{k-n}^{1/2})}C^{\lambda}_{n}(\frac{y}{\sqrt{\lambda}})C^{\lambda}_{k-n}(\frac{y}{\sqrt{\lambda}}).\\
\end{eqnarray*}
Therefore,
\begin{eqnarray*}
g_{1,\lambda}\phi(y)
&=&\chi_{[-K,K]}(y)\int_{0}^{\infty}\sum_{k=1}^{\infty}\sum_{n=1}^{k}\left\langle
\phi_{\lambda},\frac{C^{\lambda}_{n}}{\|C^{\lambda}_{n}\|^{2}}\right\rangle\left\langle
\phi_{\lambda},\frac{C^{\lambda}_{k-n}}{\|C^{\lambda}_{k-n}\|^{2}}\right\rangle t\hat{\lambda}_{n}^{1/2}\hat{\lambda}_{k-n}^{1/2}e^{-t(\hat{\lambda}_{n}^{1/2}+\hat{\lambda}_{k-n}^{1/2})}\\
&&\;\hspace{1.7cm}\times
C^{\lambda}_{n}(\frac{y}{\sqrt{\lambda}})C^{\lambda}_{k-n}(\frac{y}{\sqrt{\lambda}})(1-\frac{y^{2}}{\lambda})^{\lambda-1/2}e^{y^{2}}dt.
\end{eqnarray*}
Let us take,
\begin{eqnarray*}
H_{\lambda,K}^{N,1}(y)&=&\chi_{[-K,K]}(y)\sum_{k=N+1}^{\infty}\sum_{n=1}^{k}\left\langle
\phi_{\lambda},\frac{C^{\lambda}_{n}}{\|C^{\lambda}_{n}\|^{2}}\right\rangle\left\langle
\phi_{\lambda},\frac{C^{\lambda}_{k-n}}{\|C^{\lambda}_{k-n}\|^{2}}\right\rangle
t\hat{\lambda}_{n}^{1/2}\hat{\lambda}_{k-n}^{1/2}e^{-t(\hat{\lambda}_{n}^{1/2}+\hat{\lambda}_{k-n}^{1/2})}\\
&&\;\hspace{4cm}\times C^{\lambda}_{n}(\frac{y}{\sqrt{\lambda}})C^{\lambda}_{k-n}(\frac{y}{\sqrt{\lambda}})(1-\frac{y^{2}}{\lambda})^{\lambda/2-1/4}e^{\frac{y^{2}}{2}}.\\
\end{eqnarray*}
Thus,
\begin{eqnarray*}
g_{1,\lambda}\phi(y)&=&\chi_{[-K,K]}(y)\int_{0}^{\infty}\sum_{k=1}^{N}\sum_{n=1}^{k}\left\langle
\phi_{\lambda},\frac{C^{\lambda}_{n}}{\|C^{\lambda}_{n}\|^{2}}\right\rangle\left\langle
\phi_{\lambda},\frac{C^{\lambda}_{k-n}}{\|C^{\lambda}_{k-n}\|^{2}}\right\rangle
t\hat{\lambda}_{n}^{1/2}\hat{\lambda}_{k-n}^{1/2}e^{-t(\hat{\lambda}_{n}^{1/2}+\hat{\lambda}_{k-n}^{1/2})}\\
 &&\;\hspace{4cm}\times C^{\lambda}_{n}(\frac{y}{\sqrt{\lambda}})C^{\lambda}_{k-n}(\frac{y}{\sqrt{\lambda}})(1-\frac{y^{2}}{\lambda})^{\lambda-1/2}e^{y^{2}}dt\\
&&\;\hspace{0.3cm}+\;\;\chi_{[-K,K]}(y)\int_{0}^{\infty}\sum_{k=N+1}^{\infty}\sum_{n=1}^{k}\left\langle
\phi_{\lambda},\frac{C^{\lambda}_{n}}{\|C^{\lambda}_{n}\|^{2}}\right\rangle\left\langle
\phi_{\lambda},\frac{C^{\lambda}_{k-n}}{\|C^{\lambda}_{k-n}\|^{2}}\right\rangle
t\hat{\lambda}_{n}^{1/2}\hat{\lambda}_{k-n}^{1/2}e^{-t(\hat{\lambda}_{n}^{1/2}+\hat{\lambda}_{k-n}^{1/2})}\\
 &&\;\hspace{4cm}\times C^{\lambda}_{n}(\frac{y}{\sqrt{\lambda}})C^{\lambda}_{k-n}(\frac{y}{\sqrt{\lambda}})(1-\frac{y^{2}}{\lambda})^{\lambda-1/2}e^{y^{2}}dt\\
 \end{eqnarray*}
 \begin{eqnarray*}
 &=&\sum_{k=1}^{N}\sum_{n=1}^{k}\left\langle
\phi_{\lambda},\frac{C^{\lambda}_{n}}{\|C^{\lambda}_{n}\|^{2}}\right\rangle\left\langle
\phi_{\lambda},\frac{C^{\lambda}_{k-n}}{\|C^{\lambda}_{k-n}\|^{2}}\right\rangle
\frac{\hat{\lambda}_{n}^{1/2}\hat{\lambda}_{k-n}^{1/2}}{(\hat{\lambda}_{n}^{1/2}+\hat{\lambda}_{k-n}^{1/2})^{2}} C^{\lambda}_{n}(\frac{y}{\sqrt{\lambda}})C^{\lambda}_{k-n}(\frac{y}{\sqrt{\lambda}})\\
&&\;\hspace{2.5cm}\times(1-\frac{y^{2}}{\lambda})^{\lambda-1/2}e^{y^{2}}\;\;+\;\;\int_{0}^{\infty}H_{\lambda,K}^{N,1}(y)dt(1-\frac{y^{2}}{\lambda})^{\lambda/2-1/4}e^{\frac{y^{2}}{2}}.
\end{eqnarray*}
Now, we want to prove that for $K\in\mathbb{N}$ and $\lambda>0\;$
such that $\sqrt{\lambda}>K$ ,
\begin{equation*}
\int_{-\infty}^{\infty}\left|\int_{0}^{\infty}H_{\lambda,K}^{N,1}(y)dt\right|^{2}e^{-y^{2}}\frac{dy}{\sqrt{\pi}}\leq
\frac{\theta}{e^{(N+1)}}.
\end{equation*}
Indeed, by Minkowski integral inequality we have
\begin{eqnarray*}
\left(\int_{-\infty}^{\infty}\left|\int_{0}^{\infty}H_{\lambda,K}^{N,1}(y)dt\right|^{2}e^{-y^{2}}\frac{dy}{\sqrt{\pi}}\right)^{1/2}&\leq
&\int_{0}^{\infty}\left(\int_{-\infty}^{\infty}\left|H_{\lambda,K}^{N,1}(y)\right|^{2}e^{-y^{2}}\frac{dy}{\sqrt{\pi}}\right)^{1/2}dt,
\end{eqnarray*}
then,
\begin{eqnarray*}
&&\int_{-\infty}^{\infty}\left|H_{\lambda,K}^{N,1}(y)\right|^{2}e^{-y^{2}}\frac{dy}{\sqrt{\pi}}\\
&&\hspace{1cm} \leq\frac{1}{\sqrt{\pi}}\int_{-1}^{1}\left|\sum_{k=N+1}^{\infty}\sum_{n=1}^{k}\left\langle
 \phi_{\lambda},\frac{C^{\lambda}_{n}}{\|C^{\lambda}_{n}\|^{2}}\right\rangle\left\langle
\phi_{\lambda},\frac{C^{\lambda}_{k-n}}{\|C^{\lambda}_{k-n}\|^{2}}\right\rangle t\hat{\lambda}_{n}^{1/2}\hat{\lambda}_{k-n}^{1/2}e^{-t(\hat{\lambda}_{k}^{1/2}+\hat{\lambda}_{k-n}^{1/2})}\right.\\
&&\;\hspace{7cm}\times \left.C^{\lambda}_{n}(x)C^{\lambda}_{k-n}(x)\right|^{2}(1-x^{2})^{\lambda-1/2}dx.\\
\end{eqnarray*}
Using Gasper's  linearization of the product of Jacobi
polynomials, see \cite{gasp1},
\begin{eqnarray*}
C^{\lambda}_{n}(x)C^{\lambda}_{k-n}(x)&=&\frac{[\Gamma(\lambda+1/2)]^{2}\Gamma(n+2\lambda)\Gamma(k-n+2\lambda)}{[\Gamma(2\lambda)]^{2}\Gamma(n+\lambda+1/2)\Gamma(k-n+\lambda+1/2)}
\\&& \hspace{0.5cm} \times
P^{(\lambda-1/2,\lambda-1/2)}_{n}(1)P^{(\lambda-1/2,\lambda-1/2)}_{k-n}(1)p^{(\lambda-1/2,\lambda-1/2)}_{n}(x)p^{(\lambda-1/2,\lambda-1/2)}_{k-n}(x)\\
&=&\frac{[\Gamma(\lambda+1/2)]^{2}\Gamma(n+2\lambda)\Gamma(k-n+2\lambda)}{[\Gamma(2\lambda)]^{2}\Gamma(n+\lambda+1/2)\Gamma(k-n+\lambda+1/2)}
\\&& \hspace{0.5cm} \times
P^{(\lambda-1/2,\lambda-1/2)}_{n}(1)P^{(\lambda-1/2,\lambda-1/2)}_{k-n}(1)\sum_{i=|k-2n|}^{k}\nu(i,n,k-n)p^{(\lambda-1/2,\lambda-1/2)}_{i}(x),
\end{eqnarray*}
where,
$p^{(\alpha,\beta)}_{i}(x)=P^{(\alpha,\beta)}_{i}(x)/P^{(\alpha,\beta)}_{i}(1)\;\;$
and
\begin{equation*}
\nu(i,k-n,n)=\frac{1}{\left\|p^{(\alpha,\beta)}_{i}\right\|_{2,(\alpha,\beta)}^{2}}\int^{1}_{-1}p^{(\alpha,\beta)}_{i}(x)p^{(\alpha,\beta)}_{k-n}(x)p^{(\alpha,\beta)}_{n}(x)d\mu_{(\alpha,\beta)}(x).
\end{equation*}
Hence,
\begin{eqnarray*}
&&\int_{-\infty}^{\infty}\left|H_{\lambda,K}^{N,1}(y)\right|^{2}e^{-y^{2}}\frac{dy}{\sqrt{\pi}}\\
&& \hspace{0.5cm}\leq
\frac{1}{\sqrt{\pi}}\int_{-1}^{1}\left|\sum_{k=N+1}^{\infty}\sum_{n=1}^{k}\sum_{i=|k-2n|}^{k}\left\langle
\phi_{\lambda},\frac{C^{\lambda}_{n}}{\|C^{\lambda}_{n}\|^{2}}\right\rangle\left\langle
\phi_{\lambda},\frac{C^{\lambda}_{k-n}}{\|C^{\lambda}_{k-n}\|^{2}}\right\rangle t\hat{\lambda}_{n}^{1/2}\hat{\lambda}_{k-n}^{1/2}e^{-t(\hat{\lambda}_{k}^{1/2}+\hat{\lambda}_{k-n}^{1/2})}\right.\\
&&\hspace{0.7cm}\times
\frac{[\Gamma(\lambda+1/2)]^{2}\Gamma(n+2\lambda)\Gamma(k-n+2\lambda)}{[\Gamma(2\lambda)]^{2}\Gamma(n+\lambda+1/2)\Gamma(k-n+\lambda+1/2)}\\
&&\;\hspace{0.8cm}\times \left.P^{(\lambda-1/2,\lambda-1/2)}_{n}(1)P^{(\lambda-1/2,\lambda-1/2)}_{k-n}(1)\nu(i,n,k-n)p^{(\lambda-1/2,\lambda-1/2)}_{i}(x)\right|^{2}(1-x^{2})^{\lambda-1/2}dx.\\
\end{eqnarray*}
Then, by Parseval's identity and Cauchy-Schwartz inequality, we
have
\begin{eqnarray*}
&&\int_{-\infty}^{\infty}\left|H_{\lambda,K}^{N,1}(y)\right|^{2}e^{-y^{2}}\frac{dy}{\sqrt{\pi}}\\
&& \hspace{0.4cm}\leq
\frac{1}{\sqrt{\pi}}\frac{2\pi\Gamma(2\lambda)}{\lambda[\Gamma(\lambda)]^{2}2^{2\lambda}}\sum_{k=N+1}^{\infty}\sum_{n=1}^{k}\sum_{i=|k-2n|}^{k}\frac{\|\phi_{\lambda}\|^{2}}{\|C^{\lambda}_{k-n}\|^{2}\|C^{\lambda}_{n}\|^{2}}\left|\left\langle
\phi_{\lambda},\frac{C^{\lambda}_{n}}{\|C^{\lambda}_{n}\|}\right\rangle\right|^{2} t^{2}\hat{\lambda}_{n}\hat{\lambda}_{k-n}e^{-2t(\hat{\lambda}_{k}^{1/2}+\hat{\lambda}_{k-n}^{1/2})}\\
&&\hspace{0.5cm}\times
\left|\frac{[\Gamma(\lambda+1/2)]^{2}\Gamma(n+2\lambda)\Gamma(k-n+2\lambda)}{[\Gamma(2\lambda)]^{2}\Gamma(n+\lambda+1/2)\Gamma(k-n+\lambda+1/2)}\right|^{2}\\
&&\;\hspace{0.8cm}\times \left|P^{(\lambda-1/2,\lambda-1/2)}_{n}(1)P^{(\lambda-1/2,\lambda-1/2)}_{k-n}(1)\right|^{2}\left|\nu(i,n,k-n)\right|^{2}\left\|p^{(\lambda-1/2,\lambda-1/2)}_{i}\right\|^{2}\\
\end{eqnarray*}
Now,
\begin{eqnarray}
\nonumber\left|\nu(i,k-n,n)\left\|p^{(\alpha,\beta)}_{i}\right\|_{2,(\alpha,\beta)}\right|&\leq
&
\frac{1}{\left\|p^{(\alpha,\beta)}_{i}\right\|_{2, (\alpha,\beta)}}\int^{1}_{-1}\left|p^{(\alpha,\beta)}_{i}(x)p^{(\alpha,\beta)}_{k-n}(x)\right|\left|p^{(\alpha,\beta)}_{n}(x)\right|d\mu_{(\alpha,\beta)}(x)\\
\nonumber &\leq &
\frac{\left(%
\begin{array}{c}
  n+q \\
  \\
  n\\
\end{array}%
\right)}{\left\|p^{(\alpha,\beta)}_{i}\right\|_{2,(\alpha,\beta)}\left|P^{(\alpha,\beta)}_{n}(1)P^{(\alpha,\beta)}_{k-n}(1)\right|}\left\|p^{(\alpha,\beta)}_{i}\right\|_{2,(\alpha,\beta)}\left\|P^{(\alpha,\beta)}_{k-n}\right\|_{2,(\alpha,\beta)}\\
\nonumber&=&\frac{\left(%
\begin{array}{c}
  n+q \\
  \\
  n\\
\end{array}%
\right)\left\|P^{(\alpha,\beta)}_{k-n}\right\|_{2,(\alpha,\beta)}}{\left|P^{(\alpha,\beta)}_{n}(1)P^{(\alpha,\beta)}_{k-n}(1)\right|},
\end{eqnarray}
where $q=max(\alpha,\beta).$ Thus
\begin{equation*}
\left|P^{(\alpha,\beta)}_{n}(1)P^{(\alpha,\beta)}_{k-n}(1)\right|^{2}\left|\nu(i,k-n,n)\right|^{2}
\left\|p^{(\alpha,\beta)}_{i}\right\|^{2}_{2,(\alpha,\beta)} \leq
\left(%
\begin{array}{c}
  n+q \\
  \\
  n\\
\end{array}%
\right)^{2}\left\|P^{(\alpha,\beta)}_{k-n}\right\|_{2,
(\alpha,\beta)}^{2}.
\end{equation*}
Therefore, for  $\alpha=\beta=\lambda-1/2$, we have
\begin{eqnarray*}
&&\int_{-\infty}^{\infty}\left|H_{\lambda,K}^{N,1}(y)\right|^{2}e^{-y^{2}}\frac{dy}{\sqrt{\pi}}\\
&& \hspace{0.1cm}\leq
\frac{1}{\sqrt{\pi}}\frac{2\pi\Gamma(2\lambda)}{\lambda[\Gamma(\lambda)]^{2}2^{2\lambda}}\sum_{k=N+1}^{\infty}\sum_{n=1}^{k}\sum_{i=|k-2n|}^{k}\frac{\|\phi_{\lambda}\|^{2}}{\|C^{\lambda}_{k-n}\|^{2}\|C^{\lambda}_{n}\|^{2}}\left|\left\langle
\phi_{\lambda},\frac{C^{\lambda}_{n}}{\|C^{\lambda}_{n}\|}\right\rangle\right|^{2} t^{2}\hat{\lambda}_{n}\hat{\lambda}_{k-n}e^{-2t(\hat{\lambda}_{k}^{1/2}+\hat{\lambda}_{k-n}^{1/2})}\\
&&\hspace{0.2cm}\times
\left|\frac{[\Gamma(\lambda+1/2)]^{2}\Gamma(n+2\lambda)\Gamma(k-n+2\lambda)}{[\Gamma(2\lambda)]^{2}\Gamma(n+\lambda+1/2)\Gamma(k-n+\lambda+1/2)}\right|^{2}\left(\begin{array}{c}
  n+\lambda-1/2 \\
  \\
  n \\
\end{array}\right)^{2}\\
&&\hspace{1cm}\times\frac{[\Gamma(2\lambda)]^{2}[\Gamma(n+\lambda+1/2)]^{2}}{[\Gamma(\lambda+1/2)]^{2}[\Gamma(n+2\lambda)]^{2}}\|C^{\lambda}_{k-n}\|^{2}\\
&&\hspace{0.8cm}\times
\frac{[\Gamma(\lambda+1/2)]^{2}[\Gamma(k-n+2\lambda)]^{2}}{[\Gamma(2\lambda)]^{2}[\Gamma(k-n+\lambda+1/2)]^{2}}\frac{[\Gamma(n+\lambda+1/2)]^{2}}{[\Gamma(\lambda+1/2)]^{2}[n!]^{2}}\frac{\Gamma(2\lambda)(n+\lambda)n!}{\lambda\Gamma(n+2\lambda)}\\
&& \hspace{0.1cm}=\frac{1}{\sqrt{\pi}}\frac{2\pi\Gamma(2\lambda)}{\lambda[\Gamma(\lambda)]^{2}2^{2\lambda}}\sum_{k=N+1}^{\infty}\sum_{n=1}^{k}\sum_{i=|k-2n|}^{k}\|\phi_{\lambda}\|^{2}\left|\left\langle
\phi_{\lambda},\frac{C^{\lambda}_{n}}{\|C^{\lambda}_{n}\|}\right\rangle\right|^{2} t^{2}\hat{\lambda}_{n}\hat{\lambda}_{k-n}e^{-2t(\hat{\lambda}_{k}^{1/2}+\hat{\lambda}_{k-n}^{1/2})}\\
&&\hspace{0.8cm}\times
\frac{[\Gamma(k-n+2\lambda)]^{2}[\Gamma(n+\lambda+1/2)]^{2}(n+\lambda)}{\Gamma(2\lambda)[\Gamma(k-n+\lambda+1/2)]^{2}n!\lambda\Gamma(n+2\lambda)}.
\end{eqnarray*}
\\

On the other hand,  using the Stirling's approximation formula for
the gamma function,
$$\Gamma(az+b)\sim \sqrt{2\pi}e^{-az}(az)^{az+b-1/2}\;\;\;\;\; (|arg z|<\pi, a>0)$$
we obtain, for $\lambda $ big enough that
\begin{eqnarray*}
&&\frac{[\Gamma(k-n+2\lambda)]^{2}[\Gamma(n+\lambda+1/2)]^{2}(n+\lambda)}{\Gamma(2\lambda)[\Gamma(k-n+\lambda+1/2)]^{2}n!\lambda\Gamma(n+2\lambda)}\\
&& \hspace{0.5cm}\sim \frac{2\pi e^{-4\lambda}(2\lambda)^{4\lambda+2k-2n-1}2\pi
e^{-2\lambda}(\lambda)^{2\lambda+2n}(n+\lambda)}{\sqrt{2\pi}e^{-2\lambda}(2\lambda)^{2\lambda-1/2}2\pi
e^{-2\lambda}(\lambda)^{2\lambda+2k-2n}\sqrt{2\pi}
e^{-2\lambda}(2\lambda)^{2\lambda+n-1/2}\lambda
n!}=\frac{2^{2k}}{2^{3n}}\frac{\lambda^{n-1}(n+\lambda)}{n!}.
\end{eqnarray*}
Then,  for $\lambda $ big enough
\begin{eqnarray*}
&&\int_{-\infty}^{\infty}\left|H_{\lambda,K}^{N,1}(y)\right|^{2}e^{-y^{2}}\frac{dy}{\sqrt{\pi}}\\
& &\hspace{1cm}\leq
\frac{1}{\sqrt{\pi}}\frac{2\pi\Gamma(2\lambda)\|\phi_{\lambda}\|^{2}}{\lambda[\Gamma(\lambda)]^{2}2^{2\lambda}}\sum_{k=N+1}^{\infty}\sum_{n=1}^{k}\sum_{i=|k-2n|}^{k}\left|\left\langle
\phi_{\lambda},\frac{C^{\lambda}_{n}}{\|C^{\lambda}_{n}\|}\right\rangle\right|^{2} t^{2}\hat{\lambda}_{n}\hat{\lambda}_{k-n}\\
&&\hspace{5.5cm}\times e^{-2t(\hat{\lambda}_{k}^{1/2}+\hat{\lambda}_{k-n}^{1/2})}\frac{2^{2k}}{2^{3n}}\frac{\lambda^{n-1}(n+\lambda)}{n!}.\\
\end{eqnarray*}
Also,
\begin{eqnarray*}
t^{2}\hat{\lambda}_{n}\hat{\lambda}_{k-n}e^{-2t(\hat{\lambda}_{k}^{1/2}+\hat{\lambda}_{k-n}^{1/2})}\frac{2^{2k}}{2^{3n}}\frac{\lambda^{n-1}(n+\lambda)}{n!}&\leq &\frac{e^{-t}e^{-t\hat{\lambda}_{k-n}^{1/2}}2^{2k}\lambda^{n}(\frac{n}{\lambda}+1)}{t^{2}2^{3n}n!}.\\
\end{eqnarray*}
Taking $\lambda=(k-n)(1+2k)^{2}$, for $t\geq 1$ we get
\begin{eqnarray*}
\frac{1}{e^{t\hat{\lambda}_{k-n}^{1/2}}}\leq
\frac{1}{e^{(k-n)^{1/2}\left(k-n+2(k-n)\left(1+2k\right)^{2}\right)^{1/2}}}\leq
\frac{1}{e^{(k-n)\left(1+2\left(1+2k\right)^{2}\right)^{1/2}}}
\leq \frac{1}{e^{(k-n)}}\frac{1}{e^{k}}\frac{1}{e^{k}}.\\
\end{eqnarray*}
Now, if $0<t<1\;$ we have that for $t$ near to $0$, exists $k>0$
big enough such that $\frac{1}{k}\leq t$, i.e.
 $\;\frac{1}{t^{2}}\leq k^{2}$. Hence,
\begin{eqnarray*}
\frac{1}{e^{t\hat{\lambda}_{k-n}^{1/2}}}&\leq &
\frac{1}{e^{\frac{1}{k}(k-n)^{1/2}\left(k-n+2(k-n)\left(1+2k\right)^{2}\right)^{1/2}}}\leq \frac{1}{e^{(k-n)}}\frac{1}{e^{k}}\frac{1}{e^{k}}.\\
\end{eqnarray*}
Therefore, for $\lambda$ big enough there exists $C>0$ such that
\begin{eqnarray*}
t^{2}\hat{\lambda}_{n}\hat{\lambda}_{k-n}e^{-2t(\hat{\lambda}_{k}^{1/2}+\hat{\lambda}_{k-n}^{1/2})}\frac{2^{2k}}{2^{3n}}\frac{\lambda^{n-1}(n+\lambda)}{n!}\leq
\frac{e^{-tC}}{e^{k}}.
\end{eqnarray*}
Then, for $\lambda$ big enough we have
\begin{eqnarray*}
&&\int_{-\infty}^{\infty}\left|H_{\lambda,K}^{N,1}(y)\right|^{2}e^{-y^{2}}\frac{dy}{\sqrt{\pi}}\\
&& \hspace{0.5cm}\leq
\frac{1}{\sqrt{\pi}}\frac{2\pi\Gamma(2\lambda)}{\lambda[\Gamma(\lambda)]^{2}2^{2\lambda}}\sum_{k=N+1}^{\infty}\sum_{n=1}^{k}\sum_{i=|k-2n|}^{k}\left|\left\langle
\phi_{\lambda},\frac{C^{\lambda}_{n}}{\|C^{\lambda}_{n}\|}\right\rangle\right|^{2} t^{2}\hat{\lambda}_{n}\hat{\lambda}_{k-n}\\
&&\hspace{5.5cm}\times e^{-2t(\hat{\lambda}_{k}^{1/2}+\hat{\lambda}_{k-n}^{1/2})}\frac{2^{2k}}{2^{3n}}\frac{\lambda^{n-1}(n+\lambda)}{n!}\\
&& \hspace{0.5cm}\leq\frac{e^{-t}}{\sqrt{\pi}}\frac{2\pi\Gamma(2\lambda)}{\lambda[\Gamma(\lambda)]^{2}2^{2\lambda}}\frac{MC\|\phi_{\lambda}\|^{2}}{e^{(N+1)}}\sum_{k=N+1}^{\infty}\sum_{n=1}^{k}\sum_{i=|k-2n|}^{k}\left|\left\langle
\phi_{\lambda},\frac{C^{\lambda}_{n}}{\|C^{\lambda}_{n}\|}\right\rangle\right|^{2}
\\&& \hspace{0.5cm}\leq\frac{2\pi\Gamma(2\lambda)}{\lambda[\Gamma(\lambda)]^{2}2^{2\lambda}}\frac{e^{-t}MC\|\phi_{\lambda}\|^{4}_{2,\lambda}}{e^{(N+1)}}\leq \frac{e^{-t}MC\|\phi\|^{4}_{2,\gamma}}{e^{(N+1)}},\\
\end{eqnarray*}
Hence,
\begin{eqnarray*}
\left(\int_{-\infty}^{\infty}\left|\int_{0}^{\infty}H_{\lambda,K}^{N,1}(y)dt\right|^{2}e^{-y^{2}}\frac{dy}{\sqrt{\pi}}\right)^{1/2}&\leq
&\int_{0}^{\infty}\left(\frac{e^{-t}MC\|\phi\|^{4}_{2,\gamma}}{e^{(N+1)}}\right)^{1/2}dt\\
&=& (MC)^{1/2} \frac{\|\phi\|^{2}_{2,\gamma}}{e^{\frac{(N+1)}{2}}}\int_{0}^{\infty}e^{-\frac{t}{2}}dt= \frac{(MC)^{1/2}\|\phi\|^{2}_{2,\gamma}}{e^{\frac{(N+1)}{2}}}.\\
\end{eqnarray*}
Thus, $\left\{\int_{0}^{\infty}H_{\lambda,K}^{N,1}dt\right\}$ is a
bounded sequence on $\;L^{2}\left(\mathbb{R},\gamma\right),\;$ so
by Bourbaki-Alaoglu's theorem, there exists a sequence
$\;(\lambda_{j})_{j\in\mathbb{N}},\;$ $\lim_{j\to
\infty}\lambda_{j}=\infty$ such that, for all $N\in \mathbb{N}$,
  $\left\{\int_{0}^{\infty}H_{\lambda_{j},K}^{N,1}dt\right\}_{j\in\mathbb{N}}$ converges weakly in $\;L^{2}\left(\mathbb{R},\gamma\right)\;$ to a function
  $H_{K}^{N,1}\in L^{2}\left(\mathbb{R},\gamma\right)\;$ Moreover,
  \begin{equation}\label{norH_n}
\int_{-\infty}^{\infty}\left|H_{K}^{N,1}(y)\right|^{2}\frac{e^{-y^{2}}}{\sqrt{\pi}}dy\leq \frac{(MC)^{1/2}\|\phi\|^{2}_{2,\gamma}}{e^{\frac{(N+1)}{2}}}\\
\end{equation}
Then, there exists a non decreasing sequence
$(N_{j})_{j\in\mathbb{N}}$ such that,
\begin{equation}\label{limH_k}
H_{K}^{N_{j},1}(y)\longrightarrow 0,\hspace{0.2cm}
a.e.\hspace{0.2cm} y\in\mathbb{R}.
\end{equation}

The study  of $g_{2,\lambda}$ is essentially analogous, so fewer details will be given. For $\phi\in
L^{2}(\mathbb{R},\gamma)$ we have
\begin{eqnarray*}
\frac{\partial}{\partial
y}P^{(\alpha,\beta)}_{t}\phi_{\lambda}(\frac{y}{\sqrt{\lambda}})&=&-\sum_{k=1}^{\infty}\left\langle
\phi_{\lambda},\frac{P^{(\alpha,\beta)}_{k}}{\left\|P^{(\alpha,\beta)}_{k}\right\|_{2,(\alpha,\beta)}^{2}}\right\rangle
e^{-t\lambda_{k}^{1/2}}\frac{(k+\alpha+\beta+1)}{2}P^{(\alpha+1,\beta+1)}_{k-1}(\frac{y}{\sqrt{\lambda}}),
\end{eqnarray*}
which converges absolutely. Then, again taking the Cauchy product,
we obtain
\begin{eqnarray*}
&&\left(\frac{\partial}{\partial
y}P^{(\lambda-1/2,\lambda-1/2)}_{t}\phi_{\lambda}(\frac{y}{\sqrt{\lambda}})\right)^{2}\\
&&\hspace{0.2cm} =\sum_{k=1}^{\infty}\sum_{n=1}^{k}\left\langle
\phi_{\lambda},C^{\lambda}_{n}\right\rangle\left\langle
\phi_{\lambda},C^{\lambda}_{k-n}\right\rangle
e^{-t(\hat{\lambda}_{n}^{1/2}+\hat{\lambda}_{k-n}^{1/2})}C^{\lambda+1}_{n-1}(\frac{y}{\sqrt{\lambda}})C^{\lambda+1}_{k-n-1}(\frac{y}{\sqrt{\lambda}})\frac{(n+2\lambda)(k-n+2\lambda)}{4}\\
&&\hspace{0.3cm}\times\frac{\Gamma(2\lambda+2)\Gamma(n+\lambda+1/2)}{\Gamma(\lambda+3/2)\Gamma(n+2\lambda+1)}\frac{\Gamma(2\lambda+2)\Gamma(k-n+\lambda+1/2)}{\Gamma(\lambda+3/2)\Gamma(k-n+2\lambda+1)}
\frac{\Gamma(2\lambda)\Gamma(n+\lambda+1/2)}{\Gamma(\lambda+1/2)\Gamma(n+2\lambda)}\\
&&\hspace{0.3cm}\times\frac{\Gamma(2\lambda)\Gamma(k-n+\lambda+1/2)}{\Gamma(\lambda+1/2)\Gamma(k-n+2\lambda)}\frac{2(n+\lambda)[\Gamma(\lambda+1/2)]^{2}n!\Gamma(n+2\lambda)}{\Gamma(2\lambda+1)[\Gamma(n+\lambda+1/2)]^{2}}\\
&&\;\;\;\;\;\times\frac{2(k-n+\lambda)[\Gamma(\lambda+1/2)]^{2}(k-n)!\Gamma(k-n+2\lambda)}{\Gamma(2\lambda+1)[\Gamma(k-n+\lambda+1/2)]^{2}}\\
&&\hspace{0.2cm} =\sum_{k=1}^{\infty}\sum_{n=1}^{k}\left\langle
\phi_{\lambda},C^{\lambda}_{n}\right\rangle\left\langle
\phi_{\lambda},C^{\lambda}_{k-n}\right\rangle
e^{-t(\hat{\lambda}_{n}^{1/2}+\hat{\lambda}_{k-n}^{1/2})}C^{\lambda+1}_{n-1}(\frac{y}{\sqrt{\lambda}})C^{\lambda+1}_{k-n-1}(\frac{y}{\sqrt{\lambda}})\\
&&\;\hspace{2cm}\times\frac{4(n+\lambda)n!(k-n+\lambda)(k-n)!}
{(2\lambda)_{n}(2\lambda)_{k-n}},
\end{eqnarray*}
then
\begin{eqnarray*}
&&t(1-\frac{y^{2}}{\lambda})\left(\frac{\partial}{\partial
t}P^{(\lambda-1/2,\lambda-1/2)}_{t}\phi_{\lambda}(\frac{y}{\sqrt{\lambda}})\right)^{2}\\&=&(1-\frac{y^{2}}{\lambda})\sum_{k=1}^{\infty}\sum_{n=1}^{k}
\left\langle
\phi_{\lambda},\frac{C^{\lambda}_{n}}{\|C^{\lambda}_{n}\|^{2}}\right\rangle\left\langle
\phi_{\lambda},\frac{C^{\lambda}_{k-n}}{\|C^{\lambda}_{k-n}\|^{2}}\right\rangle
4\lambda^{2}t
e^{-t(\hat{\lambda}_{n}^{1/2}+\hat{\lambda}_{k-n}^{1/2})}
C^{\lambda+1}_{n-1}(\frac{y}{\sqrt{\lambda}})C^{\lambda+1}_{k-n-1}(\frac{y}{\sqrt{\lambda}}).
\end{eqnarray*}
Hence,
\begin{eqnarray*}
g_{2,\lambda}\phi(y)
&=&\chi_{[-K,K]}(y)\int_{0}^{\infty}\sum_{k=1}^{\infty}\sum_{n=1}^{k}\left\langle
\phi_{\lambda},\frac{C^{\lambda}_{n}}{\|C^{\lambda}_{n}\|^{2}}\right\rangle\left\langle
\phi_{\lambda},\frac{C^{\lambda}_{k-n}}{\|C^{\lambda}_{k-n}\|^{2}}\right\rangle 4\lambda^{2}t e^{-t(\hat{\lambda}_{n}^{1/2}+\hat{\lambda}_{k-n}^{1/2})}\\
&&\hspace{3cm} \times
C^{\lambda+1}_{n-1}(\frac{y}{\sqrt{\lambda}})C^{\lambda+1}_{k-n-1}(\frac{y}{\sqrt{\lambda}})(1-\frac{y^{2}}{\lambda})^{\lambda+1/2}e^{y^{2}}dt.
\end{eqnarray*}
Taking,
\begin{eqnarray*}
H_{\lambda,K}^{N,2}(y)&=&\chi_{[-K,K]}(y)\sum_{k=N+1}^{\infty}\sum_{n=1}^{k}\left\langle
\phi_{\lambda},\frac{C^{\lambda}_{n}}{\|C^{\lambda}_{n}\|^{2}}\right\rangle\left\langle
\phi_{\lambda},\frac{C^{\lambda}_{k-n}}{\|C^{\lambda}_{k-n}\|^{2}}\right\rangle 4\lambda^{2}t e^{-t(\hat{\lambda}_{n}^{1/2}+\hat{\lambda}_{k-n}^{1/2})}\\
&&\hspace{3cm} \times C^{\lambda+1}_{n-1}(\frac{y}{\sqrt{\lambda}})C^{\lambda+1}_{k-n-1}(\frac{y}{\sqrt{\lambda}})(1-\frac{y^{2}}{\lambda})^{\lambda/2+1/4}e^{\frac{y^{2}}{2}}
\end{eqnarray*}
we can write,
\begin{eqnarray*}
g_{2,\lambda}\phi(y)&=&\chi_{[-K,K]}(y)\int_{0}^{\infty}\sum_{k=1}^{N}\sum_{n=1}^{k}\left\langle
\phi_{\lambda},\frac{C^{\lambda}_{n}}{\|C^{\lambda}_{n}\|^{2}}\right\rangle\left\langle
\phi_{\lambda},\frac{C^{\lambda}_{k-n}}{\|C^{\lambda}_{k-n}\|^{2}}\right\rangle 4\lambda^{2} te^{-t(\hat{\lambda}_{n}^{1/2}+\hat{\lambda}_{k-n}^{1/2})}\\
&&\hspace{3cm} \times
C^{\lambda+1}_{n-1}(\frac{y}{\sqrt{\lambda}})C^{\lambda+1}_{k-n-1}(\frac{y}{\sqrt{\lambda}})(1-\frac{y^{2}}{\lambda})^{\lambda+1/2}e^{y^{2}}dt\\
&&\hspace{0.5cm}+\chi_{[-K,K]}(y)\int_{0}^{\infty}\sum_{k=N+1}^{\infty}\sum_{n=1}^{k}\left\langle
\phi_{\lambda},\frac{C^{\lambda}_{n}}{\|C^{\lambda}_{n}\|^{2}}\right\rangle\left\langle
\phi_{\lambda},\frac{C^{\lambda}_{k-n}}{\|C^{\lambda}_{k-n}\|^{2}}\right\rangle 4\lambda^{2} te^{-t(\hat{\lambda}_{n}^{1/2}+\hat{\lambda}_{k-n}^{1/2})}\\
&&\hspace{3cm} \times
C^{\lambda+1}_{n-1}(\frac{y}{\sqrt{\lambda}})C^{\lambda+1}_{k-n-1}(\frac{y}{\sqrt{\lambda}})(1-\frac{y^{2}}{\lambda})^{\lambda+1/2}e^{y^{2}}dt(1-\frac{y^{2}}{\lambda})^{\lambda+1/2}e^{y^{2}}dt\\
&=&\sum_{k=1}^{N}\sum_{n=1}^{k}\left\langle
\phi_{\lambda},\frac{C^{\lambda}_{n}}{\|C^{\lambda}_{n}\|^{2}}\right\rangle\left\langle
\phi_{\lambda},\frac{C^{\lambda}_{k-n}}{\|C^{\lambda}_{k-n}\|^{2}}\right\rangle \frac{4\lambda^{2}}{(\hat{\lambda}_{n}^{1/2}+\hat{\lambda}_{k-n}^{1/2})^{2}} \\
&&\hspace{3cm} \times
C^{\lambda+1}_{n-1}(\frac{y}{\sqrt{\lambda}})C^{\lambda+1}_{k-n-1}(\frac{y}{\sqrt{\lambda}})(1-\frac{y^{2}}{\lambda})^{\lambda+1/2}e^{y^{2}}\\
&&\hspace{0.5cm}+\int_{0}^{\infty} H_{\lambda,K}^{N,2}(y)dt(1-\frac{y^{2}}{\lambda})^{\lambda/2+1/4}e^{\frac{y^{2}}{2}}.\\
\end{eqnarray*}
Similarly the previous case, we want to prove that for any $K\in
\mathbb{N}$ and $\lambda>0$ such that $\sqrt{\lambda}>K,$
\begin{equation*}
\int_{-\infty}^{\infty}\left|\int_{0}^{\infty}H_{\lambda,K}^{N,2}(y)dt\right|^{2}e^{-y^{2}}\frac{dy}{\sqrt{\pi}}\leq
\frac{\theta}{e^{(N+1)}}.
\end{equation*}
Indeed, by Minkowski integral inequality we have
\begin{eqnarray*}
\left(\int_{-\infty}^{\infty}\left|\int_{0}^{\infty}H_{\lambda,K}^{N,2}(y)dt\right|^{2}e^{-y^{2}}\frac{dy}{\sqrt{\pi}}\right)^{1/2}&\leq
&\int_{0}^{\infty}\left(\int_{-\infty}^{\infty}\left|H_{\lambda,K}^{N,2}(y)\right|^{2}e^{-y^{2}}\frac{dy}{\sqrt{\pi}}\right)^{1/2}dt.
\end{eqnarray*}
Now,
\begin{eqnarray*}
&&\int_{-\infty}^{\infty}\left|H_{\lambda,K}^{N,2}(y)\right|^{2}e^{-y^{2}}\frac{dy}{\sqrt{\pi}}\\
&&\hspace{1cm} \leq
\frac{1}{\sqrt{\pi}}\int_{-1}^{1}\left|\sum_{k=N+1}^{\infty}\sum_{n=1}^{k}\left\langle
\phi_{\lambda},\frac{C^{\lambda}_{n}}{\|C^{\lambda}_{n}\|^{2}}\right\rangle\left\langle
\phi_{\lambda},\frac{C^{\lambda}_{k-n}}{\|C^{\lambda}_{k-n}\|^{2}}\right\rangle 4\lambda^{2}te^{-t(\hat{\lambda}_{k}^{1/2}+\hat{\lambda}_{k-n}^{1/2})}\right.\\
&&\;\hspace{4cm}\times \left.C^{\lambda+1}_{n-1}(x)C^{\lambda+1}_{k-n-1}(x)\right|^{2}(1-x^{2})^{\lambda-1/2}dx\\
\end{eqnarray*}
Again, using the Gasper's linearization of
the product of Jacobi polynomials,\cite{gasp1}, we may write
\begin{eqnarray*}
&&C^{\lambda+1}_{n-1}(x)C^{\lambda+1}_{k-n-1}(x)\\
&& \hspace{0.5cm} =\frac{[\Gamma(\lambda+3/2)]^{2}\Gamma(n+2\lambda+1)\Gamma(k-n+2\lambda+1)}{[\Gamma(2\lambda+2)]^{2}\Gamma(n+2\lambda+3/2)\Gamma(k-n+2\lambda+3/2)}
\\
&&\hspace{0.7cm} \times
P^{(\lambda+1/2,\lambda+1/2)}_{n-1}(1)P^{(\lambda+1/2,\lambda+1/2)}_{k-n-1}(1)p^{(\lambda+1/2,\lambda+1/2)}_{n-1}(x)p^{(\lambda+1/2,\lambda+1/2)}_{k-n-1}(x)\\
&&\hspace{0.5cm}=\frac{[\Gamma(\lambda+3/2)]^{2}\Gamma(n+2\lambda+1)\Gamma(k-n+2\lambda+1)}{[\Gamma(2\lambda+2)]^{2}\Gamma(n+2\lambda+3/2)\Gamma(k-n+2\lambda+3/2)}
\\
&&\hspace{0.7cm} \times
P^{(\lambda+1/2,\lambda+1/2)}_{n-1}(1)P^{(\lambda+1/2,\lambda+1/2)}_{k-n-1}(1)\sum_{i=|k-2n|}^{k-2}\nu(i,n-1,k-n-1)p^{(\lambda+1/2,\lambda+1/2)}_{i}(x).
\end{eqnarray*}
Thus,
\begin{eqnarray*}
&&\int_{-\infty}^{\infty}\left|H_{\lambda,K}^{N,2}(y)\right|^{2}e^{-y^{2}}\frac{dy}{\sqrt{\pi}}\\
&&\hspace{0.5cm}\leq
\frac{1}{\sqrt{\pi}}\int_{-1}^{1}\left|\sum_{k=N+1}^{\infty}\sum_{n=1}^{k}\sum_{i=|k-2n|}^{k-2}\left\langle
\phi_{\lambda},\frac{C^{\lambda}_{n}}{\|C^{\lambda}_{n}\|^{2}}\right\rangle\left\langle
\phi_{\lambda},\frac{C^{\lambda}_{k-n}}{\|C^{\lambda}_{k-n}\|^{2}}\right\rangle 4\lambda^{2}te^{-t(\hat{\lambda}_{n}^{1/2}+\hat{\lambda}_{k-n}^{1/2})}\right.\\
&&\hspace{0.7cm}\times \frac{[\Gamma(\lambda+3/2)]^{2}\Gamma(n+2\lambda+1)\Gamma(k-n+2\lambda+1)}{[\Gamma(2\lambda+2)]^{2}\Gamma(n+2\lambda+3/2)\Gamma(k-n+2\lambda+3/2)}P^{(\lambda+1/2,\lambda+1/2)}_{n-1}(1)P^{(\lambda+1/2,\lambda+1/2)}_{k-n-1}(1)\\
&&\hspace{0.8cm}\times
\left.\nu(i,n-1,k-n-1)p^{(\lambda+1/2,\lambda+1/2)}_{i}(x)\right|^{2}(1-x^{2})^{\lambda-1/2}dx.\\
\end{eqnarray*}
Then, using Parseval's identity and Cauchy-Schwartz inequality, we
have
\begin{eqnarray*}
&&\int_{-\infty}^{\infty}\left|H_{\lambda,K}^{N,2}(y)\right|^{2}e^{-y^{2}}\frac{dy}{\sqrt{\pi}}\\
&&\hspace{0.5cm}\leq
\frac{1}{\sqrt{\pi}}\frac{2\pi\Gamma(2(\lambda+1))}{(\lambda+1)[\Gamma(\lambda+1)]^{2}2^{2(\lambda+1)}}\sum_{k=N+1}^{\infty}\sum_{n=1}^{k}\sum_{i=|k-2n|}^{k}\left|\left\langle
\phi_{\lambda},\frac{C^{\lambda}_{n}}{\|C^{\lambda}_{n}\|}\right\rangle\right|^{2}\frac{\|
\phi_{\lambda}\|^{2}}{\|C^{\lambda}_{n}\|^{2}\|C^{\lambda}_{k-n}\|^{2}}\\
&&\hspace{0.6cm}\times 16\lambda^{4}t^{2}e^{-2t(\hat{\lambda}_{k}^{1/2}+\hat{\lambda}_{k-n}^{1/2})}\left|\frac{[\Gamma(\lambda+3/2)]^{2}\Gamma(n+2\lambda+1)\Gamma(k-n+2\lambda+1)}{[\Gamma(2\lambda+2)]^{2}\Gamma(n+2\lambda+3/2)\Gamma(k-n+2\lambda+3/2)}\right|^{2}\\
&&\;\hspace{0.7cm}\times \left|P^{(\lambda+1/2,\lambda+1/2)}_{n-1}(1)P^{(\lambda+1/2,\lambda+1/2)}_{k-n-1}(1)\right|^{2}\left|\nu(i,n-1,k-n-1)\right|^{2}\left\|p^{(\lambda+1/2,\lambda+1/2)}_{i}\right\|^{2}.\\
\end{eqnarray*}
Then, by by a similar argument as in the previous case, we get
\begin{eqnarray*}
&&\int_{-\infty}^{\infty}\left|H_{\lambda,K}^{N,2}(y)\right|^{2}e^{-y^{2}}\frac{dy}{\sqrt{\pi}}\\
&\leq
&\frac{1}{\sqrt{\pi}}\frac{2\pi\Gamma(2(\lambda+1))}{(\lambda+1)[\Gamma(\lambda+1)]^{2}2^{2(\lambda+1)}}\sum_{k=N+1}^{\infty}\sum_{n=1}^{k}\sum_{i=|k-2n|}^{k}\left|\left\langle
\phi_{\lambda},\frac{C^{\lambda}_{n}}{\|C^{\lambda}_{n}\|}\right\rangle\right|^{2}\frac{\|
\phi_{\lambda}\|^{2}}{\|C^{\lambda}_{n}\|^{2}\|C^{\lambda}_{k-n}\|^{2}}\\
&&\hspace{0.2cm}\times 16\lambda^{4}t^{2}e^{-2t(\hat{\lambda}_{k}^{1/2}+\hat{\lambda}_{k-n}^{1/2})}\left|\frac{[\Gamma(\lambda+3/2)]^{2}\Gamma(n+2\lambda+1)\Gamma(k-n+2\lambda+1)}{[\Gamma(2\lambda+2)]^{2}\Gamma(n+2\lambda+3/2)\Gamma(k-n+2\lambda+3/2)}\right|^{2}\\
&&\;\hspace{0.2cm}\times \left(%
\begin{array}{c}
  n+\lambda-1/2 \\
  \\
  n-1\\
\end{array}%
\right)^{2}\left[\frac{\Gamma(2\lambda)\Gamma(k-n+\lambda+1/2)}{\Gamma(\lambda+1/2)\Gamma(k-n+2\lambda)}\right]^{2}\frac{(2\lambda+2)(2\lambda+1)(k-n)}{(\lambda+1/2)^{2}(k-n+2\lambda)}\left\|C^{\lambda}_{k-n}\right\|_{\lambda}^{2}.\\
\end{eqnarray*}
On the other hand,
\begin{eqnarray*}
&&\left|\frac{[\Gamma(\lambda+3/2)]^{2}\Gamma(n+2\lambda+1)\Gamma(k-n+2\lambda+1)}{[\Gamma(2\lambda+2)]^{2}\Gamma(n+2\lambda+3/2)\Gamma(k-n+2\lambda+3/2)}\right|^{2}
\\&&\hspace{0.2cm}\times \left(%
\begin{array}{c}
 n+\lambda-1/2\\
  \\
  n-1\\
\end{array}%
\right)^{2}\left[\frac{\Gamma(2\lambda)\Gamma(k-n+\lambda+1/2)}{\Gamma(\lambda+1/2)\Gamma(k-n+2\lambda)}\right]^{2}\frac{(2\lambda+2)(2\lambda+1)(k-n)}{(\lambda+1/2)^{2}(k-n+2\lambda)}\left\|C^{\lambda}_{k-n}\right\|_{\lambda}^{2}\\
&&\hspace{0.8cm} = \frac{2[\Gamma(n+2\lambda)]^{2}
\left[\Gamma(n+\lambda+1/2)\right]^{2}\left[\Gamma(k-n+\lambda+1/2)\right]^{2}\left\|C^{\lambda}_{k-n}\right\|_{\lambda}^{2}}{(2\lambda)^{4}[\Gamma(2\lambda)]^{2}[\Gamma(n+2\lambda+1/2)]^{2}[\Gamma(k-n+2\lambda+1/2)]^{2}[\Gamma(n)]^{2}}.
\end{eqnarray*}
Hence,
\begin{eqnarray*}
&&\int_{-\infty}^{\infty}\left|H_{\lambda,K}^{N,2}(y)\right|^{2}e^{-y^{2}}\frac{dy}{\sqrt{\pi}}\\
&&\hspace{0.5cm}\leq
\frac{1}{\sqrt{\pi}}\frac{2\pi\Gamma(2(\lambda+1))}{(\lambda+1)[\Gamma(\lambda+1)]^{2}2^{2(\lambda+1)}}\sum_{k=N+1}^{\infty}\sum_{n=1}^{k}\sum_{i=|k-2n|}^{k}\left|\left\langle
\phi_{\lambda},\frac{C^{\lambda}_{n}}{\|C^{\lambda}_{n}\|}\right\rangle\right|^{2}\frac{\|
\phi_{\lambda}\|^{2}t^{2}e^{-2t(\hat{\lambda}_{k}^{1/2}+\hat{\lambda}_{k-n}^{1/2})}}{\|C^{\lambda}_{n}\|^{2}}\\
&&\hspace{0.7cm}\times \frac{2[\Gamma(n+2\lambda)]^{2}
\left[\Gamma(n+\lambda+1/2)\right]^{2}\left[\Gamma(k-n+\lambda+1/2)\right]^{2}\left\|C^{\lambda}_{k-n}\right\|_{\lambda}^{2}}{[\Gamma(2\lambda)]^{2}[\Gamma(n+2\lambda+1/2)]^{2}[\Gamma(k-n+2\lambda+1/2)]^{2}[\Gamma(n)]^{2}}\\
&&\hspace{0.5cm}=
\frac{1}{\sqrt{\pi}}\frac{2\pi\Gamma(2(\lambda+1))\|\phi_{\lambda}\|^{2}}{(\lambda+1)[\Gamma(\lambda+1)]^{2}2^{2(\lambda+1)}}\sum_{k=N+1}^{\infty}\sum_{n=1}^{k}\sum_{i=|k-2n|}^{k}\left|\left\langle
\phi_{\lambda},\frac{C^{\lambda}_{n}}{\|C^{\lambda}_{n}\|}\right\rangle\right|^{2}t^{2}e^{-2t(\hat{\lambda}_{k}^{1/2}+\hat{\lambda}_{k-n}^{1/2})}\\
&&\hspace{0.7cm}\times \frac{2[\Gamma(n+2\lambda)]
\left[\Gamma(n+\lambda+1/2)\right]^{2}\left[\Gamma(k-n+\lambda+1/2)\right]^{2}}{\Gamma(2\lambda)[\Gamma(n+2\lambda+1/2)]^{2}[\Gamma(k-n+2\lambda+1/2)]^{2}[\Gamma(n)]^{2}}\frac{(n+\lambda)n!}{\lambda}.
\end{eqnarray*}
Then, for $\lambda$ big enough we have
\begin{eqnarray*}
&&\frac{2[\Gamma(n+2\lambda)]
\left[\Gamma(n+\lambda+1/2)\right]^{2}\left[\Gamma(k-n+\lambda+1/2)\right]^{2}}{\Gamma(2\lambda)[\Gamma(n+2\lambda+1/2)]^{2}[\Gamma(k-n+2\lambda+1/2)]^{2}[\Gamma(n)]^{2}}\frac{(n+\lambda)n!}{\lambda}\\
&& \hspace{0.2cm}\sim  \frac{2
\sqrt{2\pi}e^{-2\lambda}(2\lambda)^{2\lambda+n-1/2}2\pi
e^{-2\lambda}(\lambda)^{2\lambda+2n}2\pi
e^{-2\lambda}(\lambda)^{2\lambda+2k-2n}(n+\lambda)n!}{\sqrt{2\pi}e^{-2\lambda}(2\lambda)^{2\lambda-1/2}2\pi
e^{-4\lambda}(2\lambda)^{4\lambda+2n}2\pi
e^{-4\lambda}(2\lambda)^{4\lambda+2k-2n}\lambda
[(n-1)!]^{2}}\\
&&\hspace{0.2cm} =\frac{2^{n+1}\lambda^{n-1}(n+\lambda)n!}{2^{2k+8\lambda}\lambda^{4\lambda}e^{-4\lambda}[(n-1)!]^{2}}\leq \left(\frac{n}{\lambda}+1\right)\frac{\lambda^{n}n}{\lambda^{4\lambda}(n-1)!}.\\
\end{eqnarray*}
Therefore,  for $\lambda$ big enough we have
\begin{eqnarray*}
&&\int_{-\infty}^{\infty}\left|H_{\lambda,K}^{N,2}(y)\right|^{2}e^{-y^{2}}\frac{dy}{\sqrt{\pi}}
\leq
\frac{\|\phi_{\lambda}\|^{2}}{\sqrt{\pi}}\frac{2\pi\Gamma(2(\lambda+1))}{(\lambda+1)[\Gamma(\lambda+1)]^{2}2^{2(\lambda+1)}}\\
&&\hspace{2.5cm}\times\sum_{k=N+1}^{\infty}\sum_{n=1}^{k}\sum_{i=|k-2n|}^{k}\left|\left\langle
\phi_{\lambda},\frac{C^{\lambda}_{n}}{\|C^{\lambda}_{n}\|}\right\rangle\right|^{2}
t^{2}e^{-2t(\hat{\lambda}_{k}^{1/2}+\hat{\lambda}_{k-n}^{1/2})}\left(\frac{n}{\lambda}+1\right)\frac{\lambda^{n}n}{\lambda^{4\lambda}(n-1)!}.\\
\end{eqnarray*}
Similarly the previous case, taking $\lambda=(k-n)(1+2k)^{2}$ we
have for $k$ big enough there exists $C>0$ such that
\begin{eqnarray*}
t^{2}e^{-2t(\hat{\lambda}_{k}^{1/2}+\hat{\lambda}_{k-n}^{1/2})}\left(\frac{n}{\lambda}+1\right)\frac{\lambda^{n}n}{\lambda^{4\lambda}(n-1)!}\leq
\frac{e^{-t}C}{e^{k}}.
\end{eqnarray*}
Thus, for $\lambda$ big enough we have
\begin{eqnarray*}
&&\int_{-\infty}^{\infty}\left|H_{\lambda,K}^{N,2}(y)\right|^{2}e^{-y^{2}}\frac{dy}{\sqrt{\pi}}
\leq
\frac{\|\phi_{\lambda}\|^{2}}{\sqrt{\pi}}\frac{2\pi\Gamma(2(\lambda+1))}{(\lambda+1)[\Gamma(\lambda+1)]^{2}2^{2(\lambda+1)}}\\
&&\times\sum_{k=N+1}^{\infty}\sum_{n=1}^{k}\sum_{i=|k-2n|}^{k}\left|\left\langle
\phi_{\lambda},\frac{C^{\lambda}_{n}}{\|C^{\lambda}_{n}\|}\right\rangle\right|^{2}
t^{2}e^{-2t(\hat{\lambda}_{k}^{1/2}+\hat{\lambda}_{k-n}^{1/2})}\left(\frac{n}{\lambda}+1\right)\frac{\lambda^{n}n}{\lambda^{4\lambda}(n-1)!}\\
&&\hspace{0.5cm} \leq
\frac{e^{-t}}{\sqrt{\pi}}\frac{2\pi(2\lambda+1)2\lambda\Gamma(2\lambda)}{(\lambda+1)\lambda^{2}[\Gamma(\lambda)]^{2}2^{2\lambda}2^{2}}\frac{MC\|\phi_{\lambda}\|^{2}}{e^{(N+1)}}\sum_{k=N+1}^{\infty}\sum_{n=1}^{k}\sum_{i=|k-2n|}^{k}\left|\left\langle
\phi_{\lambda},\frac{C^{\lambda}_{n}}{\|C^{\lambda}_{n}\|}\right\rangle\right|^{2}\\
&&\hspace{0.5cm}\leq
\frac{2\pi\Gamma(2\lambda)}{\lambda[\Gamma(\lambda)]^{2}2^{2\lambda}}\frac{e^{-t}MC\|\phi_{\lambda}\|^{4}_{2,\lambda}}{e^{(N+1)}}\leq \frac{e^{-t}MC\|\phi\|^{4}_{2,\gamma}}{e^{(N+1)}},\\
\end{eqnarray*}
Hence
\begin{eqnarray*}
\left(\int_{-\infty}^{\infty}\left|\int_{0}^{\infty}H_{\lambda,K}^{N,2}(y)dt\right|^{2}e^{-y^{2}}\frac{dy}{\sqrt{\pi}}\right)^{1/2}&\leq
&\int_{0}^{\infty}\left(\frac{e^{-t}MC\|\phi\|^{4}_{2,\gamma}}{e^{(N+1)}}\right)^{1/2}dt\\
&=& \frac{(MC)^{1/2}\|\phi\|^{2}_{2,\gamma}}{e^{\frac{(N+1)}{2}}}.\\
\end{eqnarray*}
 Then,
$\left\{\int_{0}^{\infty}H_{\lambda,K}^{N,2}dt\right\}$ is a
bounded sequence on $\;L^{2}\left(\mathbb{R},\gamma\right),\;$ so
by Bourbaki-Alaoglu's theorem, there exists a sequence
$\;(\lambda_{j})_{j\in\mathbb{N}}\;$ $\lim_{j\to
\infty}\lambda_{j}=\infty$ such that, for all $N\in \mathbb{N}$,
  $\left\{\int_{0}^{\infty}H_{\lambda_{j},K}^{N,2}dt\right\}_{j\in\mathbb{N}}$ converges weakly in $\;L^{2}\left(\mathbb{R},\gamma\right)\;$ to a function
  $H_{K}^{N,2}\in L^{2}\left(\mathbb{R},\gamma\right)\;$ Moreover,
   \begin{equation}\label{norH_n2}
\int_{-\infty}^{\infty}\left|H_{K}^{N,2}(y)\right|^{2}\frac{e^{-y^{2}}}{\sqrt{\pi}}dy\leq \frac{(MC)^{1/2}\|\phi\|^{2}_{2,\gamma}}{e^{\frac{(N+1)}{2}}}\\
\end{equation}
Therefore, there exists a non decreasing sequence
$(N_{j})_{j\in\mathbb{N}}$ such that,
\begin{equation}\label{limH_k}
H_{K}^{N_{j},2}(y)\longrightarrow 0,\hspace{0.2cm}
a.e.\hspace{0.2cm} y\in\mathbb{R}.
\end{equation}
On the other hand,
\begin{eqnarray}\label{F-lamdk}
\nonumber\left(F_{\lambda
,K}(y)\right)^{2}&=&g_{1,\lambda}\phi(y)+g_{2,\lambda}\phi(y)\\
\nonumber &=&\sum_{k=1}^{N}\sum_{n=1}^{k}\left\langle
\phi_{\lambda},\frac{C^{\lambda}_{n}}{\|C^{\lambda}_{n}\|^{2}}\right\rangle\left\langle
\phi_{\lambda},\frac{C^{\lambda}_{k-n}}{\|C^{\lambda}_{k-n}\|^{2}}\right\rangle
\frac{\hat{\lambda}_{n}^{1/2}\hat{\lambda}_{k-n}^{1/2}}{(\hat{\lambda}_{n}^{1/2}+\hat{\lambda}_{k-n}^{1/2})^{2}} C^{\lambda}_{n}(\frac{y}{\sqrt{\lambda}})C^{\lambda}_{k-n}(\frac{y}{\sqrt{\lambda}})\\
\nonumber &&\;\hspace{0.3cm}\times(1-\frac{y^{2}}{\lambda})^{\lambda-1/2}e^{y^{2}}\;\;+\;\;\int_{0}^{\infty}H_{\lambda,K}^{N,1}(y)dt(1-\frac{y^{2}}{\lambda})^{\lambda/2-1/4}e^{\frac{y^{2}}{2}}\\
\nonumber &&\hspace{1cm}+ \sum_{k=1}^{N}\sum_{n=1}^{k}\left\langle
\phi_{\lambda},\frac{C^{\lambda}_{n}}{\|C^{\lambda}_{n}\|^{2}}\right\rangle\left\langle
\phi_{\lambda},\frac{C^{\lambda}_{k-n}}{\|C^{\lambda}_{k-n}\|^{2}}\right\rangle \frac{4\lambda^{2}}{(\hat{\lambda}_{n}^{1/2}+\hat{\lambda}_{k-n}^{1/2})^{2}} \\
\nonumber &&\hspace{3cm} \times
C^{\lambda+1}_{n-1}(\frac{y}{\sqrt{\lambda}})C^{\lambda+1}_{k-n-1}(\frac{y}{\sqrt{\lambda}})(1-\frac{y^{2}}{\lambda})^{\lambda+1/2}e^{y^{2}}\\
\nonumber &&\hspace{0.5cm}+\int_{0}^{\infty} H_{\lambda,K}^{N,2}(y)dt(1-\frac{y^{2}}{\lambda})^{\lambda/2+1/4}e^{\frac{y^{2}}{2}}\\
\nonumber
&=&F^{N}_{\lambda,K}(y)\;+\;\int_{0}^{\infty}H_{\lambda,K}^{N,1}(y)dt(1-\frac{y^{2}}{\lambda})^{\lambda/2-1/4}e^{\frac{y^{2}}{2}}\\&&\hspace{3cm}+\;\;\int_{0}^{\infty}
H_{\lambda,K}^{N,2}(y)dt(1-\frac{y^{2}}{\lambda})^{\lambda/2+1/4}e^{\frac{y^{2}}{2}}
\end{eqnarray}
where,
\begin{eqnarray*}
F^{N}_{\lambda,K}(y)&=&\sum_{k=1}^{N}\sum_{n=1}^{k}\left\langle
\phi_{\lambda},\frac{C^{\lambda}_{n}}{\|C^{\lambda}_{n}\|^{2}}\right\rangle\left\langle
\phi_{\lambda},\frac{C^{\lambda}_{k-n}}{\|C^{\lambda}_{k-n}\|^{2}}\right\rangle
\frac{\hat{\lambda}_{n}^{1/2}\hat{\lambda}_{k-n}^{1/2}}{(\hat{\lambda}_{n}^{1/2}+\hat{\lambda}_{k-n}^{1/2})^{2}} C^{\lambda}_{n}(\frac{y}{\sqrt{\lambda}})C^{\lambda}_{k-n}(\frac{y}{\sqrt{\lambda}})\\
&&\;\hspace{0.3cm}\times(1-\frac{y^{2}}{\lambda})^{\lambda-1/2}e^{y^{2}}\;\;+\;\;\sum_{k=1}^{N}\sum_{n=1}^{k}\left\langle
\phi_{\lambda},\frac{C^{\lambda}_{n}}{\|C^{\lambda}_{n}\|^{2}}\right\rangle\left\langle
\phi_{\lambda},\frac{C^{\lambda}_{k-n}}{\|C^{\lambda}_{k-n}\|^{2}}\right\rangle \frac{4\lambda^{2}}{(\hat{\lambda}_{n}^{1/2}+\hat{\lambda}_{k-n}^{1/2})^{2}} \\
&&\hspace{0.3cm} \times
C^{\lambda+1}_{n-1}(\frac{y}{\sqrt{\lambda}})C^{\lambda+1}_{k-n-1}(\frac{y}{\sqrt{\lambda}})(1-\frac{y^{2}}{\lambda})^{\lambda+1/2}e^{y^{2}}
\end{eqnarray*}
Defining, for each $m\in\mathbb{N}$,
$F^{N_{m}}_{K}=F^{2}-H^{N_{m},1}_{K}-H^{N_{m},2}_{K},$ then since
$ F_{\lambda,K}\to  F_{K},$ as $j\to\infty$ in the weak topology
of ${L^{2}\left(\mathbb{R},\gamma\right)}$ and that  $$F_{\lambda
,K}^{N}=\left(F_{\lambda
,K}\right)^{2}-\int_{0}^{\infty}H_{\lambda,K}^{N,1}(y)dt(1-\frac{y^{2}}{\lambda})^{\lambda/2-1/4}e^{\frac{y^{2}}{2}}\;\;-\;\;\int_{0}^{\infty}
H_{\lambda,K}^{N,2}(y)dt(1-\frac{y^{2}}{\lambda})^{\lambda/2+1/4}e^{\frac{y^{2}}{2}},\;$$
we have, for all $m\in\mathbb{N},$
\begin{equation*}
F^{N_{m}}_{\lambda_{j},K}\longrightarrow
F^{N_{m}}_{K},\hspace{1cm}\mbox{weakly in}
\;L^{2}\left(\mathbb{R},\gamma\right),\hspace{0.3cm}\mbox{as}\hspace{0.3cm}j\to
\infty
\end{equation*}
also,
\begin{equation}\label{F-Herm}
F^{N_{m}}_{K}(y)\longrightarrow
F^{2}(y),\hspace{0.3cm}\mbox{as}\hspace{0.3cm}m\to
\infty\hspace{0.3cm}a.e.\hspace{0.3cm}y\in\mathbb{R}
\end{equation}

Now,
\begin{eqnarray*}
&&F^{N}_{\lambda,K}(y)\\&&\hspace{0.2cm} =\sum_{k=1}^{N}\sum_{n=1}^{k}\left\langle
\phi_{\lambda},\frac{C^{\lambda}_{n}}{\|C^{\lambda}_{n}\|^{2}}\right\rangle\left\langle
\phi_{\lambda},\frac{C^{\lambda}_{k-n}}{\|C^{\lambda}_{k-n}\|^{2}}\right\rangle
\frac{\hat{\lambda}_{n}^{1/2}\hat{\lambda}_{k-n}^{1/2}}{(\hat{\lambda}_{n}^{1/2}+\hat{\lambda}_{k-n}^{1/2})^{2}} C^{\lambda}_{n}(\frac{y}{\sqrt{\lambda}})C^{\lambda}_{k-n}(\frac{y}{\sqrt{\lambda}})\\
&&\;\hspace{0.4cm}\times(1-\frac{y^{2}}{\lambda})^{\lambda-1/2}e^{y^{2}}\;\;+\;\;\sum_{k=1}^{N}\sum_{n=1}^{k}\left\langle
\phi_{\lambda},\frac{C^{\lambda}_{n}}{\|C^{\lambda}_{n}\|^{2}}\right\rangle\left\langle
\phi_{\lambda},\frac{C^{\lambda}_{k-n}}{\|C^{\lambda}_{k-n}\|^{2}}\right\rangle \frac{4\lambda^{2}}{(\hat{\lambda}_{n}^{1/2}+\hat{\lambda}_{k-n}^{1/2})^{2}} \\
&&\hspace{0.35cm} \times
C^{\lambda+1}_{n-1}(\frac{y}{\sqrt{\lambda}})C^{\lambda+1}_{k-n-1}(\frac{y}{\sqrt{\lambda}})(1-\frac{y^{2}}{\lambda})^{\lambda+1/2}e^{y^{2}}\\
&&\hspace{0.2cm} =\sum_{k=1}^{N}\sum_{n=1}^{k}\frac{(n+\lambda)n!(k-n)!(k-n+\lambda)}{\lambda^{2}(2\lambda)_{n}(2\lambda)_{k-n}}
\frac{\hat{\lambda}_{n}^{1/2}\hat{\lambda}_{k-n}^{1/2}}{(\hat{\lambda}_{n}^{1/2}+\hat{\lambda}_{k-n}^{1/2})^{2}}\left\langle\phi_{\lambda},C^{\lambda}_{n}\right\rangle\left\langle
\phi_{\lambda},C^{\lambda}_{k-n}\right\rangle \\
&&\hspace{0.3cm} \times
C^{\lambda}_{n}(\frac{y}{\sqrt{\lambda}})C^{\lambda}_{k-n}(\frac{y}{\sqrt{\lambda}})(1-\frac{y^{2}}{\lambda})^{\lambda-1/2}e^{y^{2}}\\&&+\;\;\sum_{k=1}^{N}\sum_{n=1}^{k}\left\langle\phi_{\lambda},C^{\lambda}_{n}\right\rangle\left\langle\phi_{\lambda},C^{\lambda}_{k-n}\right\rangle
C^{\lambda+1}_{n-1}(\frac{y}{\sqrt{\lambda}})C^{\lambda+1}_{k-n-1}(\frac{y}{\sqrt{\lambda}})\\
&&\hspace{0.3cm}\times\frac{4(n+\lambda)n!(k-n+\lambda)(k-n)!}
{(2\lambda)_{n}(2\lambda)_{k-n}}
\frac{1}{(\hat{\lambda}_{k}^{1/2}+\hat{\lambda}_{k-n}^{1/2})^{2}}(1-\frac{y^{2}}{\lambda})^{\lambda+1/2}e^{y^{2}}\\
\end{eqnarray*}
\begin{eqnarray*}
&=&\sum_{k=1}^{N}\sum_{n=1}^{k}\frac{(\frac{n}{\lambda}+1)n!(k-n)!(\frac{k-n}{\lambda}+1)}{2(2+\frac{1}{\lambda})\ldots(2+\frac{n-1}{\lambda})2(2+\frac{1}{\lambda})\ldots(2+\frac{k-n-1}{\lambda})}\\
&&\;\;\hspace{0.5cm}\times\frac{(n(\frac{n}{\lambda}+2))^{1/2}((k-n)(\frac{k-n}{\lambda}+2))^{1/2}}{\{(n(\frac{n}{\lambda}+2))^{1/2}+((k-n)(\frac{k-n}{\lambda}+2))^{1/2}\}^{2}}\left\langle\phi_{\lambda},\lambda^{-\frac{n}{2}}C^{\lambda}_{n}\right\rangle\left\langle
\phi_{\lambda},\lambda^{-\frac{k-n}{2}}C^{\lambda}_{k-n}\right\rangle \\
&&\;\;\hspace{0.5cm}\times\lambda^{-\frac{n}{2}}C^{\lambda}_{n}(\frac{y}{\sqrt{\lambda}})\lambda^{-\frac{k-n}{2}}C^{\lambda}_{k-n}(\frac{y}{\sqrt{\lambda}})(1-\frac{y^{2}}{\lambda})^{\lambda-1/2}e^{y^{2}}\\\
&&\;\;+\sum_{k=1}^{N}\sum_{n=1}^{k}\left\langle
\phi_{\lambda},\lambda^{-\frac{n}{2}}C^{\lambda}_{n}\right\rangle\left\langle
\phi_{\lambda},\lambda^{-\frac{k-n}{2}}C^{\lambda}_{k-n}\right\rangle
\lambda^{-\frac{n-1}{2}}C^{\lambda+1}_{n-1}(\frac{y}{\sqrt{\lambda}})\lambda^{-\frac{k-n-1}{2}}C^{\lambda+1}_{k-n-1}(\frac{y}{\sqrt{\lambda}})\\
&&\;\hspace{0.5cm}\times\frac{4(\frac{n}{\lambda}+1)n!(\frac{k-n}{\lambda}+1)(k-n)!}
{2(2+\frac{1}{\lambda})\ldots(2+\frac{n-1}{\lambda})2(2+\frac{1}{\lambda})\ldots(2+\frac{k-n-1}{\lambda})}\\
&&\;\hspace{0.5cm}\times\frac{1}{\{(n(\frac{n}{\lambda}+2))^{1/2}+((k-n)(\frac{k-n}{\lambda}+2))^{1/2}\}^{2}}(1-\frac{y^{2}}{\lambda})^{\lambda+1/2}e^{y^{2}}.\\
\end{eqnarray*}
Then,
\begin{eqnarray*}
 &&\lim_{\lambda\to\infty }F^{N}_{\lambda,K}(y)\\
 &&\hspace{0.2cm} =\lim_{\lambda\to\infty }\sum_{k=1}^{N}\sum_{n=1}^{k}\frac{(\frac{n}{\lambda}+1)n!(k-n)!(\frac{k-n}{\lambda}+1)}{2(2+\frac{1}{\lambda})\ldots(2+\frac{n-1}{\lambda})2(2+\frac{1}{\lambda})\ldots(2+\frac{k-n-1}{\lambda})}\\
&&\;\;\hspace{0.5cm}\times\frac{(n(\frac{n}{\lambda}+2))^{1/2}((k-n)(\frac{k-n}{\lambda}+2))^{1/2}}{\{(n(\frac{n}{\lambda}+2))^{1/2}+((k-n)(\frac{k-n}{\lambda}+2))^{1/2}\}^{2}}\left\langle\phi_{\lambda},\lambda^{-\frac{n}{2}}C^{\lambda}_{n}\right\rangle\left\langle
\phi_{\lambda},\lambda^{-\frac{k-n}{2}}C^{\lambda}_{k-n}\right\rangle \\
&&\hspace{0.5cm}\times\lambda^{-\frac{n}{2}}C^{\lambda}_{n}(\frac{y}{\sqrt{\lambda}})\lambda^{-\frac{k-n}{2}}C^{\lambda}_{k-n}(\frac{y}{\sqrt{\lambda}})(1-\frac{y^{2}}{\lambda})^{\lambda-1/2}e^{y^{2}}\\\
&&\hspace{0.5cm}+\lim_{\lambda\to\infty
}\sum_{k=1}^{N}\sum_{n=1}^{k}\left\langle
\phi_{\lambda},\lambda^{-\frac{n}{2}}C^{\lambda}_{n}\right\rangle\left\langle
\phi_{\lambda},\lambda^{-\frac{k-n}{2}}C^{\lambda}_{k-n}\right\rangle
\lambda^{-\frac{n-1}{2}}C^{\lambda+1}_{n-1}(\frac{y}{\sqrt{\lambda}})\lambda^{-\frac{k-n-1}{2}}C^{\lambda+1}_{k-n-1}(\frac{y}{\sqrt{\lambda}})\\
&&\;\hspace{0.5cm}\times\frac{4(\frac{n}{\lambda}+1)n!(\frac{k-n}{\lambda}+1)(k-n)!}
{2(2+\frac{1}{\lambda})\ldots(2+\frac{n-1}{\lambda})2(2+\frac{1}{\lambda})\ldots(2+\frac{k-n-1}{\lambda})}\\
&&\;\hspace{0.5cm}\times\frac{1}{\{(n(\frac{n}{\lambda}+2))^{1/2}+((k-n)(\frac{k-n}{\lambda}+2))^{1/2}\}^{2}}(1-\frac{y^{2}}{\lambda})^{\lambda+1/2}e^{y^{2}}\\
&&\hspace{0.2cm}=\sum_{k=1}^{N}\sum_{n=1}^{k}\frac{n!(k-n)!}{2^{n}2^{k-n}}\frac{n^{1/2}(k-n)^{1/2}}{(n^{1/2}+(k-n)^{1/2})^{2}}\left\langle\phi,\frac{H_{n}}{n!}\right\rangle\left\langle
\phi,\frac{H_{k-n}}{(k-n)!}\right\rangle \frac{H_{n}(y)}{n!}\frac{H_{k-n}(y)}{(k-n)!}\\
&&\hspace{0.5cm}+\sum_{k=1}^{N}\sum_{n=1}^{k}\frac{n!(k-n)!}{2^{n}2^{k-n}}\frac{4}{2(n^{1/2}+(k-n)^{1/2})^{2}}\\
&&\hspace{1.5cm}\times\left\langle\phi,\frac{H_{n}}{n!}\right\rangle\left\langle
\phi,\frac{H_{k-n}}{(k-n)!}\right\rangle \frac{H_{n-1}(y)}{(n-1)!}\frac{H_{k-n-1}(y)}{(k-n-1)!}\\
&& \hspace{0.2cm}=\sum_{k=1}^{N}\sum_{n=1}^{k}n^{1/2}(k-n)^{1/2}\left\langle\phi,\frac{H_{n}}{\left\|H_{n}\right\|^{2}}\right\rangle\left\langle
\phi,\frac{H_{k-n}}{\left\|H_{k-n}\right\|^{2}}\right\rangle H_{n}(y)H_{k-n}(y)\\&&\;\;\;\times\int_{0}^{\infty}te^{-t(n^{1/2}+(k-n)^{1/2})}dt\\
\end{eqnarray*}
\begin{eqnarray*}
&&\hspace{0.5cm}+\sum_{k=1}^{N}\sum_{n=1}^{k}\frac{4n(k-n)}{2}\left\langle\phi,\frac{H_{n}}{\left\|H_{n}\right\|^{2}}\right\rangle\left\langle
\phi,\frac{H_{k-n}}{\left\|H_{k-n}\right\|^{2}}\right\rangle H_{n-1}(y)H_{k-n-1}(y)\\&&\;\;\;\times\int_{0}^{\infty}te^{-t(n^{1/2}+(k-n)^{1/2})}dt\\
&&\hspace{0.2cm}=\int_{0}^{\infty}\sum_{k=1}^{N}\sum_{n=1}^{k}n^{1/2}(k-n)^{1/2}te^{-t(n^{1/2}+(k-n)^{1/2})}\\&&\hspace{2cm}\times\left\langle\phi,\frac{H_{n}}{\left\|H_{n}\right\|^{2}}\right\rangle\left\langle
\phi,\frac{H_{k-n}}{\left\|H_{k-n}\right\|^{2}}\right\rangle H_{n}(y)H_{k-n}(y)dt\\
&&\hspace{0.5cm}+\int_{0}^{\infty}\sum_{k=1}^{N}\sum_{n=1}^{k}\frac{4n(k-n)}{2}te^{-t(n^{1/2}+(k-n)^{1/2})}\\&&\hspace{2cm}\times\left\langle\phi,\frac{H_{n}}{\left\|H_{n}\right\|^{2}}\right\rangle\left\langle
\phi,\frac{H_{k-n}}{\left\|H_{k-n}\right\|^{2}}\right\rangle H_{n-1}(y)H_{k-n-1}(y)dt\\
&&\hspace{0.2cm}=\int_{0}^{\infty}\left(\sqrt{t}\sum_{k=1}^{N}k^{1/2}e^{-tk^{1/2}}\left\langle\phi,\frac{H_{k}}{\left\|H_{k}\right\|^{2}}\right\rangle H_{k}(y)\right)^{2}dt\\
&&\;\;+\int_{0}^{\infty}t\left(\frac{1}{\sqrt{2}}\sum_{k=1}^{N}2ke^{-tk^{1/2}}\left\langle\phi,\frac{H_{k}}{\left\|H_{k}\right\|^{2}}\right\rangle H_{k-1}(y)\right)^{2}dt,
\end{eqnarray*}
so we have,
\begin{eqnarray*}
\lim_{\lambda\to\infty }F^{N}_{\lambda,K}(y)&=&\int_{0}^{\infty}\left(\sqrt{t}\sum_{k=1}^{N}k^{1/2}e^{-tk^{1/2}}\left\langle\phi,\frac{H_{k}}{\left\|H_{k}\right\|^{2}}\right\rangle H_{k}(y)\right)^{2}dt\\
&&\;\;+\int_{0}^{\infty}t\left(\frac{1}{\sqrt{2}}\sum_{k=1}^{N}2ke^{-tk^{1/2}}\left\langle\phi,\frac{H_{k}}{\left\|H_{k}\right\|^{2}}\right\rangle H_{k-1}(y)\right)^{2}dt.\\
\end{eqnarray*}
Thus,
\begin{eqnarray*}
F^{N_{m}}_{K}(y)&=&\int_{0}^{\infty}\left(\sqrt{t}\sum_{k=1}^{N_{m}}k^{1/2}e^{-tk^{1/2}}\left\langle\phi,\frac{H_{k}}{\left\|H_{k}\right\|^{2}}\right\rangle H_{k}(y)\right)^{2}dt\\
&&\;\;+\int_{0}^{\infty}t\left(\frac{1}{\sqrt{2}}\sum_{k=1}^{N_{m}}2ke^{-tk^{1/2}}\left\langle\phi,\frac{H_{k}}{\left\|H_{k}\right\|^{2}}\right\rangle H_{k-1}(y)\right)^{2}dt.\\
\end{eqnarray*}
On the other hand, given that $$|H_{k}(x)|\leq \eta
e^{\frac{x^{2}}{2}}2^{\frac{k}{2}}\sqrt{k!},\hspace{1.5cm}
\eta=1.086435$$ we get,
\begin{eqnarray*}
\left|k^{1/2}e^{-tk^{1/2}}\left\langle\phi,\frac{H_{k}}{\left\|H_{k}\right\|^{2}}\right\rangle|H_{k}(y)|\right|&\leq
&k^{1/2}e^{-tk^{1/2}}\left\|\phi\right\|\frac{\eta
e^{\frac{y^{2}}{2}}2^{\frac{k}{2}}\sqrt{k!}}{\pi^{1/4}2^{\frac{k}{2}}\sqrt{k!}}\\
&=&\frac{\eta\left\|\phi\right\|e^{\frac{y^{2}}{2}}}{\pi^{1/4}}k^{1/2}e^{-tk^{1/2}},
\end{eqnarray*}
then,
\begin{eqnarray*}
\left|t\left(\sum_{k=1}^{N_{m}}k^{1/2}e^{-tk^{1/2}}\left\langle\phi,\frac{H_{k}}{\left\|H_{k}\right\|^{2}}\right\rangle
H_{k}(y)\right)^{2}\right|&\leq
&\frac{\eta^{2}\left\|\phi\right\|^{2}e^{y^{2}}}{\pi^{1/2}}t\left(\sum_{k=1}^{N_{m}}k^{1/2}e^{-tk^{1/2}}\right)^{2}\\
&\leq &\frac{\eta^{2}\left\|\phi\right\|^{2}e^{y^{2}}}{\pi^{1/2}}t\left(\sum_{k=1}^{\infty}k^{1/2}e^{-tk^{1/2}}\right)^{2}.\\
\end{eqnarray*}
Now, for $t\geq 1$ we have
\begin{eqnarray*}
\sum_{k=1}^{\infty}k^{1/2}e^{-t(k^{1/2}-1/2)}&\leq
&\sum_{k=1}^{\infty}k^{1/2}e^{-k^{1/2}+1/2},
\end{eqnarray*}
which converges. Therefore
\begin{eqnarray*}
\int_{1}^{\infty}t\left(\sum_{k=1}^{\infty}k^{1/2}e^{-tk^{1/2}}\right)^{2}dt&\leq
&\int_{1}^{\infty}te^{-t}\left(\sum_{k=1}^{\infty}k^{1/2}e^{-t(k^{1/2}-1/2)}\right)^{2}dt\\
&\leq
&\int_{1}^{\infty}te^{-t}\left(\sum_{k=1}^{\infty}k^{1/2}e^{-k^{1/2}+1/2}\right)^{2}dt \leq C\int_{1}^{\infty}te^{-t}dt<\infty.\\
\end{eqnarray*}
Let $0<t<1.$ For $y\in[-K,K]$, let us define
$$\tau_{m}(y)=\sum_{k=1}^{N_{m}}k^{1/2}e^{-tk^{1/2}}\left\langle\phi,\frac{H_{k}}{\left\|H_{k}\right\|^{2}}\right\rangle
H_{k}(y).$$
 We know that for $y\in[-K,K]$,
$$\tau_{m}(y)\longrightarrow \sum_{k=1}^{\infty}k^{1/2}e^{-tk^{1/2}}\left\langle\phi,\frac{H_{k}}{\left\|H_{k}\right\|^{2}}\right\rangle
H_{k}(y), \hspace{1cm}\mbox{if}\hspace{0.3cm}m\to\infty, $$
and also,
 $$\sum_{k=1}^{\infty}k^{1/2}e^{-tk^{1/2}}\left\langle\phi,\frac{H_{k}}{\left\|H_{k}\right\|^{2}}\right\rangle
H_{k}(y)=\frac{\partial}{\partial t}P_{t}^{\gamma}(\phi(y)).$$
 Now, for $\varepsilon=1$, exists $M(\varepsilon,y)>0$ such that
 \begin{equation}\label{limPt}
|\tau_{m}(y)-\frac{\partial}{\partial
t}P_{t}^{\gamma}(\phi(y))|<1,\hspace{0.5cm}\mbox{if}\hspace{0.3cm}m\geq
M,
\end{equation}
then, for $m\geq M$ we get
$$|\tau_{m}(y)|<1+|\frac{\partial}{\partial t}P_{t}^{\gamma}(\phi(y))|<t\left(1+C(1 + |y|)e^{-t}\right).$$
Therefore, for $y\in[-K,K]$ we have
$$\int_{0}^{1}t(1+C(1 + |y|)e^{-t})^{2}dt<\infty.$$
Let us then define
\begin{eqnarray*}
\psi(t)=\left\{\begin{array}{lcl}
\frac{\eta^{2}\left\|\phi\right\|^{2}e^{y^{2}}}{\pi^{1/2}}t\left(\sum_{k=1}^{\infty}k^{1/2}e^{-tk^{1/2}}\right)^{2} & if& t\geq 1\\
   \\
 t\left(1+C(1 + |y|)e^{-t}\right)^{2} &if& 0<t<1\\
\end{array}
\right.
\end{eqnarray*}
which is integrable in $(0,\infty)$, now for all $m>0$ we have
\begin{equation*}
\left|\left(\sqrt{t}\sum_{k=1}^{N_{m}}k^{1/2}e^{-tk^{1/2}}\left\langle\phi,\frac{H_{k}}{\left\|H_{k}\right\|^{2}}\right\rangle
H_{k}(y)\right)^{2}\right|\leq \psi(t), \hspace{0.4cm}\mbox{for
all}\hspace{0.4cm}t>0,
\end{equation*}
then, by Lebesgue's dominated convergence theorem we have
\begin{eqnarray*}
&&\lim_{m\to\infty
}\int_{0}^{\infty}\left(\sqrt{t}\sum_{k=1}^{N_{m}}k^{1/2}e^{-tk^{1/2}}\left\langle\phi,\frac{H_{k}}{\left\|H_{k}\right\|^{2}}\right\rangle
H_{k}(y)\right)^{2}dt\\&=&\int_{0}^{\infty}\left(\sqrt{t}\sum_{k=1}^{\infty}k^{1/2}e^{-tk^{1/2}}\left\langle\phi,\frac{H_{k}}{\left\|H_{k}\right\|^{2}}\right\rangle
H_{k}(y)\right)^{2}dt
=\int_{0}^{\infty}t\left(\frac{\partial}{\partial
t}P_{t}^{\gamma}(\phi(y))\right)^{2}dt.
\end{eqnarray*}
Similarly, we have
\begin{equation*}
\lim_{m\to\infty
}\int_{0}^{\infty}t\left(\frac{1}{\sqrt{2}}\sum_{k=1}^{N_{m}}2ke^{-tk^{1/2}}\left\langle\phi,\frac{H_{k}}{\left\|H_{k}\right\|^{2}}\right\rangle
H_{k-1}(y)\right)^{2}dt=\int_{0}^{\infty}t\left(\frac{1}{\sqrt{2}}\frac{\partial}{\partial
y}P_{t}^{\gamma}(\phi(y))\right)^{2}dt
\end{equation*}
i.e.,
\begin{equation}\label{F-Herm2}
F^{N_{m}}_{K}(y)\to
\left(g^{\gamma}\phi(y)\right)^{2},\hspace{0.3cm}\mbox{cuando}\hspace{0.3cm}m\to
\infty\hspace{0.3cm}a.e.\hspace{0.3cm}y\in\mathbb{R}
\end{equation}
therefore, from (\ref{F-Herm}) and (\ref{F-Herm2}) we get
\begin{equation*}
F(y)=g^{\gamma}(\phi(y))
\end{equation*}

\item [ii)] The proof  is essentially analogous to i), so fewer details will be provided. Assume that the operator
$g^{(\alpha,\beta)}$ is bounded in
$L^{p}\left([-1,1],\mu_{\alpha,\beta}\right)$. Let $\phi\in
C^{\infty}_{0}(0,\infty)$ and $\beta>0,$ then we have
\begin{equation*}
\|g^{(\alpha,\beta)}\phi_{\beta}\|_{L^{p}\left([-1,1],\mu_{(\alpha,\beta)}\right)}\leq
C\|\phi_{\beta}\|_{L^{p}\left([-1,1],\mu_{(\alpha,\beta)}\right)}
\end{equation*}
where, $\phi_{\beta}(x)=\phi\left(\frac{\beta}{2}(1-x)\right),\;$
for $x\in[-1,1]$, i.e.,
\begin{eqnarray*}
\left\|\left\{\int_{0}^{\infty}t|\nabla_{(\alpha,\beta)}
P^{(\alpha,\beta)}_{t}\phi_{\beta}|^{2}dt
\right\}^{1/2}\right\|_{L^{p}\left([-1,1],\mu_{(\alpha,\beta)}\right)}\leq
C\|\phi_{\beta}\|_{L^{p}\left([-1,1],\mu_{(\alpha,\beta)}\right)}
\end{eqnarray*}
Now making the change of variable $x=1-\frac{2}{\beta}y,\;$ we
have
\begin{eqnarray*}
\left\{\eta_{\alpha,\beta}\frac{2^{\alpha+\beta+1}}{\beta^{\alpha+1}}\int_{0}^{\beta}\left|\left\{\int_{0}^{\infty}t|\nabla_{(\alpha,\beta)}
P^{(\alpha,\beta)}_{t}\phi_{\beta}(1-\frac{2}{\beta}y)|^{2}dt
\right\}^{1/2}(1-\frac{y}{\beta})^{\beta/p}e^{\frac{y}{p}}\right|^{p}y^{\alpha}e^{-y}dy\right\}^{1/p}\hspace{3cm}
\end{eqnarray*}
$$\leq
C\|\phi_{\beta}\|_{L^{p}\left([-1,1],\mu_{\lambda}\right)}$$ thus,
\begin{eqnarray*}
\left\{\int_{0}^{\beta}\left|\left\{\int_{0}^{\infty}t|\nabla_{(\alpha,\beta)}
P^{(\alpha,\beta)}_{t}\phi_{\beta}(1-\frac{2}{\beta}y)|^{2}dt
\right\}^{1/2}(1-\frac{y}{\beta})^{\beta/p}e^{\frac{y}{p}}\right|^{p}y^{\alpha}e^{-y}dy\right\}^{1/p}\hspace{3cm}
\end{eqnarray*}
$$\leq
C\left(Z_{\alpha,\beta}\right)^{-1}\|\phi_{\beta}\|_{L^{p}\left([-1,1],\mu_{\lambda}\right)}$$
where,
$Z_{\alpha,\beta}=\frac{\Gamma(\alpha+\beta+2)}{\Gamma(\alpha+1)\beta^{\alpha+2}\Gamma(\beta)}.\;$
Analogously we get,
\begin{eqnarray*}
\left\{\int_{0}^{\beta}\left|\left\{\int_{0}^{\infty}t|\nabla_{(\alpha,\beta)}
P^{(\alpha,\beta)}_{t}\phi_{\beta}(1-\frac{2}{\beta}y)|^{2}dt
\right\}^{1/2}(1-\frac{y}{\beta})^{\beta/2}e^{\frac{y}{2}}\right|^{2}y^{\alpha}e^{-y}dy\right\}^{1/2}\hspace{3cm}
\end{eqnarray*}
$$\leq
C\left(Z_{\alpha,\beta}\right)^{-1}\|\phi_{\beta}\|_{L^{2}\left([-1,1],\mu_{\lambda}\right)}$$
Now, for each $K\in\mathbb{N} $ and $\beta>0$, such that
$\beta>K,$ define the functions
\begin{eqnarray*}
F_{\beta ,K}(y)=
 \chi_{(0,K)}(y)\left\{\int_{0}^{\infty}t|\nabla_{(\alpha,\beta)}
P^{(\alpha,\beta)}_{t}\phi_{\beta}(1-\frac{2}{\beta}y)|^{2}dt
\right\}^{1/2}(1-\frac{y}{\beta})^{\beta/2}e^{\frac{y}{2}}
\end{eqnarray*}
and
\begin{eqnarray*}
f_{\beta ,K}(y)=
 \chi_{(0,K)}(y)\left\{\int_{0}^{\infty}t|\nabla_{(\alpha,\beta)}
P^{(\alpha,\beta)}_{t}\phi_{\beta}(1-\frac{2}{\beta}y)|^{2}dt
\right\}^{1/2}(1-\frac{y}{\beta})^{\beta/p}e^{\frac{y}{p}}
\end{eqnarray*}
From the previous inequalities  both series converge for all $y\in
(0,\beta)$, and $F_{\beta , K}=f_{\beta, K}\Omega_{\beta},\;$
where
$$\Omega_{\beta}(y)=e^{\frac{y}{2}-\frac{y}{p}}\left(1-\frac{y}{\beta}\right)^{\beta/2-\beta/p}$$
for all  $K\in\mathbb{N}\;\;$ and $\beta>K$. Now as,
\begin{eqnarray*}
|\Omega_{\beta}(y)|&=&e^{\frac{y}{2}-\frac{y}{p}}\left(1-\frac{y}{\beta}\right)^{\beta/2-\beta/p}
\leq e^{\frac{y}{2}-\frac{y}{p}}(e^{-y/\beta})^{\beta/2-\beta/p}=1
\end{eqnarray*}
we conclude that $\Omega_{\beta}$ is bounded in $(0,K)$. On the
other hand,
\begin{eqnarray*}
(Z_{\alpha,\beta})^{-1/p}\|\phi_{\beta}\|_{L^{p}\left([-1,1]\mu_{(\alpha,\beta)}\right)}&=&
(Z_{\alpha,\beta})^{-1/p}\left\{\eta_{\alpha,\beta}\int_{-1}^{1}\left|\phi_{\beta}(x)\right|^{p}(1-x)^{\alpha}(1+x)^{\beta}dx\right\}^{1/p},
\end{eqnarray*}
and making the change of variable $x=1-\frac{2}{\beta}y$ we get
\begin{eqnarray*}
&&(Z_{\alpha,\beta})^{-1/p}\left\{\eta_{\alpha,\beta}\int_{-1}^{1}\left|\phi_{\beta}(x)\right|^{p}(1-x)^{\alpha}(1+x)^{\beta}dx\right\}^{1/p}\hspace{5cm}
\\&&\hspace{0.5cm}=
(Z_{\alpha,\beta})^{-1/p}\left\{\eta_{\alpha,\beta}\frac{2}{\beta}\int_{0}^{\beta}\left|\phi_{\beta}(1-\frac{2}{\beta}y)\right|^{p}\left(\frac{2y}{\beta}\right)^{\alpha}\left(2-\frac{2y}{\beta}\right)^{\beta}dy\right\}^{1/p}\\
&&\hspace{0.5cm}\leq
C\left\{\int_{0}^{\beta}\left|\phi(y)\right|^{p}y^{\alpha}\left(1-\frac{y}{\beta}\right)^{\beta}dy\right\}^{1/p}\leq
C\left\{\int_{0}^{\infty}\left|\phi(y)\right|^{p}y^{\alpha}e^{-y}dy\right\}^{1/p}=C\|\phi\|_{{L^{p}\left((0,\infty),\mu_{\alpha}\right)}}.
\end{eqnarray*}
Then,
\begin{eqnarray*}
\lim_{\beta\to
\infty}(Z_{\alpha,\beta})^{-1/p}\|\phi_{\beta}\|_{{L^{p}\left([-1,1],\mu_{(\alpha,\beta)}\right)}}&\leq
& \lim_{\beta\to \infty}C\left\{\int_{0}^{\infty}\left|\phi(y)\right|^{p}y^{\alpha}e^{-y}dy\right\}^{1/p}\\
&=&C\|\phi\|_{{L^{p}\left((0,\infty),\mu_{\alpha}\right)}}.
\end{eqnarray*}
Moreover,
\begin{eqnarray}\label{desLpfiLag}
(Z_{\alpha,\beta})^{-1/p}\|\phi_{\beta}\|_{{L^{p}\left([-1,1],\mu_{(\alpha,\beta)}\right)}}\leq
C\|\phi\|_{{L^{p}\left((0,\infty),\mu_{\alpha}\right)}}.
\end{eqnarray}
Now,
\begin{eqnarray}\label{norF1gausLag}
\nonumber \|F_{\beta
,K}\|_{{L^{2}\left((0,\infty),\mu_{\alpha}\right)}}&\leq
&C(Z_{\alpha,\beta})^{-1/2}\|\phi_{\beta}\|_{L^{2}\left([-1,1],\mu_{(\alpha,\beta)}\right)}
\end{eqnarray}
and therefore, from (\ref{desLpfiLag}) for $p=2$ and
(\ref{norF1gausLag}), we get
\begin{eqnarray*}
\|F_{\beta ,K}\|_{{L^{2}\left((0,\infty),\mu_{
\alpha}\right)}}\leq
C\|\phi\|_{{L^{2}\left((0,\infty),\mu_{\alpha}\right)}}.
\end{eqnarray*}
Analogously, using that $ \Omega_{\beta}$ is bounded in $(0,K)$,
\begin{eqnarray}\label{norf1gausLag}
\nonumber \|F_{\beta
,K}\|_{{L^{p}\left((0,\infty),\mu_{\alpha}\right)}} &\leq
&C\left\{\frac{1}{\Gamma(\alpha+1)}\int_{0}^{k}|f_{\beta
,K}(y)|^{p}y^{\alpha}e^{-y}dy\right\}^{1/p}\\
&\leq
&C(Z_{\alpha,\beta})^{-1/p}\|\phi_{\beta}\|_{L^{p}\left([-1,1],\mu_{(\alpha,\beta)}\right)}.
\end{eqnarray}
Then, from (\ref{desLpfiLag}) and (\ref{norf1gausLag}),
\begin{eqnarray*}
\|F_{\beta ,K}\|_{{L^{p}\left((0,\infty),\mu_{\alpha}\right)}}
\leq C\|\phi\|_{{L^{p}\left((0,\infty),\mu_{\alpha}\right)}},
\end{eqnarray*}
for all $\beta>K$. Therefore $\{F_{\beta,K}\}$ is a bounded
subsequence in $\;{L^{2}\left((0,\infty),\mu_{\alpha}\right)}\;$
with $\;{L^{p}\left((0,\infty),\mu_{\alpha}\right)}.$ Thus by the
 Bourbaki-Alaoglu's theorem, there exists an increasing sequence $\{\beta_{j}\}_{j\in\mathbb{N}}\;$ with $\lim_{j\to
\infty}\beta_{j}=\infty$ and functions
$\;F_{K}\in{{L^{2}\left((0,\infty),\mu_{\alpha}\right)}}$ and
$f_{K}\in{{L^{p}\left((0,\infty),\mu_{\alpha}\right)}}$ satisfying
that

\begin{itemize}
\item $ F_{\beta_j,K}\to  F_{K},$ as $j\to\infty,$ in the weak
topology of ${L^{2}\left((0,\infty),\mu_{\alpha}\right)}$ \item
$F_{\beta_j,K}\to  f_{K},$ as $j\to\infty,$ in the weak topology
of ${L^{p}\left((0,\infty),\mu_{\alpha}\right)}.$
\end{itemize}
Then, as in i), we can conclude that there exists an increasing
sequence
 $(\beta_{j})_{j\in
\mathbb{N}}\subset (0,\infty)$ such that $\lim_{j\to
\infty}\beta_{j}=\infty,$ and a function $F\in
L^{p}\left((0,\infty),\mu_{\alpha}\right)\cap
L^{2}\left((0,\infty),\mu_{\alpha}\right),$ such that

\begin{itemize}
\item for each $K\in \mathbb{N},\;$  $ F_{\beta_{j},K}\to F,$ as
$j\to\infty,$ in the weak topology on
${L^{2}\left((0,\infty),\mu_{\alpha}\right)}$ and in the weak
topology on ${L^{p}\left((0,\infty),\mu_{\alpha}\right)}$.
 \item $\|F\|_{L^{p}\left((0,\infty),\mu_{\alpha}\right)}\leq
C\|\phi\|_{L^{p}\left((0,\infty),\mu_{\alpha}\right)}$.
\end{itemize}

Analogously to the Hermite case, we have
\begin{eqnarray*}
F_{\beta ,K}(y)&=&\chi_{(0,K)}(y)
\left\{\int_{0}^{\infty}t\left(\frac{\partial}{\partial
t}P^{(\alpha,\beta)}_{t}\phi_{\beta}(1-\frac{2}{\beta}y)\right)^{2}dt(1-\frac{y}{\beta})^{\beta}e^{y}\right.\\&&\hspace{0.5cm}+\left.\int_{0}^{\infty}\frac{t}{\beta}\left(\frac{\partial}{\partial
y}P^{(\alpha,\beta)}_{t}\phi_{\beta}(1-\frac{2}{\beta}y)\right)^{2}dt4y(1-\frac{y}{\beta})^{\beta+1}e^{y}
\right\}^{1/2}.\\
\end{eqnarray*}
Now. let us define
\begin{equation*}
g_{1,\beta}\phi(y)=\chi_{(0,K)}(y)\int_{0}^{\infty}t\left(\frac{\partial}{\partial
t}P^{(\alpha,\beta)}_{t}\phi_{\beta}(1-\frac{2}{\beta}y)\right)^{2}dt(1-\frac{y}{\beta})^{\beta}e^{y}
\end{equation*}
and
\begin{equation*}
g_{2,\beta}\phi(y)=\chi_{(0,K)}(y)\int_{0}^{\infty}\frac{t}{\beta}\left(\frac{\partial}{\partial
y}P^{(\alpha,\beta)}_{t}\phi_{\beta}(1-\frac{2}{\beta}y)\right)^{2}dt4y(1-\frac{y}{\beta})^{\beta+1}e^{y}
\end{equation*}
Let us first $g_{1,\beta}$.\\
For a function $\phi\in L^{2}((0,\infty),\mu_{\alpha})$
\begin{eqnarray*}
&&\left(\sqrt{t}\frac{\partial}{\partial
t}P^{(\alpha,\beta)}_{t}\phi_{\beta}(1-\frac{2}{\beta}y)\right)^{2}=\sum_{k=1}^{\infty}\sum_{n=1}^{k}\left\langle
\phi_{\beta},\frac{P^{(\alpha,\beta)}_{n}}{\left\|P^{(\alpha,\beta)}_{n}\right\|_{2}^{2}}\right\rangle\left\langle
\phi_{\beta},\frac{P^{(\alpha,\beta)}_{k-n}}{\left\|P^{(\alpha,\beta)}_{k-n}\right\|_{2}^{2}}\right\rangle
t\lambda_{n}^{1/2}\lambda_{k-n}^{1/2}\\&&\hspace{6cm}\times
\;e^{-t(\lambda_{k}^{1/2}+\lambda_{k-n}^{1/2})}P^{(\alpha,\beta)}_{n}(1-\frac{2}{\beta}y)P^{(\alpha,\beta)}_{k-n}(1-\frac{2}{\beta}y).
\end{eqnarray*}
Therefore, as in the previous case
\begin{eqnarray*}
g_{1,\beta}\phi(y)&=&\sum_{k=1}^{N}\sum_{n=1}^{k}\left\langle
\phi_{\beta},\frac{P^{(\alpha,\beta)}_{n}}{\left\|P^{(\alpha,\beta)}_{n}\right\|_{2}^{2}}\right\rangle\left\langle
\phi_{\beta},\frac{P^{(\alpha,\beta)}_{k-n}}{\left\|P^{(\alpha,\beta)}_{k-n}\right\|_{2}^{2}}\right\rangle
\frac{\lambda_{n}^{1/2}\lambda_{k-n}^{1/2}}{(\lambda_{k}^{1/2}+\lambda_{k-n}^{1/2})^{2}}\\
&&\hspace{0.3cm}\times\;P^{(\alpha,\beta)}_{n}(1-\frac{2}{\beta}y)P^{(\alpha,\beta)}_{k-n}(1-\frac{2}{\beta}y)(1-\frac{y}{\beta})^{\beta}e^{y}\\
&&\hspace{0.5cm}+\;\int_{0}^{\infty}H_{\beta,K}^{N,1}(y)dt\;(1-\frac{y}{\beta})^{\beta/2}e^{y/2},\\
\end{eqnarray*}
where,
\begin{eqnarray*}
H_{\beta,K}^{N,1}(y)&=&\chi_{(0,K)}(y)\sum_{k=N+1}^{\infty}\sum_{n=1}^{k}\left\langle
\phi_{\beta},\frac{P^{(\alpha,\beta)}_{n}}{\left\|P^{(\alpha,\beta)}_{n}\right\|_{2}^{2}}\right\rangle\left\langle
\phi_{\beta},\frac{P^{(\alpha,\beta)}_{k-n}}{\left\|P^{(\alpha,\beta)}_{k-n}\right\|_{2}^{2}}\right\rangle
t\lambda_{n}^{1/2}\lambda_{k-n}^{1/2}\\&&\hspace{0.7cm}\times\;
e^{-t(\lambda_{k}^{1/2}+\lambda_{k-n}^{1/2})}P^{(\alpha,\beta)}_{n}(1-\frac{2}{\beta}y)P^{(\alpha,\beta)}_{k-n}(1-\frac{2}{\beta}y)(1-\frac{y}{\beta})^{\beta/2}e^{y/2}.\\
\end{eqnarray*}
We want to prove that for $K\in \mathbb{N}$ and $\beta>K,$
\begin{equation*}
\int_{0}^{\infty}\left|\int_{0}^{\infty}H_{\beta,K}^{N,1}(y)dt\right|^{2}
y^{\alpha}e^{-y}\frac{dy}{\Gamma(\alpha+1)}\leq
\frac{C\left\|\phi\right\|_{2,\alpha}^{4}}{e^{N+1}}.
\end{equation*}
Now,
\begin{eqnarray*}
&&\int_{0}^{\infty}\left|H_{\beta,K}^{N,1}(y)\right|^{2}y^{\alpha}e^{-y}\frac{dy}{\Gamma(\alpha+1)}\\
&\leq
&\frac{\beta^{\alpha}}{2^{\alpha+\beta}\Gamma(\alpha+1)}\int_{0}^{\infty}\left|\sum_{k=N+1}^{\infty}\sum_{n=1}^{k}\left\langle
\phi_{\beta},\frac{P^{(\alpha,\beta)}_{n}}{\left\|P^{(\alpha,\beta)}_{n}\right\|_{2}^{2}}\right\rangle\left\langle
\phi_{\beta},\frac{P^{(\alpha,\beta)}_{k-n}}{\left\|P^{(\alpha,\beta)}_{k-n}\right\|_{2}^{2}}\right\rangle
t\lambda_{n}^{1/2}\lambda_{k-n}^{1/2}e^{-t(\lambda_{k}^{1/2}+\lambda_{k-n}^{1/2})} \right.\\
&&\;\hspace{4cm}\times\left. P^{(\alpha,\beta)}_{n}(1-\frac{2}{\beta}y)P^{(\alpha,\beta)}_{k-n}(1-\frac{2}{\beta}y)\right|^{2}\left(\frac{2y}{\beta}\right)^{\alpha}2^{\beta}(1-\frac{y}{\beta})^{\beta}dy\\
&=&\frac{\beta^{\alpha+1}}{2^{\alpha+\beta+1}\Gamma(\alpha+1)}\int_{-1}^{1}\left|\sum_{k=N+1}^{\infty}\sum_{n=1}^{k}\left\langle
\phi_{\beta},\frac{P^{(\alpha,\beta)}_{n}}{\left\|P^{(\alpha,\beta)}_{n}\right\|_{2}^{2}}\right\rangle\left\langle
\phi_{\beta},\frac{P^{(\alpha,\beta)}_{k-n}}{\left\|P^{(\alpha,\beta)}_{k-n}\right\|_{2}^{2}}\right\rangle
t\lambda_{n}^{1/2}\lambda_{k-n}^{1/2}e^{-t(\lambda_{k}^{1/2}+\lambda_{k-n}^{1/2})} \right.\\
&&\;\hspace{4cm}\times\left. P^{(\alpha,\beta)}_{n}(x)P^{(\alpha,\beta)}_{k-n}(x)\right|^{2}(1-x)^{\alpha}(1+x)^{\beta}dx.\\
\end{eqnarray*}
Then, again by Gasper's linearization of the product of Jacobi
polynomials \cite{gasp1} and using Parseval's identity, we have
\begin{eqnarray*}
&&\int_{0}^{\infty}\left|H_{\beta,K}^{N,1}(y)\right|^{2}y^{\alpha}e^{-y}\frac{dy}{\Gamma(\alpha+1)}\\
&\leq
&\frac{\beta^{\alpha+1}}{\eta_{\alpha,\beta}2^{\alpha+\beta+1}\Gamma(\alpha+1)}\sum_{k=N+1}^{\infty}\sum_{n=1}^{k}\sum_{i=k-2n}^{k}\left|\left\langle
\phi_{\beta},\frac{P^{(\alpha,\beta)}_{n}}{\left\|P^{(\alpha,\beta)}_{n}\right\|_{2}^{2}}\right\rangle\right|^{2}\left|\left\langle
\phi_{\beta},\frac{P^{(\alpha,\beta)}_{k-n}}{\left\|P^{(\alpha,\beta)}_{k-n}\right\|_{2}^{2}}\right\rangle\right|^{2}
t^{2}\lambda_{n}\lambda_{k-n} \\
&&\;\hspace{1cm}\times \left|P^{(\alpha,\beta)}_{n}(1)P^{(\alpha,\beta)}_{k-n}(1)\right|^{2}e^{-2t(\lambda_{k}^{1/2}+\lambda_{k-n}^{1/2})}\left|\nu(i,k-n,n)\right|^{2} \left\|p^{(\alpha,\beta)}_{i}\right\|^{2}_{2}.\\
\end{eqnarray*}
Using analogous boundedness argument as in the Hermite case, we get
\begin{eqnarray*}
&&\int_{0}^{\infty}\left|H_{\beta,K}^{N,1}(y)\right|^{2}y^{\alpha}e^{-y}\frac{dy}{\Gamma(\alpha+1)}\\
& \leq &
\frac{\beta^{\alpha+1}\left\|\phi_{\beta}\right\|^{2}}{\eta_{\alpha,\beta}2^{\alpha+\beta+1}\Gamma(\alpha+1)}\sum_{k=N+1}^{\infty}\sum_{n=1}^{k}\left|\left\langle
\phi_{\beta},\frac{P^{(\alpha,\beta)}_{n}}{\left\|P^{(\alpha,\beta)}_{n}\right\|_{2}}\right\rangle\right|^{2} \frac{t^{2}\lambda_{n}\lambda_{k-n}2n\left(%
\begin{array}{c}
  n+q \\
  \\
  n\\
\end{array}%
\right)^{2}e^{-2t(\lambda_{k}^{1/2}+\lambda_{k-n}^{1/2})}}{\left\|P^{(\alpha,\beta)}_{n}\right\|_{2}^{2}}\\
& \leq &
\frac{e^{-t}\beta^{\alpha+1}\left\|\phi_{\beta}\right\|^{2}}{\eta_{\alpha,\beta}2^{\alpha+\beta+1}\Gamma(\alpha+1)}\sum_{k=N+1}^{\infty}\sum_{n=1}^{k}\left|\left\langle
\phi_{\beta},\frac{P^{(\alpha,\beta)}_{n}}{\left\|P^{(\alpha,\beta)}_{n}\right\|_{2}}\right\rangle\right|^{2} \frac{2n\left(%
\begin{array}{c}
  n+q \\
  \\
  n\\
\end{array}%
\right)^{2}e^{-t\lambda_{k-n}^{1/2}}}{t^{2}\left\|P^{(\alpha,\beta)}_{n}\right\|_{2}^{2}}.\\
\end{eqnarray*}
Now, for $\beta $ big enough
\begin{eqnarray*}
\frac{1}{\left\|P^{(\alpha,\beta)}_{n}\right\|_{2}^{2}}&=&\frac{(2n+\alpha
+\beta +1) \Gamma(\alpha+1)\Gamma(\beta+1)\Gamma \left( n+1\right)
\Gamma \left( n+\alpha +\beta +1\right)
}{\Gamma(\alpha+\beta+2)\Gamma \left( n+\alpha +1\right) \Gamma
\left( n+\beta +1\right) }\leq \frac{C n!}{(\alpha+1)_{n}}
\end{eqnarray*}
and

\begin{equation}\label{comb}
\left(%
\begin{array}{c}
  n+q \\
  \\
  n\\
\end{array}%
\right)=\left(%
\begin{array}{c}
  n+\beta \\
  \\
  n\\
\end{array}%
\right)=\frac{\Gamma(n+\beta+1)}{\Gamma(\beta+1)n!}\sim
\frac{\beta^{n}}{n!}.\hspace{4cm}
\end{equation}
Thus, for $\beta $ big enough,
\begin{eqnarray*}
&&\int_{0}^{\infty}\left|H_{\beta,K}^{N,1}(y)\right|^{2}y^{\alpha}e^{-y}\frac{dy}{\Gamma(\alpha+1)}\\
&\leq
&\frac{Ce^{-t}\beta^{\alpha+1}\left\|\phi_{\beta}\right\|^{2}}{\eta_{\alpha,\beta}2^{\alpha+\beta+1}\Gamma(\alpha+1)}\sum_{k=N+1}^{\infty}\sum_{n=1}^{k}\left|\left\langle
\phi_{\beta},\frac{P^{(\alpha,\beta)}_{n}}{\left\|P^{(\alpha,\beta)}_{n}\right\|_{2}}\right\rangle\right|^{2}
\frac{2n\beta^{2n} e^{-t\lambda_{k-n}^{1/2}}}{t^{2}(\alpha+1)_{n}}
\end{eqnarray*}
Taking $\beta=(k-n)\left(1+2k^{2}\right)^{2},$ then, $k$ is
big enough if $\beta$ is big enough and we have two cases.\\
\begin{itemize}
\item Case $t\geq 1$: In this case, we have
\begin{eqnarray*}
\frac{1}{e^{t\lambda_{k-n}^{1/2}}}&\leq &
\frac{1}{e^{(k-n)^{1/2}\left(k-n+(k-n)\left(1+2k^{2}\right)^{2}+\alpha+1\right)^{1/2}}}\leq \frac{1}{e^{(k-n)}}\frac{1}{e^{k}}.\\
\end{eqnarray*}
Therefore, for $k$ big enough exist $C>0$ such that
\begin{eqnarray*}
\sum_{k=N+1}^{\infty}\sum_{n=1}^{k}\left|\left\langle
\phi_{\beta},\frac{P^{(\alpha,\beta)}_{n}}{\left\|P^{(\alpha,\beta)}_{n}\right\|_{2}}\right\rangle\right|^{2}
\frac{2n\beta^{2n}
e^{-t\lambda_{k-n}^{1/2}}}{t^{2}(\alpha+1)_{n}}&\leq
&\frac{C}{e^{N}}\sum_{k=N+1}^{\infty}\sum_{n=1}^{k}\left|\left\langle
\phi_{\beta},\frac{P^{(\alpha,\beta)}_{n}}{\left\|P^{(\alpha,\beta)}_{n}\right\|_{2}}\right\rangle\right|^{2},
\end{eqnarray*}
thus,
\begin{equation}\label{int-H1N}
\int_{0}^{\infty}\left|H_{\beta,K}^{N,1}(y)\right|^{2}y^{\alpha}e^{-y}\frac{dy}{\Gamma(\alpha+1)}\leq
\frac{Ce^{-t}\left\|\phi\right\|_{2,\alpha}^{4}}{e^{N+1}}
\end{equation}
\item Case $0<t<1$: for $t$ near to $0$, exists $k>0$ big enough such that
$\frac{1}{k}\leq t$, i.e.
 $\;\frac{1}{t^{2}}\leq k^{2}$.Then, analogously to the Hermite case
 \begin{eqnarray*}
\frac{1}{e^{t\lambda_{k-n}^{1/2}}}&\leq &\frac{1}{e^{(k-n)}}\frac{1}{e^{k}}\frac{1}{e^{k}}.\\
\end{eqnarray*}
Hence,
\begin{equation*}
\int_{0}^{\infty}\left|H_{\beta,K}^{N,1}(y)\right|^{2}y^{\alpha}e^{-y}\frac{dy}{\Gamma(\alpha+1)}\leq
\frac{Ce^{-t}\left\|\phi\right\|_{2,\alpha}^{4}}{e^{N+1}}.
\end{equation*}
Therefore,
\begin{eqnarray*}
\left(\int_{0}^{\infty}\left|\int_{0}^{\infty}H_{\beta,K}^{N,1}(y)dt\right|^{2}y^{\alpha}e^{-y}\frac{dy}{\Gamma(\alpha+1)}\right)^{1/2}&\leq
&\int_{0}^{\infty}\left(\frac{e^{-t}C\|\phi\|_{2,\alpha}^{4}}{e^{(N+1)}}\right)^{1/2}dt\\
&=& \frac{(C)^{1/2}\|\phi\|_{2,\alpha}^{2}}{e^{\frac{(N+1)}{2}}}\\
\end{eqnarray*}
Thus, $\left\{\int_{0}^{\infty}H_{\beta,K}^{N,1}dt\right\}$ is a
bounded sequence on
$\;L^{2}\left((0,\infty),\mu_{\alpha}\right)\;$ so by
Bourbaki-Alaoglu's theorem, there exists a sequence
$\;(\lambda_{j})_{j\in\mathbb{N}}\;$ $\lim_{j\to
\infty}\lambda_{j}=\infty$ such that, for all $N\in \mathbb{N}$,
  $\left\{\int_{0}^{\infty}H_{\beta_{j},K}^{N,1}dt\right\}_{j\in\mathbb{N}}$ converges weakly in $\;L^{2}\left((0,\infty),\mu_{\alpha}\right)\;$ to a function
  $H_{K}^{N,1}\in L^{2}\left((0,\infty),\mu_{\alpha}\right)\;$ Moreover,
  \begin{equation}\label{norL_n}
\int_{0}^{\infty}\left|H_{K}^{N,1}(y)\right|^{2}y^{\alpha}e^{-y}\frac{dy}{\Gamma(\alpha+1)}\leq \frac{C}{e^{\frac{(N+1)}{2}}}\\
\end{equation}
Then, there exists a non decreasing sequence
$(N_{j})_{j\in\mathbb{N}}$ such that,
\begin{equation}\label{limL_k}
H_{K}^{N_{j},1}(y)\longrightarrow 0,\hspace{0.2cm}
a.e.\hspace{0.2cm} y\in\mathbb{R}.
\end{equation}

\end{itemize}

For the function $g_{2,\beta}$ we get similar estimates. Given $\phi\in L^{2}((0,\infty),\mu_{\alpha}),$
\begin{eqnarray*}
\left(\frac{\sqrt{t}}{\sqrt{\beta}}\frac{\partial}{\partial
y}P^{(\alpha,\beta)}_{t}\phi_{\beta}(1-\frac{2}{\beta}y)\right)^{2}
&=&\sum_{k=1}^{\infty}\sum_{n=1}^{k}\left\langle
\phi_{\beta},\frac{P^{(\alpha,\beta)}_{n}}{\left\|P^{(\alpha,\beta)}_{n}\right\|_{2}^{2}}\right\rangle\left\langle
\phi_{\beta},\frac{P^{(\alpha,\beta)}_{k-n}}{\left\|P^{(\alpha,\beta)}_{k-n}\right\|_{2}^{2}}\right\rangle\frac{\lambda_{n}\lambda_{k-n}}{4n(k-n)}\\
&&\times
P^{(\alpha+1,\beta+1)}_{n-1}(1-\frac{2}{\beta}y)P^{(\alpha+1,\beta+1)}_{k-n-1}(1-\frac{2}{\beta}y)
\frac{t}{\beta}e^{-t(\lambda_{k}^{1/2}+\lambda_{k-n}^{1/2})}.\\
\end{eqnarray*}
Thus,
\begin{eqnarray*}
g_{2,\beta}\phi(y) &=&\sum_{k=1}^{N}\sum_{n=1}^{k}\left\langle
\phi_{\beta},\frac{P^{(\alpha,\beta)}_{n}}{\left\|P^{(\alpha,\beta)}_{n}\right\|_{2}^{2}}\right\rangle\left\langle
\phi_{\beta},\frac{P^{(\alpha,\beta)}_{k-n}}{\left\|P^{(\alpha,\beta)}_{k-n}\right\|_{2}^{2}}\right\rangle
P^{(\alpha+1,\beta+1)}_{n-1}(1-\frac{2}{\beta}y)P^{(\alpha+1,\beta+1)}_{k-n-1}(1-\frac{2}{\beta}y)
\\&&\hspace{0.5cm}\times
\;\frac{y}{\beta n(k-n)}\frac{\lambda_{n}\lambda_{k-n}}{(\lambda_{k}^{1/2}+\lambda_{k-n}^{1/2})^{2}}(1-\frac{y}{\beta})^{\beta+1}e^{y}+\;\int_{0}^{\infty}H_{\beta,K}^{N,2}(y)dt
y^{1/2}(1-\frac{y}{\beta})^{(\beta+1)/2}e^{y/2},
\end{eqnarray*}
where,
\begin{eqnarray*}
H_{\beta,K}^{N,2}(y)&=&\chi_{(0,K)}(y)\sum_{k=N+1}^{\infty}\sum_{n=1}^{k}\left\langle
\phi_{\beta},\frac{P^{(\alpha,\beta)}_{n}}{\left\|P^{(\alpha,\beta)}_{n}\right\|_{2}^{2}}\right\rangle\left\langle
\phi_{\beta},\frac{P^{(\alpha,\beta)}_{k-n}}{\left\|P^{(\alpha,\beta)}_{k-n}\right\|_{2}^{2}}\right\rangle
\frac{t}{\beta}e^{-t(\lambda_{k}^{1/2}+\lambda_{k-n}^{1/2})}\frac{\lambda_{n}\lambda_{k-n}}{n(k-n)}
\\&&\hspace{1cm}\times
P^{(\alpha+1,\beta+1)}_{n-1}(1-\frac{2}{\beta}y)P^{(\alpha+1,\beta+1)}_{k-n-1}(1-\frac{2}{\beta}y)y^{1/2}(1-\frac{y}{\beta})^{(\beta+1)/2}e^{y/2}.
\;\\
\end{eqnarray*}
Again, we want to prove that for $K\in
\mathbb{N}$ and $\beta>K,$
\begin{equation*}
\int_{0}^{\infty}\left|\int_{0}^{\infty}H_{\beta,K}^{N,2}(y)dt\right|^{2}
y^{\alpha}e^{-y}\frac{dy}{\Gamma(\alpha+1)}\leq
\frac{(C)^{1/2}\|\phi\|_{2,\alpha}^{2}}{e^{\frac{(N+1)}{2}}}.
\end{equation*}
Let us see this,
\begin{eqnarray*}
&&\int_{0}^{\infty}\left|H_{\beta,K}^{N,2}(y)\right|^{2}y^{\alpha}e^{-y}\frac{dy}{\Gamma(\alpha+1)}\\
&& \quad \leq
\frac{\beta^{\alpha+2}}{2^{\alpha+\beta+3}\Gamma(\alpha+1)}\int_{-1}^{1}\left|\sum_{k=N+1}^{\infty}\sum_{n=1}^{k}\sum_{i=k-2n}^{k-2}\left\langle
\phi_{\beta},\frac{P^{(\alpha,\beta)}_{n}}{\left\|P^{(\alpha,\beta)}_{n}\right\|_{2}^{2}}\right\rangle\left\langle
\phi_{\beta},\frac{P^{(\alpha,\beta)}_{k-n}}{\left\|P^{(\alpha,\beta)}_{k-n}\right\|_{2}^{2}}\right\rangle
\frac{t}{\beta} \right.\\&&\hspace{1cm}\times
\left.e^{-t(\lambda_{k}^{1/2}+\lambda_{k-n}^{1/2})}\frac{\lambda_{n}\lambda_{k-n}}{n(k-n)}P^{(\alpha+1,\beta+1)}_{n-1}(1)P^{(\alpha+1,\beta+1)}_{k-n-1}(1)\nu(i,k-n-1,n-1)p^{(\alpha+1,\beta+1)}_{i}(x)\right|^{2}\\
&&\hspace{3cm}\times\left(1-x\right)^{\alpha+1}(1+x)^{\beta+1}dx.\\
\end{eqnarray*}
Thus, using Parseval's identity
\begin{eqnarray*}
&&\int_{0}^{\infty}\left|H_{\beta,K}^{N,2}(y)\right|^{2}y^{\alpha}e^{-y}\frac{dy}{\Gamma(\alpha+1)}\\
&& \hspace{0.6cm}\leq
\frac{\beta^{\alpha+2}}{\eta_{\alpha+1,\beta+1}2^{\alpha+\beta+3}\Gamma(\alpha+1)}\sum_{k=N+1}^{\infty}\sum_{n=1}^{k}\sum_{i=k-2n}^{k-2}
\left|\left\langle
\phi_{\beta},\frac{P^{(\alpha,\beta)}_{n}}{\left\|P^{(\alpha,\beta)}_{n}\right\|_{2}^{2}}\right\rangle\right|^{2}\left|\left\langle
\phi_{\beta},\frac{P^{(\alpha,\beta)}_{k-n}}{\left\|P^{(\alpha,\beta)}_{k-n}\right\|_{2}^{2}}\right\rangle\right|^{2}
\left|\frac{t}{\beta}\right|^{2}\\&&\hspace{1cm}\times
e^{-2t(\lambda_{k}^{1/2}+\lambda_{k-n}^{1/2})}\left|\frac{\lambda_{n}\lambda_{k-n}}{n(k-n)}P^{(\alpha+1,\beta+1)}_{n-1}(1)P^{(\alpha+1,\beta+1)}_{k-n-1}(1)\right|^{2}\\
&&\hspace{2cm}\times\left|\nu(i,k-n-1,n-1)\right|^{2}\left\|p^{(\alpha+1,\beta+1)}_{i}\right\|^{2}_{(\alpha+1,\beta+1)}.\\
\end{eqnarray*}
Analogously the previous case,
\begin{eqnarray*}
&&\int_{0}^{\infty}\left|H_{\beta,K}^{N,2}(y)\right|^{2}y^{\alpha}e^{-y}\frac{dy}{\Gamma(\alpha+1)}\\
&\leq
&\frac{\beta^{\alpha+2}\|\phi_{\beta}\|^{2}}{\eta_{\alpha+1,\beta+1}2^{\alpha+\beta+3}\Gamma(\alpha+1)}\sum_{k=N+1}^{\infty}\sum_{n=1}^{k}
\left|\left\langle
\phi_{\beta},\frac{P^{(\alpha,\beta)}_{n}}{\left\|P^{(\alpha,\beta)}_{n}\right\|_{2}}\right\rangle\right|^{2}
\frac{t^{2}2(n-1)}{\beta^{2}\left\|P^{(\alpha,\beta)}_{n}\right\|_{2}^{2}\left\|P^{(\alpha,\beta)}_{k-n}\right\|_{2}^{2}}\\&&\hspace{1cm}\times
e^{-2t(\lambda_{k}^{1/2}+\lambda_{k-n}^{1/2})}\left|\frac{\lambda_{n}\lambda_{k-n}}{n(k-n)}\right|^{2}.
\left(%
\begin{array}{c}
  n-1+q \\
  \\
  n-1\\
\end{array}%
\right)^{2}\left\|P^{(\alpha+1,\beta+1)}_{k-n-1}\right\|_{(\alpha+1,\beta+1)}^{2}.\\
\end{eqnarray*}
Also, given that
\begin{equation*}
\left\|P^{(\alpha+1,\beta+1)}_{k-n-1}\right\|_{(\alpha+1,\beta+1)}^{2}=\frac{(\alpha+\beta+3)(\alpha+\beta+2)(k-n)}{(\alpha+1)(\beta+1)(k-n+\alpha+\beta+1)}
\left\|P^{(\alpha,\beta)}_{k-n}\right\|_{(\alpha,\beta)}^{2},
\end{equation*}
we get,
\begin{eqnarray*}
&&\int_{0}^{\infty}\left|H_{\beta,K}^{N,2}(y)\right|^{2}y^{\alpha}e^{-y}\frac{dy}{\Gamma(\alpha+1)}\\
&\leq
&\frac{\beta^{\alpha+2}\|\phi_{\beta}\|^{2}}{\eta_{\alpha+1,\beta+1}2^{\alpha+\beta+3}\Gamma(\alpha+1)}\sum_{k=N+1}^{\infty}\sum_{n=1}^{k}
\left|\left\langle
\phi_{\beta},\frac{P^{(\alpha,\beta)}_{n}}{\left\|P^{(\alpha,\beta)}_{n}\right\|_{2}}\right\rangle\right|^{2}
\frac{t^{2}2e^{-2t(\lambda_{k}^{1/2}+\lambda_{k-n}^{1/2})}\lambda_{n}^{2}\lambda_{k-n}
}{\beta^{2}n}\hspace{2cm}\\&&\hspace{1cm}\times
\frac{\left(%
\begin{array}{c}
  n-1+q \\
  \\
  n-1\\
\end{array}%
\right)^{2}(\alpha+\beta+3)(\alpha+\beta+2)}{(\alpha+1)(\beta+1)\left\|P^{(\alpha,\beta)}_{n}\right\|_{2}^{2}}\\
&\leq
&\frac{e^{-t}\beta^{\alpha+1}\|\phi_{\beta}\|^{2}}{\eta_{\alpha,\beta}2^{\alpha+\beta+1}\Gamma(\alpha+1)}\sum_{k=N+1}^{\infty}\sum_{n=1}^{k}
\left|\left\langle
\phi_{\beta},\frac{P^{(\alpha,\beta)}_{n}}{\left\|P^{(\alpha,\beta)}_{n}\right\|_{2}}\right\rangle\right|^{2}
\frac{2e^{-t\lambda_{k-n}^{1/2}}(n+\alpha+\beta+1)}{t^{2}\beta}
\frac{\left(%
\begin{array}{c}
  n-1+q \\
  \\
  n-1\\
\end{array}%
\right)^{2}}{\left\|P^{(\alpha,\beta)}_{n}\right\|_{2}^{2}}.\\
\end{eqnarray*}
Then we have for $\beta>0$ big enough
\begin{equation*}
\frac{\left(%
\begin{array}{c}
  n-1+q \\
  \\
  n-1\\
\end{array}%
\right)^{2}(n+\alpha+\beta+1)}{\beta\left\|P^{(\alpha,\beta)}_{n}\right\|_{2}^{2}}\leq
\frac{Cn\beta^{2(n-1)}}{(\alpha+1)_{n}}\leq
\frac{Cn\beta^{2n}}{(\alpha+1)_{n}}.
\end{equation*}
Hence, for $\beta>0$ big enough
\begin{eqnarray*}
&&\int_{0}^{\infty}\left|H_{\beta,K}^{N,2}(y)\right|^{2}y^{\alpha}e^{-y}\frac{dy}{\Gamma(\alpha+1)}\\
&\leq &
\frac{Ce^{-t}\beta^{\alpha+1}\|\phi_{\beta}\|^{2}}{\eta_{\alpha,\beta}2^{\alpha+\beta+1}\Gamma(\alpha+1)}\sum_{k=N+1}^{\infty}\sum_{n=1}^{k}
\left|\left\langle
\phi_{\beta},\frac{P^{(\alpha,\beta)}_{n}}{\left\|P^{(\alpha,\beta)}_{n}\right\|_{2}}\right\rangle\right|^{2}\frac{2n\beta^{2n}e^{-t\lambda_{k-n}^{1/2}}}{t^{2}(\alpha+1)_{n}}
\end{eqnarray*}
then, analogous to (\ref{int-H1N}) we have
\begin{equation*}
\int_{0}^{\infty}\left|H_{\beta,K}^{N,2}(y)\right|^{2}y^{\alpha}e^{-y}\frac{dy}{\Gamma(\alpha+1)}\leq
\frac{Ce^{-t}\left\|\phi\right\|_{2,\alpha}^{4}}{e^{N+1}},
\end{equation*}
therefore,
\begin{eqnarray*}
\left(\int_{0}^{\infty}\left|\int_{0}^{\infty}H_{\beta,K}^{N,2}(y)dt\right|^{2}y^{\alpha}e^{-y}\frac{dy}{\Gamma(\alpha+1)}\right)^{1/2}&\leq
&\frac{(C)^{1/2}\|\phi\|_{2,\alpha}^{2}}{e^{\frac{(N+1)}{2}}}.\\
\end{eqnarray*}
Thus, $\left\{\int_{0}^{\infty}H_{\beta,K}^{N,2}dt\right\}$ is a
bounded sequence on
$\;L^{2}\left((0,\infty),\mu_{\alpha}\right),\;$ so by
Bourbaki-Alaoglu's theorem, there exists a sequence
$\;(\lambda_{j})_{j\in\mathbb{N}}\;$ $\lim_{j\to
\infty}\lambda_{j}=\infty$ such that, for all $N\in \mathbb{N}$,
  $\left\{\int_{0}^{\infty}H_{\beta_{j},K}^{N,2}dt\right\}_{j\in\mathbb{N}}$ converges weakly in $\;L^{2}\left((0,\infty),\mu_{\alpha}\right)\;$ to a function
  $H_{K}^{N,2}\in L^{2}\left((0,\infty),\mu_{\alpha}\right)\;$ Moreover,
  \begin{equation}\label{norLg}
\int_{0}^{\infty}\left|H_{K}^{N,2}(y)\right|^{2}y^{\alpha}e^{-y}\frac{dy}{\Gamma(\alpha+1)}\leq \frac{C}{e^{\frac{(N+1)}{2}}}\\
\end{equation}
Then, there exists a non decreasing sequence
$(N_{j})_{j\in\mathbb{N}}$ such that,
\begin{equation}\label{limLg}
H_{K}^{N_{j},2}(y)\longrightarrow 0,\hspace{0.2cm}
a.e.\hspace{0.2cm} y\in\mathbb(0,\infty).
\end{equation}
Therefore, similarly to (\ref{F-lamdk})
\begin{eqnarray*}
\left(F_{\beta
,K}(y)\right)^{2}&=&g_{1,\beta}\phi(y)+g_{2,\beta}\phi(y)\\
&=&F^{N}_{\beta,K}(y)\;+\;\int_{0}^{\infty}H_{\beta,K}^{N,1}(y)dt(1-\frac{y}{\beta})^{\beta/2}e^{y/2}\;+\;\int_{0}^{\infty} H_{\beta,K}^{N,2}(y)dty^{1/2}(1-\frac{y}{\beta})^{(\beta+1)/2}e^{y/2},\\
\end{eqnarray*}
where
\begin{eqnarray*}
F^{N}_{\beta,K}(y)&=&\sum_{k=1}^{N}\sum_{n=1}^{k}\left\langle
\phi_{\beta},\frac{P^{(\alpha,\beta)}_{n}}{\left\|P^{(\alpha,\beta)}_{n}\right\|_{2}^{2}}\right\rangle\left\langle
\phi_{\beta},\frac{P^{(\alpha,\beta)}_{k-n}}{\left\|P^{(\alpha,\beta)}_{k-n}\right\|_{2}^{2}}\right\rangle
\frac{\lambda_{n}^{1/2}\lambda_{k-n}^{1/2}}{(\lambda_{k}^{1/2}+\lambda_{k-n}^{1/2})^{2}}\\
&&\hspace{1cm}\times\;P^{(\alpha,\beta)}_{n}(1-\frac{2}{\beta}y)P^{(\alpha,\beta)}_{k-n}(1-\frac{2}{\beta}y)(1-\frac{y}{\beta})^{\beta}e^{y}\\
&&\hspace{0.2cm}+\;\sum_{k=1}^{N}\sum_{n=1}^{k}\left\langle
\phi_{\beta},\frac{P^{(\alpha,\beta)}_{n}}{\left\|P^{(\alpha,\beta)}_{n}\right\|_{2}^{2}}\right\rangle\left\langle
\phi_{\beta},\frac{P^{(\alpha,\beta)}_{k-n}}{\left\|P^{(\alpha,\beta)}_{k-n}\right\|_{2}^{2}}\right\rangle\frac{\lambda_{n}\lambda_{k-n}}{(\lambda_{k}^{1/2}+\lambda_{k-n}^{1/2})^{2}}
\\&&\hspace{1cm}\times
\;\frac{y}{\beta n(k-n)}P^{(\alpha+1,\beta+1)}_{n-1}(1-\frac{2}{\beta}y)P^{(\alpha+1,\beta+1)}_{k-n-1}(1-\frac{2}{\beta}y)(1-\frac{y}{\beta})^{\beta+1}e^{y}.\\
\end{eqnarray*}
Defining, for each $m\in\mathbb{N}$,
$F^{N_{m}}_{K}=F^{2}-H^{N_{m},1}_{K}-H^{N_{m},2}_{K},$ then since
$ F_{\beta,K}\to  F_{K},$ as $j\to\infty$ in the weak topology of
${L^{2}\left((0,\infty),\mu_{\alpha}\right)}$ and that
$$F_{\beta ,K}^{N}=\left(F_{\beta
,K}\right)^{2}-\int_{0}^{\infty}H_{\beta,K}^{N,1}(y)dt(1-\frac{y}{\beta})^{\beta/2}e^{y/2}\;-\;\int_{0}^{\infty}
H_{\beta,K}^{N,2}(y)dty^{1/2}(1-\frac{y}{\beta})^{(\beta+1)/2}e^{y/2},\;$$
we have, for all $m\in\mathbb{N},$
\begin{equation*}
F^{N_{m}}_{\beta_{j},K}\longrightarrow
F^{N_{m}}_{K},\hspace{1cm}\mbox{weakly in}
\;L^{2}\left((0,\infty),\mu_{\alpha}\right),\hspace{0.2cm}\mbox{as}\hspace{0.2cm}j\to
\infty,
\end{equation*}
and also,
\begin{equation}\label{F-Lague}
F^{N_{m}}_{K}(y)\longrightarrow
F^{2}(y),\hspace{0.3cm}\mbox{as}\hspace{0.3cm}m\to
\infty\hspace{0.3cm}a.e.\hspace{0.3cm}y\in(0,\infty).
\end{equation}
Now, given that
\begin{equation*}
\lim_{\beta\to\infty}P^{(\alpha+1,\beta+1)}_{n-1}(1-\frac{2}{\beta}y)=\lim_{\beta\to\infty}P^{(\alpha+1,\beta+1)}_{n-1}(1-\frac{2y}{\beta+1}\frac{\beta+1}{\beta})=L_{n-1}^{\alpha+1}(y)
\end{equation*}
and
\begin{eqnarray*}
\lim_{\beta\to\infty}\frac{1}{\left\|P^{(\alpha,\beta)}_{n}\right\|_{(\alpha,\beta)}^{2}}&=&\lim_{\beta\to\infty}(\frac{2n}{\beta}+\frac{\alpha}{\beta}+1+\frac{1}{\beta})n!
\frac{\Gamma(\alpha+1)}{\Gamma(n+\alpha+1)}\frac{\beta^{\alpha+1}\Gamma(\beta+1)}{\Gamma(\alpha+\beta+2)}\frac{\Gamma(\beta+n+\alpha+1)}{\beta^{\alpha}\Gamma(n+\beta+1)}\\
&=&\frac{n!\Gamma(\alpha+1)}{\Gamma(n+\alpha+1)}=\frac{1}{\left\|L_{n}^{\alpha}\right\|_{\alpha}^{2}};
\end{eqnarray*}
we get,
\begin{eqnarray*}
&&\lim_{\beta\to\infty}F^{N}_{\beta,K}(y)=\lim_{\beta\to\infty}\sum_{k=1}^{N}\sum_{n=1}^{k}\left\langle
\phi_{\beta},\frac{P^{(\alpha,\beta)}_{n}}{\left\|P^{(\alpha,\beta)}_{n}\right\|_{2}^{2}}\right\rangle\left\langle
\phi_{\beta},\frac{P^{(\alpha,\beta)}_{k-n}}{\left\|P^{(\alpha,\beta)}_{k-n}\right\|_{2}^{2}}\right\rangle
\frac{\lambda_{n}^{1/2}\lambda_{k-n}^{1/2}}{(\lambda_{k}^{1/2}+\lambda_{k-n}^{1/2})^{2}}\\
&&\hspace{1cm}\times\;P^{(\alpha,\beta)}_{n}(1-\frac{2}{\beta}y)P^{(\alpha,\beta)}_{k-n}(1-\frac{2}{\beta}y)(1-\frac{y}{\beta})^{\beta}e^{y}\\
&&\hspace{0.2cm}+\;\lim_{\beta\to\infty}\sum_{k=1}^{N}\sum_{n=1}^{k}\left\langle
\phi_{\beta},\frac{P^{(\alpha,\beta)}_{n}}{\left\|P^{(\alpha,\beta)}_{n}\right\|_{2}^{2}}\right\rangle\left\langle
\phi_{\beta},\frac{P^{(\alpha,\beta)}_{k-n}}{\left\|P^{(\alpha,\beta)}_{k-n}\right\|_{2}^{2}}\right\rangle\frac{\lambda_{n}\lambda_{k-n}}{(\lambda_{k}^{1/2}+\lambda_{k-n}^{1/2})^{2}}
\\&&\hspace{1cm}\times
\;\frac{y}{\beta n(k-n)}P^{(\alpha+1,\beta+1)}_{n-1}(1-\frac{2}{\beta}y)P^{(\alpha+1,\beta+1)}_{k-n-1}(1-\frac{2}{\beta}y)(1-\frac{y}{\beta})^{\beta+1}e^{y}\\
&=&\sum_{k=1}^{N}\sum_{n=1}^{k}\left\langle
\phi,\frac{L_{n}^{\alpha}}{\left\|L_{n}^{\alpha}\right\|_{\alpha}^{2}}\right\rangle\left\langle
\phi,\frac{L_{k-n}^{\alpha}}{\left\|L_{k-n}^{\alpha}\right\|_{\alpha}^{2}}\right\rangle
\frac{n^{1/2}(k-n)^{1/2}}{(n^{1/2}+(k-n)^{1/2})^{2}}L_{n}^{\alpha}(y)L_{k-n}^{\alpha}(y)\\
&&\hspace{0.2cm}+\;\sum_{k=1}^{N}\sum_{n=1}^{k}\left\langle
\phi,\frac{L_{n}^{\alpha}}{\left\|L_{n}^{\alpha}\right\|_{\alpha}^{2}}\right\rangle\left\langle
\phi,\frac{L_{k-n}^{\alpha}}{\left\|L_{k-n}^{\alpha}\right\|_{\alpha}^{2}}\right\rangle\frac{y}{(n^{1/2}+(k-n)^{1/2})^{2}}L^{\alpha+1}_{n-1}(y)L^{\alpha+1}_{k-n-1}(y)\\
&=&\int_{0}^{\infty}t\left(\sum_{k=1}^{N}\left\langle
\phi,\frac{L_{n}^{\alpha}}{\left\|L_{n}^{\alpha}\right\|_{\alpha}^{2}}\right\rangle
k^{1/2}e^{-tk^{1/2}}L_{n}^{\alpha}(y)\right)^{2}dt\\
&&\hspace{0.2cm}+\;\int_{0}^{\infty}t\left(\sum_{k=1}^{N}\left\langle
\phi,\frac{L_{n}^{\alpha}}{\left\|L_{n}^{\alpha}\right\|_{\alpha}^{2}}\right\rangle \sqrt{y} e^{-tk^{1/2}} L^{\alpha+1}_{n-1}(y)\right)^{2}dt.\\
\end{eqnarray*}
Thus,
\begin{eqnarray*}
F^{N_{m}}_{K}(y)&=&\int_{0}^{\infty}t\left(\sum_{k=1}^{N_{m}}\left\langle
\phi,\frac{L_{n}^{\alpha}}{\left\|L_{n}^{\alpha}\right\|_{\alpha}^{2}}\right\rangle
k^{1/2}e^{-tk^{1/2}}L_{n}^{\alpha}(y)\right)^{2}dt\\
&&\hspace{0.2cm}+\;\int_{0}^{\infty}t\left(\sum_{k=1}^{N_{m}}\left\langle
\phi,\frac{L_{n}^{\alpha}}{\left\|L_{n}^{\alpha}\right\|_{\alpha}^{2}}\right\rangle \sqrt{y} e^{-tk^{1/2}} L^{\alpha+1}_{n-1}(y)\right)^{2}dt.\\
\end{eqnarray*}
We want to prove that
\begin{equation}\label{F-Lag}
F^{N_{m}}_{K}(y)\to
\left(g^{\alpha}\phi(y)\right)^{2},\hspace{0.3cm}\mbox{as}\hspace{0.3cm}m\to
\infty\hspace{0.3cm}a.e.\hspace{0.3cm}y\in\mathbb{R}
\end{equation}
Indeed, initially we have
\begin{eqnarray*}
\left|\sum_{k=1}^{N_{m}}\left\langle\phi
,\frac{L_{k}^{\alpha}}{\left\|L_{k}^{\alpha}\right\|^{2}_{\alpha}}\right\rangle
k^{1/2}e^{-tk^{1/2}}L_{k}^{\alpha}(x)\right|&\leq &
\|\phi\|\sum_{k=1}^{\infty}k^{1/2}e^{-tk^{1/2}}\frac{\left|L_{k}^{\alpha}(x)\right|}{\left\|L_{k}^{\alpha}\right\|_{\alpha}}.
\end{eqnarray*}
From (8.22.1) of  \cite{sz}, for  $\alpha$ y $x>0$, we have
$$L_{k}^{\alpha}(x)=\pi^{-1/2}e^{x/2} x^{-\alpha/2-1/4}cos\{2(kx)^{1/2}-\alpha\pi/2-\pi/4\}+O(k^{\alpha/2-3/4}).$$
Thus, there exist $M>0$ such that
\begin{equation*}
\left|L_{k}^{\alpha}(x)-\pi^{-1/2}e^{x/2}x^{-\alpha/2-1/4}k^{\alpha/2-1/4}cos\{2(kx)^{1/2}-\alpha\pi/2-\pi/4\}\right|\leq
Mk^{\alpha/2-3/4},
\end{equation*}
hence,
\begin{equation*}
\left|L_{k}^{\alpha}(x)\right|\leq
k^{\alpha/2-1/4}\left(M+\pi^{-1/2}e^{x/2}x^{-\alpha/2-1/4}\right).
\end{equation*}
On the other hand, sin for  $k$ big enough
$
\left\|L_{k}^{\alpha}\right\|^{2}_{\alpha}\approx
\frac{k^{\alpha}}{\Gamma(\alpha+1)},
$
then, there exist $N_{1}>0$ such that
\begin{eqnarray*}
\frac{\left|L_{k}^{\alpha}(x)\right|}{\left\|L_{k}^{\alpha}\right\|_{\alpha}}&\leq
&
\frac{\left(M+\pi^{-1/2}e^{x/2}x^{-\alpha/2-1/4}\right)\left[\Gamma(k+1)\right]^{1/2}}{k^{1/4}},\quad \mbox{for any} \quad k\geq N_{1}.
\end{eqnarray*}
Thus,
\begin{eqnarray*}
\sum_{k=1}^{\infty}k^{1/2}e^{-tk^{1/2}}\frac{\left|L_{k}^{\alpha}(x)\right|}{\left\|L_{k}^{\alpha}\right\|_{\alpha}}&\leq
&
\left(M+\pi^{-1/2}e^{x/2}x^{-\alpha/2-1/4}\right)\left[\Gamma(\alpha+1)\right]^{1/2}\sum_{k=1}^{N_{1}}k^{1/4}e^{-tk^{1/2}}\frac{k^{\alpha/2}}{\left\|L_{k}^{\alpha}\right\|_{\alpha}}
\\&&
\;\;\;+\left(M+\pi^{-1/2}e^{x/2}x^{-\alpha/2-1/4}\right)\left[\Gamma(\alpha+1)\right]^{1/2}\sum_{k=1}^{\infty}k^{1/2}e^{-tk^{1/2}}.
\end{eqnarray*}
Then, for $t\geq 1$ we have
\begin{eqnarray}\label{cotaLg1}
\nonumber&&\left|t\left(\sum_{k=1}^{N_{m}}\left\langle\phi
,\frac{L_{k}^{\alpha}}{\left\|L_{k}^{\alpha}\right\|^{2}_{\alpha}}\right\rangle
k^{1/2}e^{-tk^{1/2}}L_{k}^{\alpha}(x)\right)^{2}\right|\\
\nonumber &\leq
&\left(M+\pi^{-1/2}e^{x/2}x^{-\alpha/2-1/4}\right)^{2}\Gamma(\alpha+1)\left(\sum_{k=1}^{N_{1}}k^{1/4}e^{-tk^{1/2}}\frac{k^{\alpha/2}}{\left\|L_{k}^{\alpha}\right\|_{\alpha}}\right)^{2}t
\\\nonumber&&\;\;\;+\left(M+\pi^{-1/2}e^{x/2}x^{-\alpha/2-1/4}\right)^{2}\Gamma(\alpha+1)\left(\sum_{k=1}^{N_{1}}k^{1/4}e^{-k^{1/2}}\frac{k^{\alpha/2}}{\left\|L_{k}^{\alpha}\right\|_{\alpha}}\right)t\left(\sum_{k=1}^{\infty}k^{1/2}e^{-tk^{1/2}}\right)\\\nonumber &&
\;\;\;+\left(M+\pi^{-1/2}e^{x/2}x^{-\alpha/2-1/4}\right)^{2}\Gamma(\alpha+1)t\left(\sum_{k=1}^{\infty}k^{1/2}e^{-tk^{1/2}}\right)^{2}.
\end{eqnarray}
Now, analogously to the Hermite case, we have
\begin{eqnarray*}
&&\int_{1}^{\infty}t\left\{\left(M+\pi^{-1/2}e^{x/2}x^{-\alpha/2-1/4}\right)^{2}\Gamma(\alpha+1)\left(\sum_{k=1}^{N_{1}}k^{1/4}e^{-tk^{1/2}}\frac{k^{\alpha/2}}{\left\|L_{k}^{\alpha}\right\|_{\alpha}}\right)^{2}
\right.\\&&\;\;\;+\left(M+\pi^{-1/2}e^{x/2}x^{-\alpha/2-1/4}\right)^{2}\Gamma(\alpha+1)\left(\sum_{k=1}^{N_{1}}k^{1/4}e^{-k^{1/2}}\frac{k^{\alpha/2}}{\left\|L_{k}^{\alpha}\right\|_{\alpha}}\right)\left(\sum_{k=1}^{\infty}k^{1/2}e^{-tk^{1/2}}\right)\\&&
\;\;\;\left.+\left(M+\pi^{-1/2}e^{x/2}x^{-\alpha/2-1/4}\right)^{2}\Gamma(\alpha+1)\left(\sum_{k=1}^{\infty}k^{1/2}e^{-tk^{1/2}}\right)^{2}\right\}dt\\
&\leq
&\left(M+\pi^{-1/2}e^{x/2}x^{-\alpha/2-1/4}\right)^{2}\Gamma(\alpha+1)\left\{\left(\sum_{k=1}^{N_{1}}\frac{e^{-(k^{1/2}-1/2)}k^{\alpha/2+1/4}}{\left\|L_{k}^{\alpha}\right\|_{\alpha}}\right)^{2}
\right.\\&&+\left(\sum_{k=1}^{N_{1}}\frac{e^{-k^{1/2}}k^{\alpha/2+1/4}}{\left\|L_{k}^{\alpha}\right\|_{\alpha}}\right)\left(\sum_{k=1}^{\infty}k^{1/2}e^{-(k^{1/2}-1)}\right)\;+\;\left.\left(\sum_{k=1}^{\infty}k^{1/2}e^{-(k^{1/2}-1/2)}\right)^{2}\right\}\int_{1}^{\infty}te^{-t}dt\;<\infty .
\end{eqnarray*}
For $\; 0<t<1,\;$ given that
$|\frac{\partial}{\partial t}P_{t}^{\alpha}(\phi(y))|\leq C(1 +
|y|)e^{-t}\;$ we get for $y\in(0,K)$
\begin{equation*}
t\left|\sum_{k=1}^{N_{m}}\left\langle\phi
,\frac{L_{k}^{\alpha}}{\left\|L_{k}^{\alpha}\right\|^{2}_{\alpha}}\right\rangle
k^{1/2}e^{-tk^{1/2}}L_{k}^{\alpha}(x)\right|^{2}<t\left(1+|\frac{\partial}{\partial
t}P_{t}^{\alpha}(\phi(y))|\right)^{2}<t\left(1+C(1 +
|y|)e^{-t}\right)^{2},
\end{equation*}
 where,
$$\int_{0}^{1}t(1+C(1 + |y|)e^{-t})^{2}dt<\infty,$$
 then, by Lebesgue's dominated convergence theorem we have
\begin{eqnarray*}
&&\lim_{m\to\infty
}\int_{0}^{\infty}t\left(\sum_{k=1}^{N_{m}}\left\langle\phi
,\frac{L_{k}^{\alpha}}{\left\|L_{k}^{\alpha}\right\|^{2}_{\alpha}}\right\rangle
k^{1/2}e^{-tk^{1/2}}L_{k}^{\alpha}(x)\right)^{2}dt\\
&=&\int_{0}^{\infty}t\left(\sum_{k=1}^{\infty}\left\langle\phi
,\frac{L_{k}^{\alpha}}{\left\|L_{k}^{\alpha}\right\|^{2}_{\alpha}}\right\rangle
k^{1/2}e^{-tk^{1/2}}L_{k}^{\alpha}(x)\right)^{2}dt=\int_{0}^{\infty}t\left(\frac{\partial}{\partial
t}P_{t}^{\alpha}(\phi(y))\right)^{2}dt.
\end{eqnarray*}
Similarly as before,
\begin{equation*}
\lim_{m\to\infty
}\int_{0}^{\infty}t\left(\sum_{k=1}^{N_{m}}\left\langle
\phi,\frac{L_{n}^{\alpha}}{\left\|L_{n}^{\alpha}\right\|_{\alpha}^{2}}\right\rangle
\sqrt{y}e^{-tk^{1/2}}L^{\alpha+1}_{n-1}(y)\right)^{2}dt=\int_{0}^{\infty}t\left(\sqrt{y}\frac{\partial}{\partial
y}P_{t}^{\alpha}(\phi(y))\right)^{2}dt.
\end{equation*}
Thus,
\begin{equation*}
F^{N_{m}}_{K}(y)\to
\left(g^{\alpha}\phi(y)\right)^{2},\hspace{0.3cm}\mbox{as}\hspace{0.3cm}m\to
\infty\hspace{0.3cm}a.e.\hspace{0.3cm}y\in\mathbb{R},
\end{equation*}
hence, from (\ref{F-Lague}) and (\ref{F-Lag}) we have $
F(y)=g^{\alpha}(\phi(y)).\;\blacksquare
$

\end{enumerate}


\end{document}